%% file: main_firstdraft.tex
\documentclass[preprint,12pt]{elsarticle}

\usepackage{amssymb,amsmath,amsfonts,bm}
\usepackage[margin=1in]{geometry}
\usepackage{cuted}
\usepackage[utf8]{inputenc}
\usepackage[T1]{fontenc}
\usepackage[english]{babel}
\usepackage{fancyhdr}
\usepackage{anyfontsize}
\usepackage[active]{srcltx}
\usepackage{soul}
\usepackage{xcolor}
\usepackage[colorlinks=true, allcolors=blue]{hyperref}
\usepackage{url}
\usepackage{color}
\usepackage{array}
\usepackage{textcomp}
\usepackage[linesnumbered,ruled,vlined]{algorithm2e}
\usepackage{pgfplots}
\pgfplotsset{compat=1.17}
\input{newcommands}

\usepackage{subfiles}
\usepackage{subcaption}
\usepackage{smartdiagram}
\usesmartdiagramlibrary{additions}
\usepackage{lineno}

\usepackage{upgreek}
\usepackage{textgreek}
\usepackage{nicefrac}

\usepackage{cancel}
\usepackage[capitalise]{cleveref}
\Crefname{equation}{Eq.}{Eqs.}
\Crefname{appendix}{Appendix}{Appendices}
\usepackage{comment}
\usepackage{svg}
\usepackage{float}
\allowdisplaybreaks
\usepackage{comment}

\modulolinenumbers[5]

\journal{}

\begin{document}

\begin{frontmatter}



\title{Reduced-order modelling of parameter-dependent systems with invariant manifolds: application to Hopf bifurcations in follower force problems
}

\author[1]{André de F. Stabile\corref{cor1}}
\ead{andre.de-figueiredo-stabile@ensta-paris.fr}
\author[2]{Alessandra Vizzaccaro}
\author[3]{Loïc Salles}
\author[4]{Alessio Colombo}
\author[4]{Attilio Frangi}
\author[1]{Cyril Touzé}

\affiliation[1]{organization={Institute of Mechanical Sciences and Industrial Applications (IMSIA), ENSTA Paris - CNRS - EDF - CEA, Institut Polytechnique de Paris},
            city={Palaiseau},
            country={France}}

\affiliation[2]{organization={College of Engineering, Mathematics and Physical Sciences, University of Exeter},
            city={Exeter},
            country={UK}}

\affiliation[3]{organization={Department of Aerospace \& Mechanical Engineering, University of Liège},
            city={Liège},
            country={Belgium}}

\affiliation[4]{organization={Department of Civil and Environmental Engineering, Politecnico di Milano},
            city={Milan},
            country={Italy}}

\cortext[cor1]{Corresponding author}

\begin{abstract}
The direct parametrisation method for invariant manifolds is adjusted to consider a varying parameter. More specifically, the case of systems experiencing a Hopf bifurcation in the parameter range of interest are investigated, and the ability to predict the amplitudes of the limit cycle oscillations after the bifurcation is demonstrated. The cases of the Ziegler pendulum and Beck's column, both of which have a follower force, are considered for applications. By comparison with the eigenvalue trajectories in the conservative case, it is advocated that using two master modes to derive the ROM, instead of only considering the unstable one, should give more accurate results. Also, in the specific case where an exceptional bifurcation point is met, a numerical strategy enforcing the presence of Jordan blocks in the Jacobian matrix during the procedure, is devised. The ROMs are constructed for the Ziegler pendulum having two and three degrees of freedom, and then Beck's column is investigated, where a finite element procedure is used to space discretize the problem. The numerical results show the ability of the ROMs to correctly predict the amplitude of the limit cycles up to a certain range, and it is shown that computing the ROM after the Hopf bifurcation gives the most satisfactory results. This feature is analyzed in terms of phase space representations, and the two proposed adjustments are shown to improve the validity range of the ROMs.
\end{abstract}



\end{frontmatter}

\section{Introduction}

Nonlinear techniques for simulation-free model order reduction (MOR) based on the parametrisation method for invariant manifolds, show several successful applications in recent years, especially for vibrating systems~\cite{PONSIOEN2018,PONSIOEN2020,BreunungHaller18,artDNF2020,JAIN2021How,vizza21high,opreni22high,Mereles:rotor,martin2019,li2021periodic,opreniPiezo,Frangi:electromech,Pinho:shells}, but also for different domains like fluid dynamics with the Navier-Stokes equation~\cite{Buza:NS} or astronomy~\cite{LeBihan_2017}. In the field of nonlinear vibrations, all these works complement and improve the earlier developments defining and using Nonlinear Normal Modes (NNMs) as invariant manifolds for accurate derivations of reduced order models (ROMs)~\cite{ShawPierre91,Shaw94,touze03-NNM,TOUZE:JSV:2006,Haller2016,ReviewROMGEOMNL,TouzeCISM2} thanks to the following decisive advances. First, the uniqueness of the invariant manifold has been theoretically proven in~\cite{Haller2016} thanks to the parametrisation method, yielding spectral submanifolds (SSMs) as the unique smoothest subspace, which is thus the sought and already used NNM of earlier works. Second, arbitrary order expansions have been devised, allowing automated powerful results with a fine convergence as long as the validity limit of the local theory is fulfilled~\cite{PONSIOEN2018,vizza21high,LamarqueUP}. Finally, methods that are directly applicable to a Finite Element (FE) discretisation, of broad use for engineering structures, have been developed~\cite {artDNF2020,JAIN2021How,vizza2024superharm}. 

In the previously cited contributions, most of the developments were concerned with computing the invariant manifolds of a fixed point representing the structure at rest. In this realm, NNMs offer a clear and direct continuation of the idea of linear normal modes. While the majority of the developments for vibrating structures were concerned with geometric nonlinearity, recent applications extend the approach and consider coupled systems or different physics. For instance, friction is considered in~\cite{Mereles:rotor,Mereles:bif}, a weak piezo-electric coupling is addressed in~\cite{opreniPiezo}, while a strong electromechanical coupling with application to Micro-Electro-Mechanical Systems (MEMS), is derived in~\cite{Frangi:electromech}.

The treatment of parameter-dependent ROMs using either the center manifold approach or the normal form theory has already been considered in the past. For example, the idea of adding the parameter as an additional state variable with trivial dynamics is already addressed in~\cite{NayfehBala} in the context of center manifold reduction of a simple system experiencing a Hopf bifurcation, and the technique has been applied to friction-induced vibration in braking systems in~\cite{Sinou2003}. In the context of forced systems, the external forcing has been added as an extra oscillator for the computation of NNMs with the center manifold technique in~\cite{ShawShock}. Considering now normal form theory, parameter-dependent cases in nonlinear vibrations, including bifurcations leading either to divergence or flutter, were already investigated in~\cite{Hsu:NFa,Hsu:NFb}, where the particular case of the Ziegler pendulum was considered. Additionally, other treatments in the realm of parameter-dependent systems can also be found for example in~\cite{IoossForcing,Jezequel91}, including the case of external forcing.

Now focusing on the more general context of the parametrisation method for invariant manifolds, parameter-dependent problems have been less frequently addressed. For example, in the non-autonomous case, where the excitation frequency can be seen as a varying parameter, different treatments have been proposed. In~\cite{vizza2024superharm}, the ROM is computed for a single value and then used to compute bifurcation diagrams with numerical continuation for slight variations in the vicinity of the expansion point. On the other hand, the invariant manifolds are computed for each forcing frequency in~\cite{JIANG2005H,Thurnher:resparam}. For rotating systems, interpolations between different ROMs computed at selected rotating speeds have proven effective in~\cite{Martin:rotation}. On the other hand, MOR techniques using manifolds and embeddings that are especially concerned with large parameter variations have also been exploited in different fields, see {\em e.g.}~\cite{SilkeGlas:MorMan,oulghelou:grassmann,Barnett2022} and references therein.

The previously cited approaches to deal with parameter-dependent ROMs rarely consider the cases where the fixed point, in the vicinity of which the NNMs are computed, experiences a bifurcation in the considered parameter range. However, bifurcations of the fixed point are commonly encountered in diverse mechanical situations~\cite{Bolotin1964}, for instance in the case of buckling, where the position of the structure at rest becomes unstable for an increasing value of the load. Another case of interest is that of the flutter, where a Hopf bifurcation occurs with the birth of limit cycle oscillations. This case is typically encountered for a wing subjected to a uniform flow with increasing velocity~\cite{Dowell,Dimitriadis}, but also for pipes conveying fluid~\cite{GregoryPaidoussis,DOARE2002}, or structures actuated with follower forces as the Ziegler pendulum~\cite{Ziegler1952,Bentvelsen2018,Luongo:minZieg} or Beck's column~\cite{Beck1952,Migli:Beck} for instance.

In such a case, the parametrisation method for invariant manifolds needs to be revisited to account for this major change in the dynamics. For a system encountering a Hopf bifurcation, a first important step has been proposed in~\cite{MingwuLi2024} with application to a pipe conveying fluid. The parametrisation method has been computed by selecting the unstable mode as the master one, and adjusted by incorporating the bifurcation parameter as an additional variable with trivial dynamics, in a treatment resembling the so-called ``suspension trick'' largely used in dynamical systems, see {\em e.g.}~\cite{NayfehBala,ShawShock,vizza2024superharm}. In particular, the parametrisation method is computed for a given value of the flow velocity, and the ROM is then used to predict the amplitudes of the limit cycles after the Hopf bifurcation. Interestingly, it has been found that the best results are obtained when the parametrisation point is selected after the bifurcation, meaning that the unstable invariant manifold of the unstable fixed point provides a better approximation of the limit cycles than those obtained using the center manifold computed exactly at the Hopf bifurcation parameter value.

This contribution aims to elaborate further on the use of MOR techniques for vibrating systems experiencing a Hopf bifurcation, using the direct parametrisation method for invariant manifolds (DPIM). The bifurcation parameter is introduced as an added variable with trivial dynamics, and two distinctive features are addressed as compared to~\cite{MingwuLi2024}. First, it is shown how considering two master modes can improve the predictions. The two modes are selected from the inspection of the undamped or lightly damped problem since this framework uncovers the coalescence of frequencies, at the heart of the instability, and the presence of a 1:1 resonance between the two bifurcating modes. Even if this information is lost when damping is added, it is shown to improve the results of the ROM, underlining that the near 1:1 resonance affects the quality of the ROM. 
Second, in the case of a perfect frequency coalescence, degenerate eigenvalues occur and Jordan blocks appear in the linear part of the dynamics. In such a case, a numerical strategy that enforces the presence of Jordan blocks is tested to improve the results. 
Finally, the findings are numerically assessed in three cases: the Ziegler pendulum with two and three degrees of freedom (DOF), and a cantilever beam with a follower force (Beck's column). In the latter case, the problem is discretised by the finite element method to illustrate an application to a numerical problem featuring a large number of  DOF.

\section{Models and methods}

In this section, a general methodology for embedding parameter variations in the framework of the DPIM is first introduced, following~\cite{MingwuLi2024,vizza2024superharm}. Then, the mechanical models experiencing Hopf bifurcations are presented, with a special emphasis on their specificities in terms of eigenvalue trajectories as a function of the bifurcation parameter. The Ziegler pendulum is here selected as a prototypical case. Finally, the adjustments proposed to enhance the ROM's predictive capacity are detailed based on the understanding of the linear behaviour.

\subsection{Direct parametrisation method with added bifurcation parameter} \label{subsec:Reduction}

This section is concerned with the inclusion of a parameter in the formalism of the DPIM. The procedure strictly follows the idea presented in~\cite{MingwuLi2024,vizza2024superharm}, such that the presentation is kept minimal. The reader is referred to~\cite{Haro} for a general presentation of the parametrisation method, and~\cite{vizza2024superharm,TouzeCISM2,JAIN2021How} for its direct application to vibrating systems.

For the sake of generality, let us consider a generic parameter-dependent problem in the form of a differential-algebraic equation (DAE) for an unknown state space vector $\bfy \in \mathbb{C}^D$, and a scalar control parameter $\mu$. The ROM is constructed from the parametrisation method at a given value of the parameter $\mu_0$, and aims at predicting the behaviour for varying parameter values in the vicinity. In order to accurately embed the parameter variations in the formulation, the parameter is taken as an additional state variable having trivial dynamics, since being time-independent. Note that the main interest is in bifurcating systems, thus a bifurcation might occur for $\mu=\mu_0$. However, the method as presented in this section is general and can be applied with or without specifying that a bifurcation occurs at $\mu_0$, neither indicating the type of bifurcation. Since solutions are searched for in the vicinity of $\mu_0$, the system is written as:
\begin{subequations} \label{eq:DAE}
\begin{align}
    \bfB \dot{\bfy} &= \bfA \left( \bfy_0 + \bfy \right) + \bfQ_1(\bfy_0 + \bfy,\bfy_0 + \bfy) + \bfQ_2(\bfy_0 + \bfy,\mu_0 + \mu) + \bfQ_3(\mu_0 + \mu,\mu_0 + \mu), \label{eq:DAEa} \\
    \dot{\mu} &= 0. \label{eq:DAEb}
\end{align}
\end{subequations}
The following features have been considered. First, the parameter $\mu$ is added to the original problem by adding~\cref{eq:DAEb} which has trivial dynamics. This step is key to applying the DPIM to an extended problem where the bifurcation parameter will be properly dealt with. Second, the nonlinearities are supposed to be quadratic and described by the smooth and analytical bilinear functions $\bfQ_1 : \mathbb{C}^{D} \times \mathbb{C}^{D} \rightarrow \mathbb{C}^{D}$, $\bfQ_2 : \mathbb{C}^{D} \times \mathbb{R} \rightarrow \mathbb{C}^{D}$ and $\bfQ_3 : \mathbb{R}^2 \rightarrow \mathbb{C}^D$. The nonlinear terms have been split into three parts to take into account quadratic nonlinearities involving only the state ($\bfQ_1$), or only the parameter ($\bfQ_3$), or mixed ($\bfQ_2$). The assumption of only quadratic nonlinearity is not restrictive, as any smooth nonlinearity can be transformed into a quadratic one thanks to a quadratic recast~\cite{COCHELIN2009,KARKAR2013,Guillot2019}.

The same assumptions as in \cite{vizza2024superharm} relative to the matrices are considered: $\bfA$ and $\bfB$ are real-valued, with $\bfA$ of rank $D$ and $\bfB$ possibly rank-deficient, where $D$ is the dimension of the state-space. The fixed point $\bfy_0$ for $\mu=\mu_0$ is such that
\begin{equation} \label{eq:StaticEq}
    \bfA \bfy_0 + \bfQ_1(\bfy_0,\bfy_0) + \bfQ_2(\bfy_0,\mu_0) + \bfQ_3(\mu_0,\mu_0) = \mathbf{0},
\end{equation}
with $\mu$ imposing an increment to this $\mu_0$ value and $\bfy$ a (finite) perturbation around the equilibrium state. 

\cref{eq:DAEa} can be simplified by noticing that the following equations hold for the quadratic tensors:
\begin{subequations} \label{eq:Bilinear}
\begin{align}
    \bfQ_1(\bfy_0+\bfy,\bfy_0+\bfy) &= \bfQ_1(\bfy_0,\bfy_0) + \bfQ_1(\bfy_0,\bfI) \bfy + \bfQ_1(\bfI,\bfy_0) \bfy + \bfQ_1(\bfy,\bfy) \\
    \bfQ_2(\bfy_0+\bfy,\mu_0+\mu) &= \bfQ_2(\bfy_0,\mu_0) + \bfQ_2(\bfy_0,1) \mu + \bfQ_2(\bfI,\mu_0) \bfy + \bfQ_2(\bfy,\mu) \\
    \bfQ_3(\mu_0+\mu,\mu_0+\mu) &= \bfQ_3(\mu_0,\mu_0) + \bfQ_3(\mu_0,1) \mu + \bfQ_3(1,\mu_0) \mu + \bfQ_3(\mu,\mu),
\end{align}
\end{subequations}
with $\bfI \in \mathbb{R}^{D\times D}$ the identity matrix. In the above equations, the identity matrix is given as an input to the quadratic tensors, which constitutes a slight abuse of notation. In order to avoid any misinterpretation, the meaning of this operation is detailed in \cref{ap:NonlinearTensors} with indicial expressions. Inserting \cref{eq:Bilinear} into \cref{eq:DAEa}, canceling terms with the fixed point~\cref{eq:StaticEq}, and introducing the augmented state variables vector $\bfytil = \left[\bfy \quad \mu \right]^T$, the following system is found:
\begin{equation} \label{eq:DAEaugmented}
    \underbrace{
    \begin{bmatrix}
        \bfB & \mathbf{0} \\
        \mathbf{0} & 1
    \end{bmatrix}
    }_{\bfBtil}
    \underbrace{
    \dot{
    \begin{bmatrix}
        \bfy \\
        \mu
    \end{bmatrix}
    }
    }_{\dot{\bfytil}}
    =
    \underbrace{
    \begin{bmatrix}
        \bfA_t & \bfA_0 \\
        \mathbf{0} & 0
    \end{bmatrix}
    }_{\bfAtil_t}
    \underbrace{
    \begin{bmatrix}
        \bfy \\
        \mu
    \end{bmatrix}
    }_{\bfytil}
    +
    \underbrace{
    \begin{bmatrix}
        \bfQ_1(\bfy,\bfy) + \bfQ_2(\bfy,\mu) + \bfQ_3(\mu,\mu) \\
        0
    \end{bmatrix}
    }_{\bfQtil(\bfytil,\bfytil)},
\end{equation}
with $\bfA_t$ defined by
\begin{equation} \label{eq:At}
    \bfA_t = \bfA + \bfQ_1(\bfy_0,\bfI) \bfy + \bfQ_1(\bfI,\bfy_0) \bfy + \bfQ_2(\bfI,\mu_0),
\end{equation}
and 
\begin{equation}
    \bfA_0 = \bfQ_2(\bfy_0,1) + \bfQ_3(\mu_0,1) \mu + \bfQ_3(1,\mu_0) \mu.
\end{equation}

Interestingly, \cref{eq:DAEaugmented} has exactly the same shape as the starting equations used in~\cite{vizza2024superharm}. Consequently, the DPIM can be applied to~\cref{eq:DAEaugmented} with minor modifications. The only point needing further attention is the matrix $\bfAtil_t$, which is singular because of its last line, a point that was not addressed in~\cite{vizza2024superharm}. The assumption that the linear part of the right-hand side is non-singular is only used to solve the associated eigenvalue problem. Thus, it is possible to extend the method and consider a singular matrix $\bfAtil_t$ by treating differently the physical eigenvectors and the one associated with the bifurcation parameter. 

Let us assume that the master modes selected to perform the DPIM are the first $d$ ones. It is here recalled that the method only requires the computation of the $d\ll D$ master modes. The right eigenproblem associated to~\cref{eq:DAEaugmented} reads
\begin{equation} \label{eq:eigenvalues}
    \left( 
    \begin{bmatrix}
        \bfA_t & \bfA_0 \\
        \mathbf{0} & 0
    \end{bmatrix} 
    - \lambda_s 
    \begin{bmatrix}
        \bfB & \mathbf{0} \\
        \mathbf{0} & 1
    \end{bmatrix}
    \right) 
    \begin{bmatrix}
         \bfY_s \\
         \text{Y}_s^\mu
    \end{bmatrix} 
    = \mathbf{0}, \quad s=1,\ldots,d+1,
\end{equation}
where the right eigenvector has been partitioned into two parts, respectively $\bfY_s$ related to the usual state variables, and $\text{Y}_s^\mu$ to the bifurcation parameter. When $s \leq d$, the eigenvalue $\lambda_s$ is supposed to be non-vanishing. Thus, the last line of \cref{eq:eigenvalues} implies $\text{Y}_s^\mu = 0$, while the other lines reduce to the usual eigenvalue problem
\begin{equation} \label{eq:usualEigenvalues1storder}
    \left( \bfA_t - \lambda_s \bfB \right) \bfY_s = \mathbf{0} \quad s=1,\ldots,d.
\end{equation}
Instead, when the case $s=d+1$ is considered, $\lambda_{d+1} = 0$ since the control parameter does not evolve with time. Then, since the eigenvalues normalization is arbitrary, by choosing $\text{Y}_{d+1}^\mu = 1$, the following equation needs to be solved 
\begin{equation} \label{eq:pEigenvalue1storder}
    \bfA_t \bfY_{d+1} = -\bfA_0.
\end{equation}
The solution to this equation exists as long as the matrix $\bfA_t$ is non-singular, which is always supposed to be the case for the applications treated in this contribution. It should be noticed that, if desired, the normalization of $\bfY_{d+1}$ can be altered afterwards,  {\em e.g.} to respect mass orthonormality. Thanks to~\cref{eq:pEigenvalue1storder}, the eigenvector relative to the bifurcation parameter added as a state variable can be easily computed and the linear master eigenvectors are known. Note that this added eigenvector has a physical meaning, see 
\cref{ap:2ordersimpl} where this is further detailed for the specific case of second-order mechanical systems.

The DPIM can thus be applied directly to \cref{eq:DAEaugmented}. The method assumes first that a nonlinear mapping exists, relating the physical variables $\bfy$ to the so-called {\em normal} coordinate $\bfztil$ as
\begin{equation}\label{eq:ziNLmap00}
    \bfy = \bfW(\bfztil).
\end{equation}
The vector $\bfztil$ is of dimension $d+1$, and contains the $d \ll D$ usual normal coordinates $\bfz$, appended with the parameter $\mu$ as added variable:
\begin{equation}
    \bfztil = 
    \begin{bmatrix}
        \bfz \\
        \mu
    \end{bmatrix}.
\end{equation}
Interestingly, $\mu$ is discarded from the left-and side of~\cref{eq:ziNLmap00} following the idea given in~\cite{vizza2024superharm}, since being a non-physical coordinate with trivial dynamics, which must not be confused with the original problem. Finally, the reduced dynamics, governing the evolution onto the selected $(d+1)$-dimensional invariant manifold, is introduced as an unknown as
\begin{equation}\label{eq:zireducedyn00}
    \dot{\bfztil} = \bff(\bfztil).
\end{equation}
Again, the method is derived such that the last line of~\cref{eq:zireducedyn00} is left unmodified and always equal to~\cref{eq:DAEb}, {\em i.e.} $\dot{\mu}=0$.

\cref{eq:ziNLmap00,eq:zireducedyn00} can be inserted into~\cref{eq:DAEaugmented} to eliminate time, giving rise to the so-called invariance equation, which is then solved recursively, by considering polynomial expansions for~\cref{eq:ziNLmap00,eq:zireducedyn00}. More specifically, one introduces:
\begin{subequations}
    \begin{align}
        \bfW(\bfztil) &= \sum_{p=1}^o \sum_{k=1}^{m_p} \bfW^{(p,k)} \bfztil^{\alphavec(p,k)}, \\
        \bff(\bfztil) &= \sum_{p=1}^o \sum_{k=1}^{m_p} \bff^{(p,k)} \bfztil^{\alphavec(p,k)},
    \end{align}
\end{subequations}
with the unknown coefficients $\bfW^{(p,k)}$ and $\bff^{(p,k)}$. The multi-index notation is used, and $\alphavec (p,k)= \left\{ \alpha_1 \quad \alpha_2 \quad \ldots \quad \alpha_{d+1} \right\}$ refers to the $k$-th monomial of order $p$, $k\in[1,m_p]$, $m_p$ being the number of monomials of order $p$ in $d+1$ coordinates. The monomial associated to $\alphavec (p,k)$ simply reads $\bfztil^{\alphavec (p,k)} = z_1^{\alpha_1}z_2^{\alpha_2} \hdots \mu^{\alpha_{d+1}}$. Each $\alpha_j$ is such that $0\leq \alpha_j\leq p$, and collects the power associated to $z_j$, such that $\sum_{j=1}^{d+1} \alpha_j = p$. The maximal order of the polynomial expansion is $o$, and the numerical solutions will be referred to as order-$o$ solutions in the remainder.

Introducing the polynomial expansions in the invariance equations leads to the so-called homological equations of order $p$, which are solved recursively. Since~\cref{eq:DAEaugmented} is exactly in the same format as the systems treated in \cite{vizza2024superharm}, the methodology addressed therein, and specifically the algorithm and expressions for the different terms composing the homological equations, can be directly applied to the present case. For this reason, neither the assembly of the terms in the homological equations nor the solution process of the system of equations will be further detailed in this contribution, the interested reader is directly referred to~\cite{vizza2024superharm} for algorithmic details.

An important remark however stands, since the eigenvalue associated with the added variable representing the control parameter vanishes. As a consequence, numerous new trivial resonances will be fulfilled at the nonlinear level, giving rise to much more resonant monomials in a normal form style solution when compared to the case of lightly damped mechanical systems. This has important consequences on the nature of the ROMs produced, but does not impair the DPIM algorithm, which can easily handle such cases thanks to the resonance sets used to fill the resonance relationships, see~{\em e.g.}~\cite{vizza2024superharm,vizza21high} for more details. 

The mechanical systems under study and their eigenvalues' dependence upon the bifurcation parameter are now detailed to justify the specific choices to consider for a vibrating system experiencing a flutter-like instability.

\subsection{Hopf bifurcations and Ziegler pendulum} \label{subsec:Systems}

The proposed reduction technique aims to provide an accurate estimate of the amplitudes of the limit cycles for a mechanical system undergoing a Hopf bifurcation.  In order to assess and develop the method on a simple example featuring these main characteristics, the case of a Ziegler pendulum is selected~\cite{Ziegler1952,Luongo2015,Bentvelsen2018}. This section introduces the equations of motion of the Ziegler problem, and recalls some basic results related to the stability analysis, which are key to understanding the choices operated to derive the ROM.

\subsubsection{The Ziegler pendulum}\label{subsubsec:zieg00}

\begin{figure}[h!]
\centering
\begin{subfigure}{0.49\textwidth}
    \centering
    \includegraphics[scale=0.8]{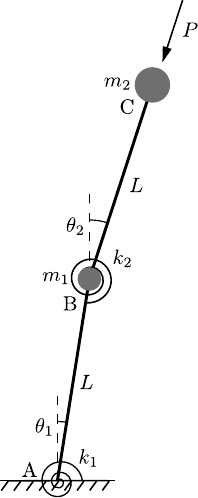}
    \caption{}
    \label{fig:2DOFZiegler}
\end{subfigure}
\begin{subfigure}{0.49\textwidth}
    \centering
    \includegraphics[scale=0.8]{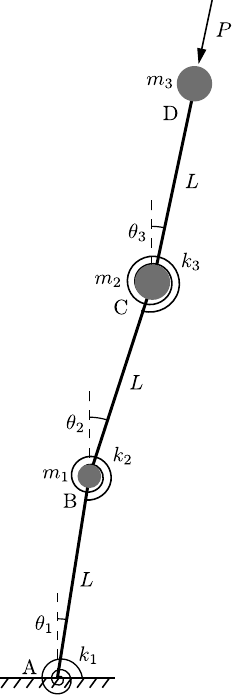}
    \caption{}
    \label{fig:3DOFZiegler}
\end{subfigure}
\caption{Two models of Ziegler pendulum. (a)  classical 2-DOF system.  (b) A 3-DOF version.}
\end{figure}

The Ziegler pendulum, initially introduced in~\cite{Ziegler1952} and shown in~\cref{fig:2DOFZiegler}, is a 2-DOF inverted pendulum excited by a follower force. Two rigid bars of length $L$ with concentrated masses $m_1$ and $m_2$ at their ends are connected by springs of stiffness $k_1$ and $k_2$. A follower force of magnitude $P$, always aligned with the second bar $BC$, is the external excitation, and $P$ is the bifurcation parameter. Angles $\theta_1$ and $\theta_2$ of the bars with the vertical are used as generalized coordinates for the model. A 3-DOF version of this system will also be considered in~\cref{subsec:3DOFZiegler}, and is shown in~\cref{fig:3DOFZiegler}. 

The equations of motion for the 2-DOF Ziegler pendulum are well-established in the literature, see e.g. \cite{Ziegler1952,Bentvelsen2018,Luongo2015}, and up to third-order, read
\begin{equation} \label{eq:EOMZiegler2DOF}
    \bfM \ddot{\thetavec} + \bfC \dot{\thetavec} + \left( \bfK + \bfKg \right) \thetavec = \mathbf{F}_{nl},
\end{equation}
with
\begin{equation}
\begin{gathered}
    \bfM = L^2
    \begin{bmatrix}
        m_1+m_2 & m_2 \\
        m_2 & m_2 
    \end{bmatrix},
    \quad 
    \bfK =
    \begin{bmatrix}
        k_1+k_2 & -k_2 \\
        -k_2 & \phantom{-}k_2 
    \end{bmatrix},
    \quad 
    \bfKg = PL 
    \begin{bmatrix}
        -1 & 1 \\
        \phantom{-}0 & 0 
    \end{bmatrix},
    \\
    \thetavec =
    \begin{bmatrix}
        \theta_1 \\
        \theta_2 
    \end{bmatrix},
    \quad
    \bfF_{nl} = 
    -\frac{PL}{6}
    \begin{bmatrix}
        (\theta_1-\theta_2)^3 \\
        0
    \end{bmatrix}.
\end{gathered}
\end{equation}
An arbitrary damping matrix $\bfC$ is introduced in~\cref{eq:EOMZiegler2DOF}. In many applications, a structural Rayleigh damping of the form $\bfC = 2 \left( \xi_k \bfK + \xi_m \bfM \right) $ is considered, with $\xi_k$ and $\xi_m$ the amplitudes of stiffness and mass-proportional terms.

\subsubsection{Linear stability analysis} \label{subsec:LinearStability}

By monitoring the evolution of the eigenvalues with the control parameter~$P$, a Hopf bifurcation is seen to occur once two complex conjugate eigenvalues cross simultaneously the imaginary axis. 
This instability is generally referred to as the flutter instability in aeroelastic problems. When damping is not considered, it is characterized by a frequency coalescence~\cite{Rocard49,fung2008introduction,Dimitriadis,Kirillov:stab}: two of the eigenfrequencies of the system merge for a given value of the bifurcation parameter. Without damping, this point is also exactly that of the Hopf bifurcation since a positive real part appears just after the coalescence.

Following the terminology introduced in~\cite{Heiss2000,Seyranian2005,Kirillov:stab}, the point where the eigenfrequencies coalesce is called an {\em exceptional point} (EP): both eigenvalues and eigenvectors merge, forming a Jordan block. In other words, the algebraic multiplicity of the merging eigenvalues is 2 while the geometric multiplicity is 1. The other possible case is that of a {\em diabolic point} (DP): algebraic and geometric multiplicities are equal to two, meaning that eigenvalues are merging but the eigenvectors remain different and stay linearly independent~\cite{Seyranian2005}. A complete study of the possible cases is investigated in~\cite{Seyranian2005}. Interestingly, the appearance of diabolic points is not possible
for codimension 1 bifurcations, such as the Hopf bifurcation, when the Jacobian matrix
describing the dynamics is real and asymmetric, as in the present case. This implies that, as
long there is a coalescence of eigenfrequencies, an exceptional point is necessarily at hand, and
no further check on the dimension of the subspace generated by the eigenvectors is needed: the Jordan block is present.

\begin{figure}[h!]
\centering
\begin{subfigure}{0.49\textwidth}
    \centering
    \includegraphics[width=\textwidth]{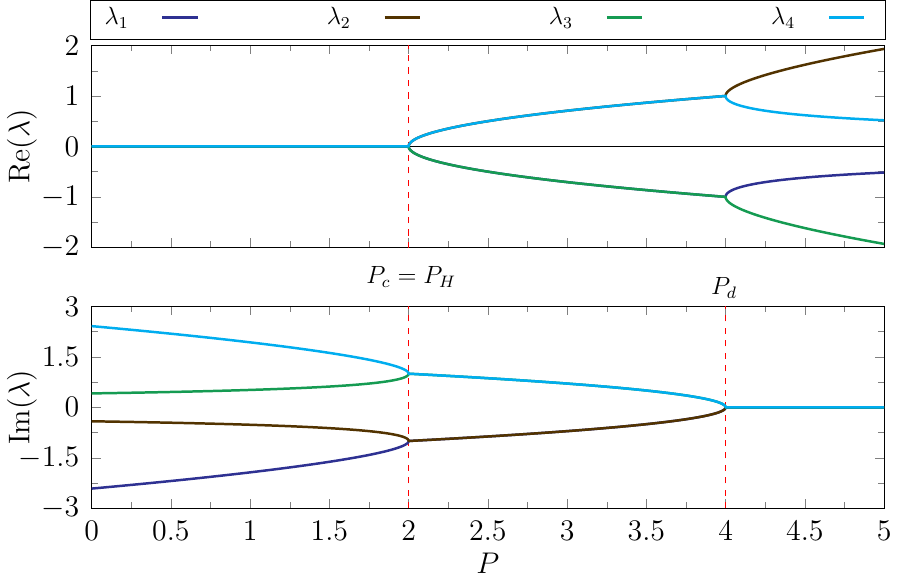}
    \caption{}
    \label{fig:2DOFZiegler_UniParam_eig_undamped}
\end{subfigure}
\begin{subfigure}{0.49\textwidth}
    \centering
    \includegraphics[width=\textwidth]{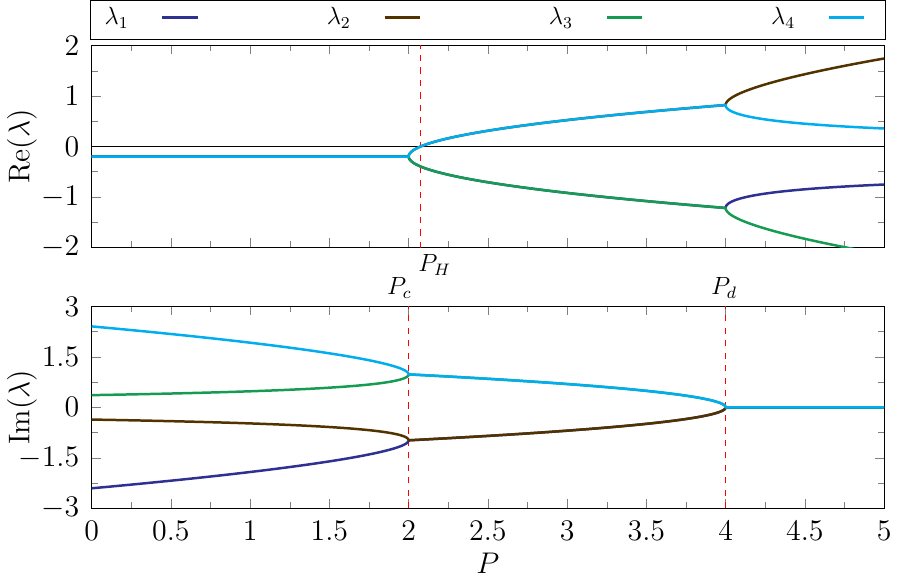}
    \caption{}
    \label{fig:2DOFZiegler_UniParam_eig_mass}
\end{subfigure}
\par\bigskip
\begin{subfigure}{0.49\textwidth}
    \centering
    \includegraphics[width=\textwidth]{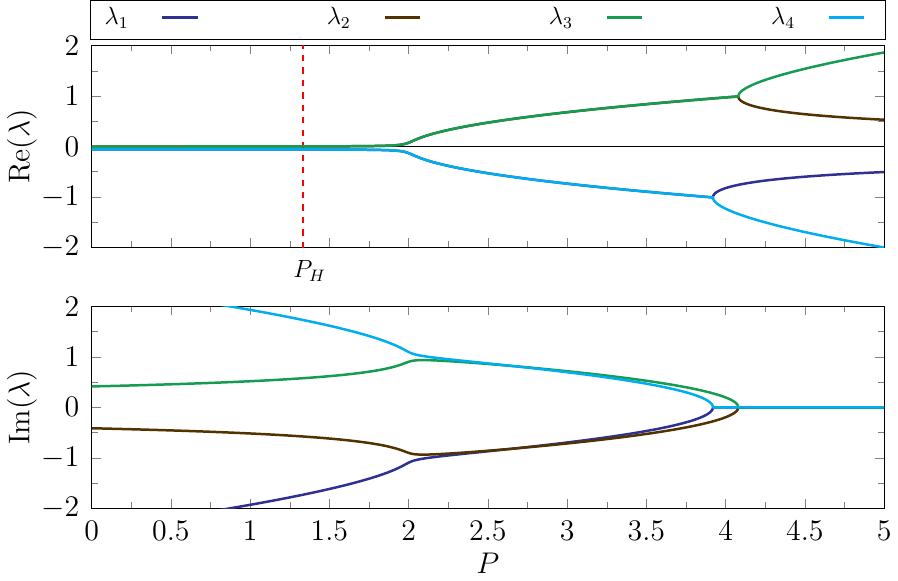}
    \caption{}
    \label{fig:2DOFZiegler_UniParam_eig_stiff_smaller}
\end{subfigure}
\begin{subfigure}{0.49\textwidth}
    \centering
    \includegraphics[width=\textwidth]{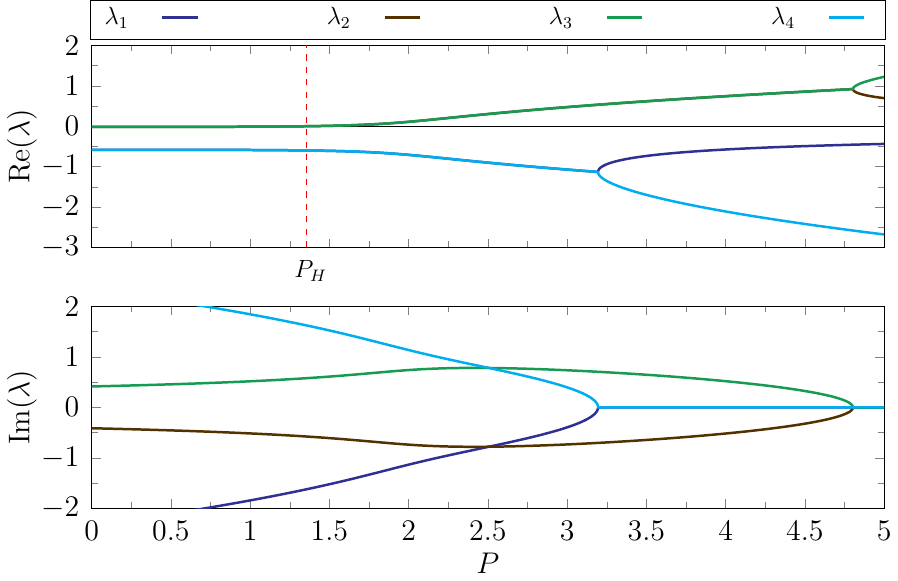}
    \caption{}
    \label{fig:2DOFZiegler_UniParam_eig_stiff_small}
\end{subfigure}
\caption{Eigenvalue trajectories for the 2-DOF Ziegler pendulum. Three distinct scenarios are considered:  without damping, with mass-proportional damping,  and with stiffness-proportional damping. Parameter values set as $m_1=m_2=k_1=k_2=L=1$. (a) undamped system (b) mass-proportional damping with $\xi_m = 0.2$ (c)-(d) cases with stiffness-proportional damping, and $\xi_k=0.01$ or $\xi_k=0.1$. The black, solid line on the top plots indicates zero.}
\label{fig:2DOFZiegler_UniParam_eig}
\end{figure}

Let us illustrate the possible cases and the influence of the damping in the linear stability analysis of the Ziegler pendulum. \cref{fig:2DOFZiegler_UniParam_eig} shows four possible cases with different damping scenarios, where all the free parameters have been set to unity: $m_1=m_2=1$; $k_1=k_2=1$ and $L=1$. \cref{fig:2DOFZiegler_UniParam_eig_undamped} shows the conservative case obtained with $\bfC={\mathbf 0}$. Without the external follower force, $P=0$, the eigenspectrum is composed of a couple of purely imaginary complex eigenvalues ${\pm i\omega_1, \pm i \omega_2}$, with $(\omega_1,\omega_2)$ the two eigenfrequencies of the vibrating system. Increasing the load of the follower force, one can observe that the two eigenfrequencies coalesce at the point $P=P_c$. This point is also the Hopf bifurcation, denoted as $P_H$, since a pair of positive real parts appears. This is the classical scenario for the flutter instability, characterized by an exceptional point at $P_c=P_H$, and the appearance of a Jordan block. It is also interesting to note that the two merged imaginary parts stay exactly equal until the next bifurcation at $P=4$. This underlines a strong 1:1 resonance between the two modes which have equal imaginary (oscillatory) parts, but opposite real parts: one of the two modes is unstable, while the other one is stable. At $P=4$, a divergence instability occurs, and this point is referred to as $P_d$ in the following. Here, the imaginary parts vanish and there is no oscillatory motion anymore, typical of a divergence instability. The range of parameter of interest for our study is thus the whole interval where stable limit cycles develop, {\em i.e.} for $P_H \leq P \leq P_d$ (here: $2\leq P\leq 4$).

Let us now illustrate how this classical scenario is modified when damping is taken into account. \cref{fig:2DOFZiegler_UniParam_eig_mass} first shows how mass-proportional damping affects the bifurcation scenario. Since the control parameter is only involved in the stiffness term $\bfKg$ in~\cref{eq:EOMZiegler2DOF}, in this particular case adding mass-proportional damping just shifts the real parts of the eigenvalues, without fundamentally modifying the bifurcation scenario. The main consequence is that the frequency coalescence is still exactly verified, but the Hopf bifurcation point $P_H$ now occurs for a slightly larger value than $P_c$. All other characteristics are left unchanged: an EP is found at $P=P_c$ and a Jordan block appears, and a 1:1 resonance occurs since the imaginary parts are equal. The divergence is also observed for larger values at $P=P_d$.

This slight modification is only possible for the case of mass-proportional damping, but considering an arbitrary damping matrix will substantially modify the bifurcation scenario. This is illustrated in~\cref{fig:2DOFZiegler_UniParam_eig_stiff_smaller,fig:2DOFZiegler_UniParam_eig_stiff_small}, where increasing values of stiffness-proportional damping are considered. The first case in~\cref{fig:2DOFZiegler_UniParam_eig_stiff_smaller} is that of a small damping value $\xi_k$, selected to stay close to the conservative case.
As compared to~\cref{fig:2DOFZiegler_UniParam_eig_undamped}, one can observe that the coalescence of the eigenfrequencies is lost, and the EP does not exist anymore. The Hopf bifurcation point $P_H$ of course persists and is located for these specific parameter values at $P_H=1.6$. It is interesting to note that only one unstable mode appears at $P=P_H$. However, looking at the eigenvalue trajectories and comparing them to the conservative case, one can see that the 1:1 resonance between the unstable mode (in green) and the stable one emanating from the largest eigenvalue $\omega_2$ for $P=0$, is still almost fulfilled. A near 1:1 resonance exists on the parameter range of interest where fluttering oscillations develop. 

Finally, a large value of the stiffness-proportional damping is shown in~\cref{fig:2DOFZiegler_UniParam_eig_stiff_small}. In that case, the frequency coalescence is completely lost. The 1:1 resonance is not fulfilled anymore, but the imaginary parts still stay near each other with values in the same range. Nevertheless, the Hopf and divergence points are still present. 

Since the parametrisation method for invariant manifold makes important use of the linear characteristics in order to build a simulation-free nonlinear ROM, all the features displayed in these eigenvalue trajectories will be further discussed and analyzed to construct an efficient ROM. This is detailed in the next section.

\subsection{Adjustements to the general reduction technique}\label{subsec:adjustz}

\subsubsection{Two-mode strategy} \label{subsubsec:TwoModes}

In the conservative case, and as illustrated with the Ziegler pendulum, the point where the eigenfrequencies coalesce is the Hopf bifurcation point. At this EP, all the eigenfrequencies are equal with values $\omega_c$ and vanishing real parts. Consequently, the center manifold, which captures the dynamics in the neighbourhood of the nonhyperbolic fixed point~\cite{Wiggins}, is here four-dimensional. After the bifurcation, the four eigenvalues involved in the coalescence scenario have opposite real and imaginary parts and read $\lambda_j =\pm r_c \pm i \omega_c$, $j \in \{1,4\}$, see {\em e.g.}~\cref{fig:2DOFZiegler_UniParam_eig_undamped}.
One strategy to derive a ROM is to select only the mode corresponding to the pair of eigenvalues with positive real parts $\lambda_j = + r_c \pm i \omega_c$. This strategy has been proposed in~\cite{MingwuLi2024}, and will be referred to as the one-mode ROM in the remainder. However, a clear 1:1 resonance exists between the oscillatory parts, and this feature persists, in the conservative case, as long as the next bifurcation point is encountered. Therefore, our proposal is to consider the 2 modes (corresponding to the four eigenvalues  $\lambda_j =\pm r_c \pm i \omega_c$
) to construct the ROM with the DPIM. It is expected that taking the 1:1 internally resonant behaviour will lead to more accurate predictions of the post-critical regime as well as the bifurcation point.

This statement is only strictly valid as long as the system is conservative. However, for lightly damped systems, near-resonances occur and the 1:1 resonance between the imaginary parts is not completely lost but approximately true. Hence, the two-mode strategy could still be helpful in the construction of the ROM for general situations. Nevertheless, It should be noted that its applicability for damped scenarios depends on the underlying conservative system at hand, in order to identify the two coalescing modes, which might not always be possible. Therefore, in the next sections, this scheme will be tested and compared with its one-mode counterpart. 

In the remainder of the study, the complex normal form style in the parametrisation method, is selected~\cite{vizza21high,TouzeCISM2,Stabile:morfesym}. In order to ensure that the monomials corresponding to the 1:1 resonance are present even in the case of quasi-resonance, the exact 1:1 relationship between the eigenfrequencies is strictly enforced in the calculation. Another strategy could be to use the graph style parametrisation.

\subsubsection{Inclusion of Jordan blocks} \label{subsubsec:Jordan}

As mentioned in \cref{subsec:LinearStability}, once exceptional points are at hand, the linear part of the dynamics can no longer be diagonalized, and Jordan blocks appear in the Jacobian matrix.

From a numerical standpoint, when usual routines for eigenproblems are employed, there will never be two perfectly identical eigenvalues, so the resulting eigenvalues matrix is still diagonal, with two of its diagonal elements almost equal (for the sake of simplicity, we will assume that the maximum eigenvalues multiplicity is 2). These almost identical eigenvalues have a two-dimensional eigenspace, and correspond to almost identical eigenvectors for the exceptional point case. However, working directly with them would render the following computations ill-conditioned since they are almost perfectly aligned. 

Therefore, an algorithm to create a better-conditioned eigensubspace starting from the left and right eigenvectors obtained by usual numerical routine, should be devised and expected to improve the numerical outputs of the reduction strategy. The main idea is to artificially impose Jordan blocks in the Jacobian matrix and find modified eigenvectors orthogonal to one another for each eigenspace. 

The details of this technical derivation are reported in~\cref{app:camillejordan}, for the sake of brevity. Such a treatment is awaited to provide better results in the reduction technique only in the case where an exceptional point is met. Numerical illustrations are provided in the next section.

\subsubsection{Selection of the expansion point to compute the ROM}\label{subsubsec:expoint}

As explained in~\cref{subsec:Reduction}, the ROM is constructed with the DPIM at the expansion point $\mu=\mu_0$. In the case of bifurcating systems like the Ziegler pendulum investigated before, different options are possible to select the $\mu_0$ value at which the ROM is computed. It could be before the instability, exactly at the Hopf bifurcation point, or after the instability. 

In the numerical examples, the consequences of this choice will be investigated and further commented on. In particular, it will be shown that, thanks to the two-mode strategy, the expansion point can be selected before the bifurcation and the ROM is then able to predict the value of the bifurcation point. If one is interested in an accurate prediction of the amplitudes of the limit cycles developing after the Hopf bifurcation, it will be shown that the best results are obtained when the ROM is computed for $\mu_0$ values that are slightly larger than the bifurcation point, {\em i.e.} when the fixed point is unstable. This result, already commented on in~\cite{MingwuLi2024}, will be here analyzed in terms of phase space representations.

\section{The Ziegler pendulum} \label{sec:Ziegler}

In this Section, the results of the different reduction strategies are compared for the case of the Ziegler pendulum. Since the two-mode strategy proposed to improve the results leads to keeping the same number of coordinates in the ROM for the 2-DOF Ziegler, the results are then extended to a 3-DOF pendulum. A continuous problem is investigated in~\cref{sec:Beck}.

\subsection{A 2-DOF Ziegler pendulum} \label{subsec:2DOFZiegler}

In this section, the 2-DOF Ziegler pendulum, shown in \cref{fig:2DOFZiegler} and introduced previously to illustrate some key concepts, is used to test the proposed approach. The length of the bars is fixed to $L=1$. The stiffness and masses are not taken here as unity, contrary to the results shown in~\cref{subsec:LinearStability}. Instead, their values have been optimized to repel the divergence stability point far from the Hopf bifurcation, to ensure the largest possible range of control parameter values for which limit cycle oscillations develop. This leads to select the following parameter values
\begin{equation} \label{eq:2DOFZieglerAleParam}
     k_1 = \delta^2 k_2, \quad m_1 = \gamma^2 m_2, \quad k_2 = m_2 = 1; \;  \mbox{with}\quad \gamma^2 = \nicefrac{25}{4}, \; \mbox{and} \quad \delta^2 = \nicefrac{41}{4}.
\end{equation}

\begin{figure}[h!]
\centering
\begin{subfigure}{0.325\textwidth}
    \centering
    \includegraphics[width=\textwidth]{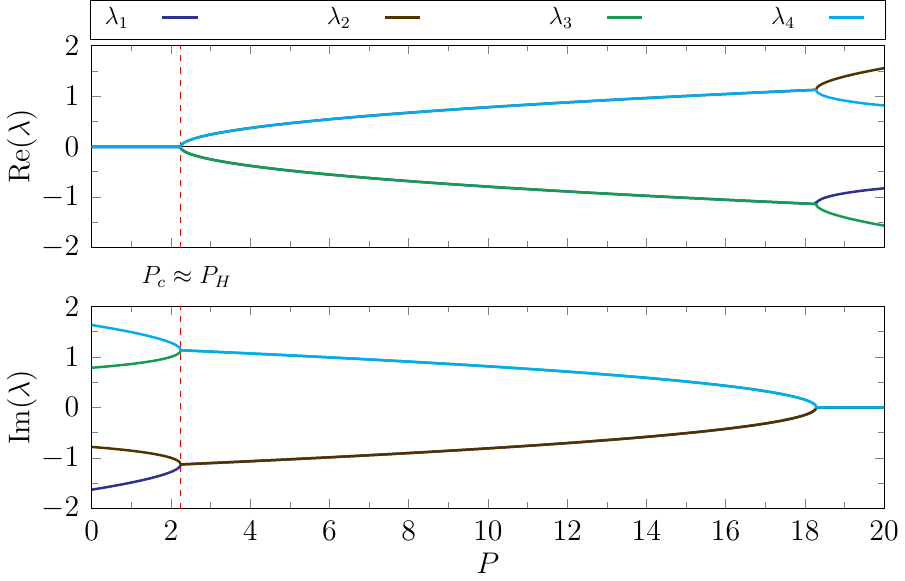}
    \caption{$\xi_m = 0.01$ and $\xi_k = 0$.}
    \label{fig:2DOFZiegler_AleParam_eig_mass}
\end{subfigure}
\begin{subfigure}{0.325\textwidth}
    \centering
    \includegraphics[width=\textwidth]{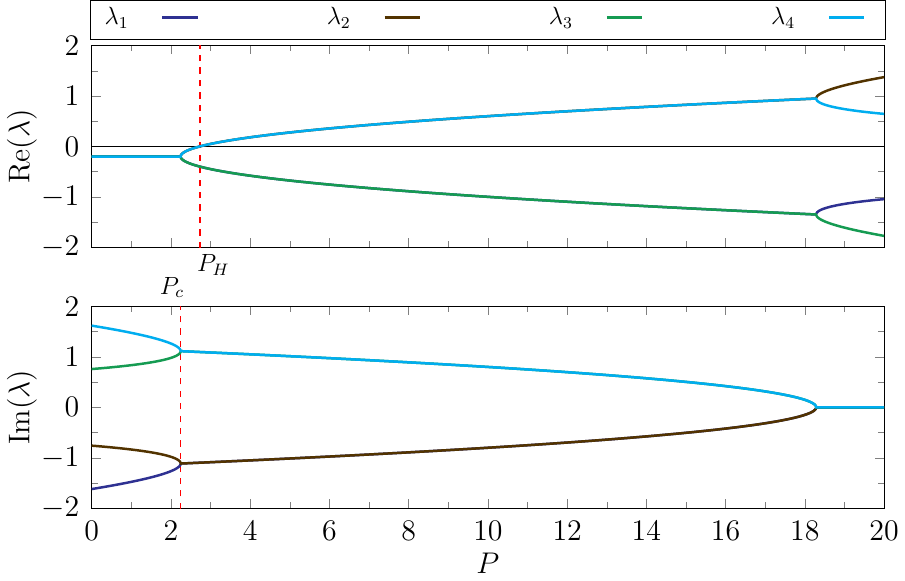}
    \caption{$\xi_m = 0.2$ and $\xi_k = 0$.}
    \label{fig:2DOFZiegler_AleParam_eig_massHigh}
\end{subfigure}
\begin{subfigure}{0.325\textwidth}
    \centering
    \includegraphics[width=\textwidth]{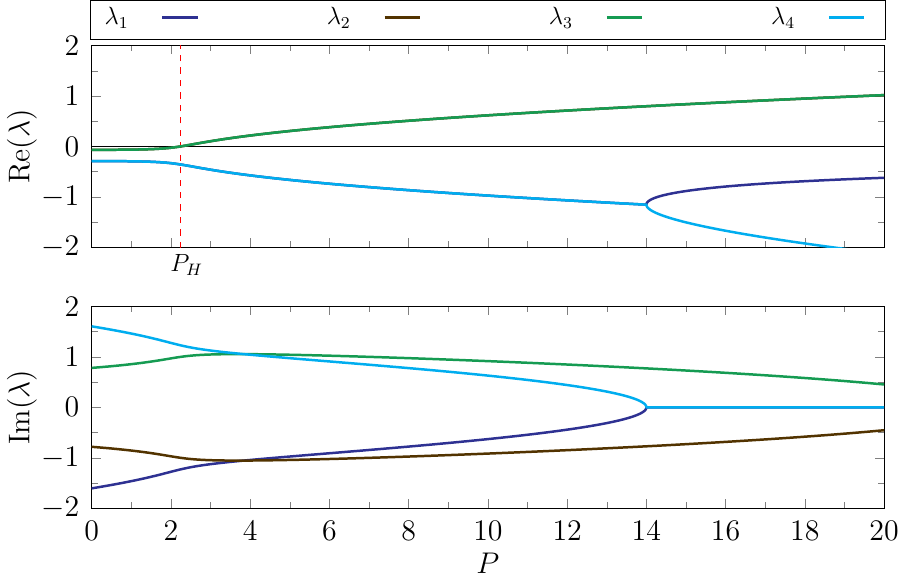}
    \caption{$\xi_m = 0$ and $\xi_k = 0.1$.}
    \label{fig:2DOFZiegler_AleParam_eig_stiff}
\end{subfigure}
\caption{Eigenvalue trajectories for the 2-DOF Ziegler pendulum with $L=1$ and other parameters as given by \cref{eq:2DOFZieglerAleParam}.}
\label{fig:2DOFZiegler_AleParam_eig}
\end{figure}

Regarding damping, three scenarios are considered: in two of them, mass-proportional damping is selected, with $\xi_m = 0.01$ and $\xi_m = 0.2$, whereas in the third stiffness-proportional damping with $\xi_k=0.1$ is considered. The eigenvalue trajectories for these three cases are illustrated in \cref{fig:2DOFZiegler_AleParam_eig}. The same general comments made on \cref{subsec:LinearStability} apply to this example. It should be noted that from \cref{fig:2DOFZiegler_AleParam_eig_mass} to \cref{fig:2DOFZiegler_AleParam_eig_massHigh} the increase in damping moves the points where the eigenvalues coalescence and the Hopf bifurcation happens further apart. Also, it should be noted that the divergence bifurcation in \cref{fig:2DOFZiegler_AleParam_eig} occurs at much larger values of the control parameter than for the systems in \cref{fig:2DOFZiegler_UniParam_eig}. 

\begin{figure}[H]
    \centering
    \begin{subfigure}{0.49\textwidth}
        \centering
        \includegraphics[width=\textwidth]{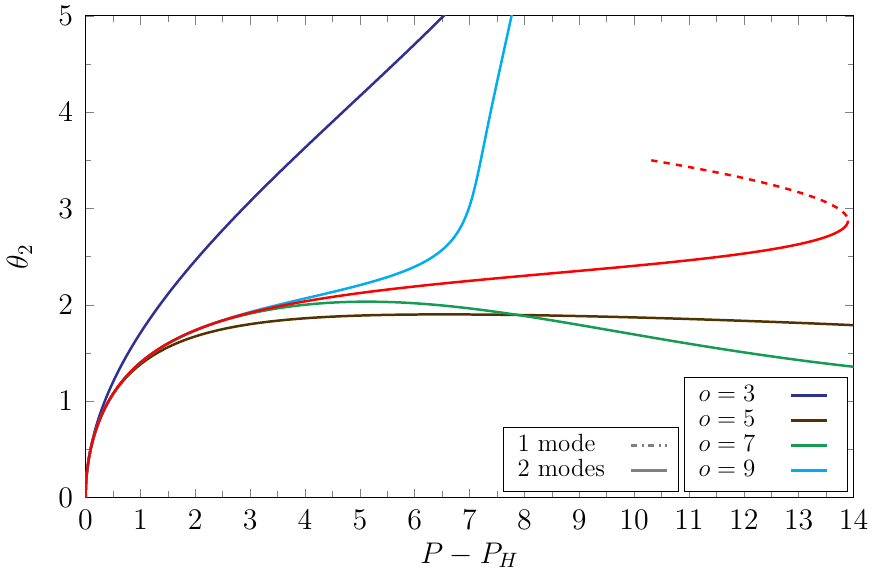}
        \caption{}
        \label{fig:2DOFZiegler_BifDiagrams_Mass_BP}
    \end{subfigure}
    \begin{subfigure}{0.49\textwidth}
        \centering
        \includegraphics[width=\textwidth]{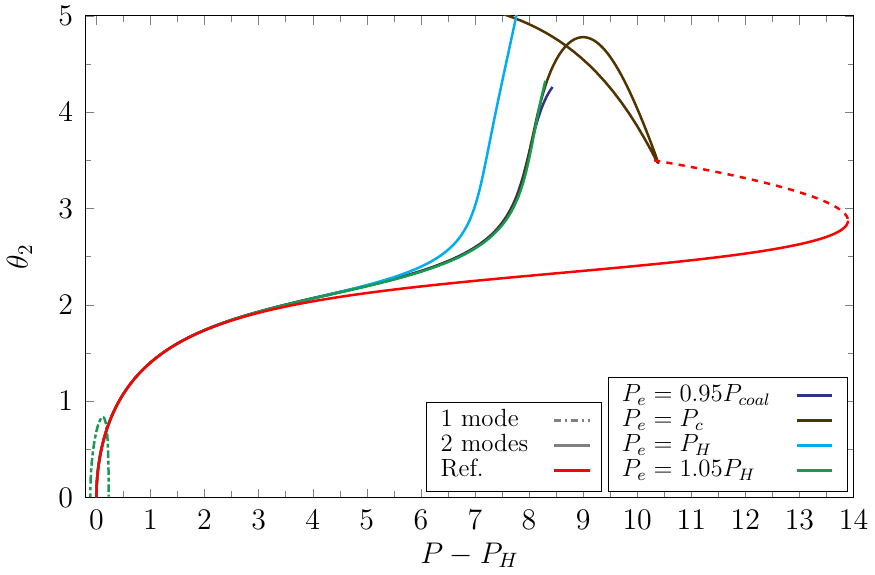}
        \caption{}
        \label{fig:2DOFZiegler_BifDiagrams_Mass_o9}
    \end{subfigure}
    \vspace{-15pt} \vphantom{a} {\scriptsize $\xi_m = 0.01, \xi_k = 0$.} \vspace{6pt}

    \begin{subfigure}{0.49\textwidth}
        \centering
        \includegraphics[width=\textwidth]{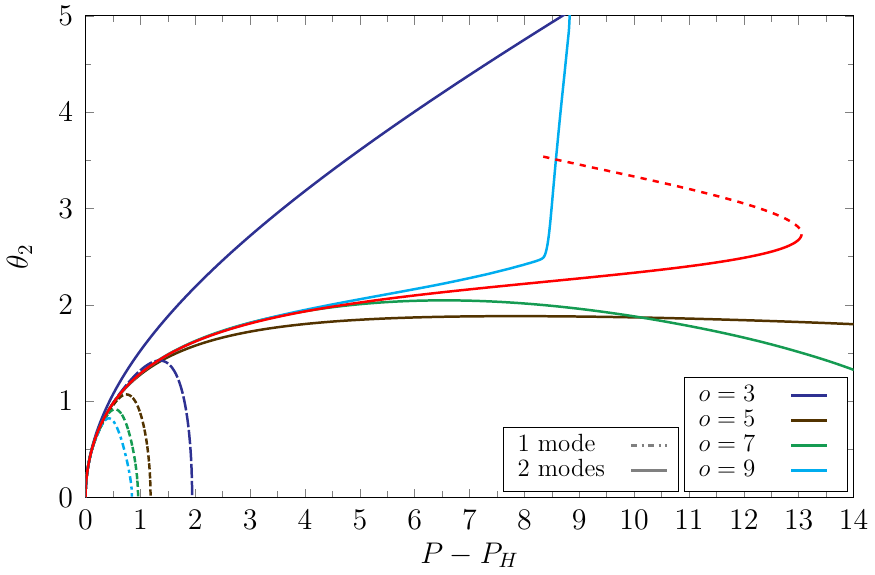}
        \caption{}
        \label{fig:2DOFZiegler_BifDiagrams_MassHigh_BP}
    \end{subfigure}
    \begin{subfigure}{0.49\textwidth}
        \centering
        \includegraphics[width=\textwidth]{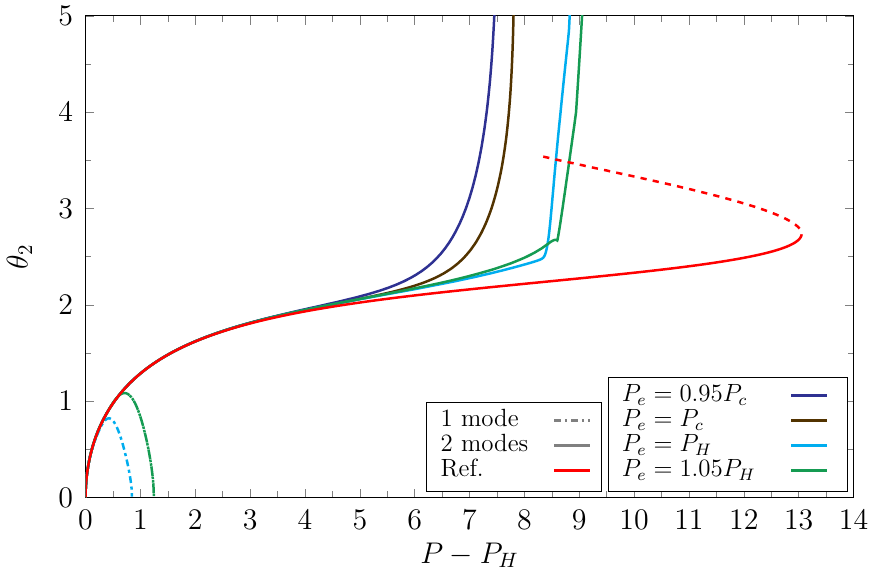}
        \caption{}
        \label{fig:2DOFZiegler_BifDiagrams_MassHigh_o9}
    \end{subfigure}
    \vspace{-15pt} \vphantom{a} {\scriptsize $\xi_m = 0.2, \xi_k = 0$.} \vspace{6pt}
    
    \begin{subfigure}{0.49\textwidth}
        \centering
        \includegraphics[width=\textwidth]{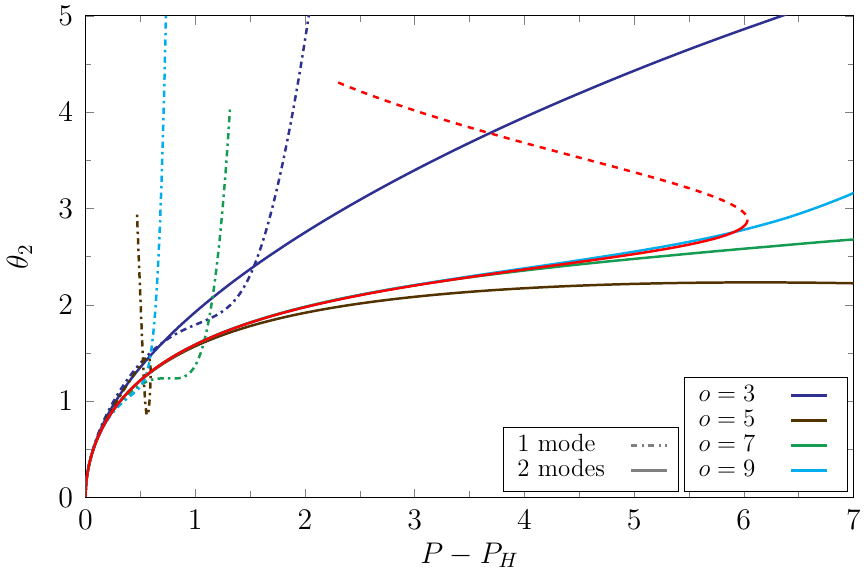}
        \caption{}
        \label{fig:2DOFZiegler_BifDiagrams_Stiff_BP}
    \end{subfigure}
    \begin{subfigure}{0.49\textwidth}
        \centering
        \includegraphics[width=\textwidth]{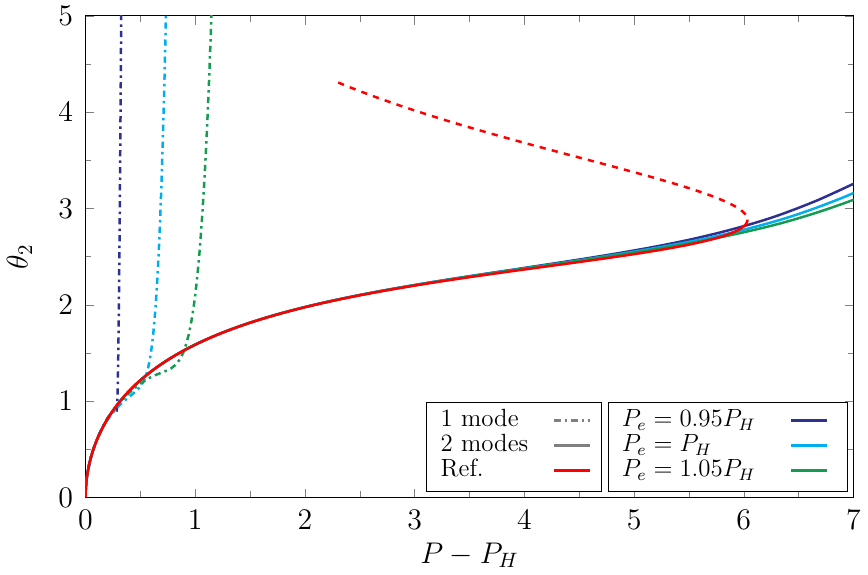}
        \caption{}
        \label{fig:2DOFZiegler_BifDiagrams_Stiff_o9}
    \end{subfigure}
    \vspace{-15pt} \vphantom{a} {\scriptsize $\xi_m = 0, \xi_k = 0.1$.} \vspace{6pt}

    \caption{Bifurcation diagrams for the 2-DOF Ziegler pendulum with $L=1$ and other parameters as given by \cref{eq:2DOFZieglerAleParam}. The amplitude of the limit cycle for $\theta_2$ is given as a function of $P-P_H$. $P_e$ denotes the parametrisation point, $P_c$ the point of eigenfrequencies coalescence and $P_H$ the Hopf bifurcation point. (a), (c) and (e): solutions for increasing orders when the parametrisation is computed for $P_e=P_H$. (b), (d) and (f): solutions for a fixed order 9 of the parametrisation and different values for the expansion point $P_e$.  Each line corresponds to a different damping scenario. Full-order solutions computed by numerical continuation implemented in the package Matcont \cite{dhooge2004matcont} in red. Dashed line for unstable solutions. Stability is reported for the full-order model only.}
    \label{fig:2DOFZiegler_BifDiagrams}
\end{figure}

For the analysis, parametrisations with increasing orders from 3 to 9 are considered. Furthermore, the expansion point $P_e$ around which the parametrisation is computed, is also varied. 
Note that this expansion point was denoted as $\mu_0$ in Section~\ref{subsec:Reduction} presenting the general methodology, but is now referred to as $P_e$ since $P$ is the bifurcation parameter in this problem.
For systems with mass-proportional damping, the coalescence point $P_c$ and the Hopf bifurcation point $P_H$ are chosen, as well as expansions around $0.95 P_c$ and $1.05 P_H$. As for the system with stiffness-proportional damping, points $0.95 P_H$, $P_H$ and $1.05 P_H$ are also selected to compute the ROM. Bifurcation diagrams for the different scenarios are computed thanks to a numerical continuation procedure embedded in the package Matcont \cite{dhooge2004matcont}. \cref{fig:2DOFZiegler_BifDiagrams} summarizes the obtained results where a full-order solution taken as reference is shown with a red line. \cref{fig:2DOFZiegler_BifDiagrams_Mass_BP,fig:2DOFZiegler_BifDiagrams_MassHigh_BP,fig:2DOFZiegler_BifDiagrams_Stiff_BP} show ROMs computed around the bifurcation point with different expansion orders, while \cref{fig:2DOFZiegler_BifDiagrams_Mass_o9,fig:2DOFZiegler_BifDiagrams_MassHigh_o9,fig:2DOFZiegler_BifDiagrams_Stiff_o9} show ROMs computed at different expansion points but with a fixed expansion order, equal to 9.

The main point of interest is to investigate if the two adjustments to the method proposed in~\cref{subsec:adjustz} improve the ROM's ability to provide an accurate prediction of the limit cycles' amplitudes. As a first remark, the ROM computed by retaining only the unstable mode has been found unable to detect the Hopf bifurcation point and compute the bifurcated branch if the expansion point is selected before or at the coalescence point. On the other hand, the two-mode strategy allows finding the Hopf bifurcation point in any case, which thus constitutes a first advantage. Second, one can observe as a general trend for all the results reported in~\cref{fig:2DOFZiegler_BifDiagrams} that adding the companion mode in the ROM leads to an important improvement. Indeed, in most cases, using only the unstable mode in the reduction leads either to no bifurcated branches found, or to approximations that are valid for small parameter variations. It should however be noted that when two master modes are selected, no reduction of the system actually takes place, which might be relevant to the results. This fact will be commented further in \cref{subsec:3DOFZiegler} and~\cref{sec:Beck}, where it will be highlighted that the results found here still hold for systems with larger dimension where a reduction occurs when selecting two master modes.

Regarding the strategy to enforce Jordan blocks proposed in~\cref{subsubsec:Jordan}, numerical results underline that this adjustment is meaningful and needed only when an exact coalescence is found and the parametrisation is performed around the EP. Indeed, in such cases, ROMs computed without enforcing Jordan blocks could not compute the bifurcated branch. On the other hand, in all other cases where the expansion point is not taken to be the EP, enforcing the Jordan blocks has no influence on the results and is not needed.

Considering now the different tested orders for the ROMs, it can be seen from \cref{fig:2DOFZiegler_BifDiagrams_Mass_BP,fig:2DOFZiegler_BifDiagrams_MassHigh_BP,fig:2DOFZiegler_BifDiagrams_Stiff_BP} that an order 9 parametrisation seems sufficient to provide a comfortable range of parameter variation where the limit cycles' amplitudes are correctly predicted. Also, from \cref{fig:2DOFZiegler_BifDiagrams_Mass_o9,fig:2DOFZiegler_BifDiagrams_MassHigh_o9,fig:2DOFZiegler_BifDiagrams_Stiff_o9}, it is possible to see that when 2 master modes are kept the choice of the expansion does not influence importantly the quality of the ROM. When only 1 mode is retained, it seems more beneficial to parametrise after the bifurcation point, a result that is in agreement with \cite{MingwuLi2024}.

\begin{figure}
    \centering
    \begin{subfigure}{0.49\textwidth}
        \centering
        \includegraphics[width=\textwidth]{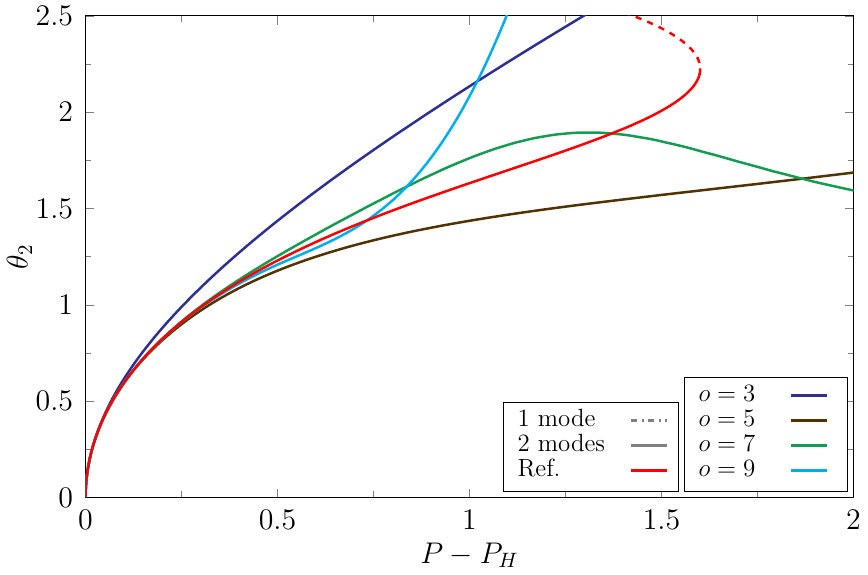}
        \caption{}
        \label{fig:2DOFZiegler_BifDiagrams_UniParam_MassHigh_BP}
    \end{subfigure}
    \begin{subfigure}{0.49\textwidth}
        \centering
        \includegraphics[width=\textwidth]{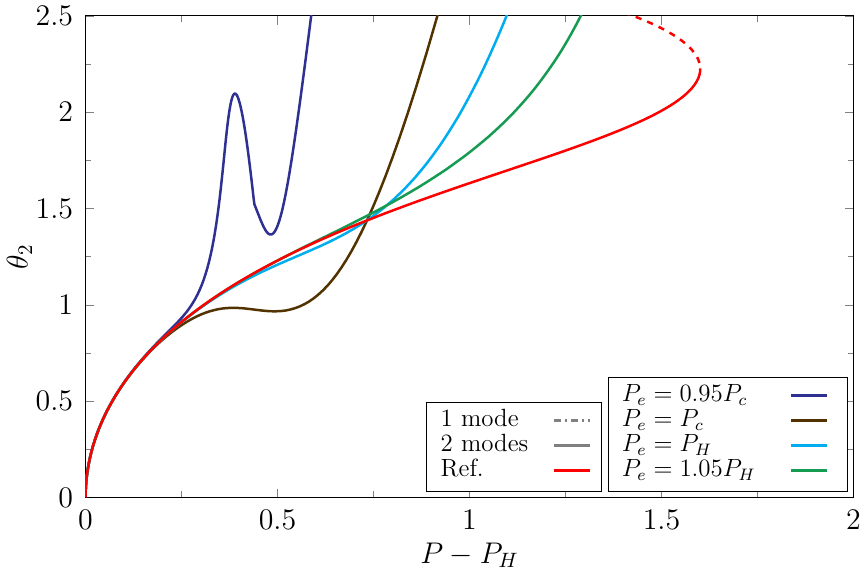}
        \caption{}
        \label{fig:2DOFZiegler_BifDiagrams_UniParam_MassHigh_o9}
    \end{subfigure}
    \caption{Bifurcation diagrams for the 2-DOF Ziegler pendulum with $m_1=m_2=k_1=k_2=L=1$, $\xi_k = 0$ and $\xi_m = 0.2$. The amplitudes of the limit cycles for $\theta_2$ are shown as a function of $P-P_H$. $P_e$ denotes the parametrisation point, $P_c$ the point of eigenfrequencies coalescence and $P_H$ the Hopf bifurcation point. (a) Solutions for increasing orders when $P_e=P_H$. (b) Solutions for fixed order 9 and different expansion points.}
    \label{fig:2DOFZiegler_BifDiagrams_UniParam}
\end{figure}

Another interesting comment concerns the effect of other bifurcation points on the validity limit of the ROMs. Since the method is fundamentally local in nature, and an expansion point around the Hopf bifurcation is selected, the quality of the approximation should worsen as the control parameter approaches values at which other bifurcations occur. An example of this statement can be seen for the Ziegler pendulum, as the system approaches the divergence bifurcation occurring for larger values of the load. In order to illustrate this, the case depicted by \cref{fig:2DOFZiegler_UniParam_eig_mass}, where all parameters have unitary values, is considered. The bifurcation diagrams obtained for this case are reported in~\cref{fig:2DOFZiegler_BifDiagrams_UniParam}. Comparing this case with the previous one shown in~\cref{fig:2DOFZiegler_BifDiagrams}, where the parameters have been selected to repel the divergence bifurcation far from the Hopf one, the validity limit appears to be smaller in~\cref{fig:2DOFZiegler_BifDiagrams_UniParam}. Indeed, denoting by $P_v$ the maximum load value for which a good agreement between the ROM and the reference solution can be obtained, one finds by visual inspection that $P_v-P_H \approx 0.5$ in~\cref{fig:2DOFZiegler_BifDiagrams_UniParam_MassHigh_o9} and $P_v-P_H \approx 4$ in~\cref{fig:2DOFZiegler_BifDiagrams_MassHigh_o9}. However, these values should only be meaningful when compared to the size of the interval between the divergence and the Hopf bifurcations. Defining $P_d$ as the divergence load, this length is equal to $P_d-P_H \approx 2$ in the first case and $P_d-P_H \approx 16$ in the second. Therefore, the ratio $\nicefrac{P_{v}-P_H}{P_d-P_H}$ remains approximately constant for both situations. For the other examples in this contribution, parameters such that the divergence bifurcation is repelled will be chosen whenever possible, in order to have a larger absolute range of validity of the expansion. 

\subsection{Geometrical interpretation in phase space}

This section aims to justify the choice of expansion point to compute the ROM {\it after} the bifurcation point. With such a choice, the method computes the unstable manifold of the unstable fixed point, and it is found that this ROM correctly predicts the limit cycle location in phase space after the Hopf bifurcation. 

A geometrical interpretation for the 2-DOF Ziegler pendulum is provided using partial 3-dimensional views of the complete 4-dimensional phase space. A key idea in dynamical systems theory is that the unstable manifold of a saddle point is generally connected to another stable fixed point via a heteroclinic orbit (see {\it e.g.}~\cite{gucken83,Wiggins} or~\cite{DauchotMann,DOEDEL2011222} for illustrations) or spirals to a limit cycle if the associated stable attractor is of that kind (see {\it e.g.} the case of van der Pol oscillator or the predator-prey model). In the same line and in the case of chaotic dynamics, the strange attractor is contained in the closure of the unstable manifold~\cite{BowenRuelle,EckmanRuelle,Shilnikov}.

If an exact computation of the unstable manifold is achievable, in our present case it will always converge to the stable limit cycle and thus provide an excellent prediction for the amplitudes of the bifurcated solutions. However, since asymptotic expansions are used in the calculations, and since the limit cycle moves further away from the unstable fixed point when the bifurcation parameter is varied, the validity limit of the approximation will, at some point, be smaller than the distance from the fixed point to the limit cycle, and the method will then fail in providing an accurate prediction.

\begin{figure}[h]
    \centering
    \begin{subfigure}{0.49\textwidth}
        \includegraphics[width=\textwidth]{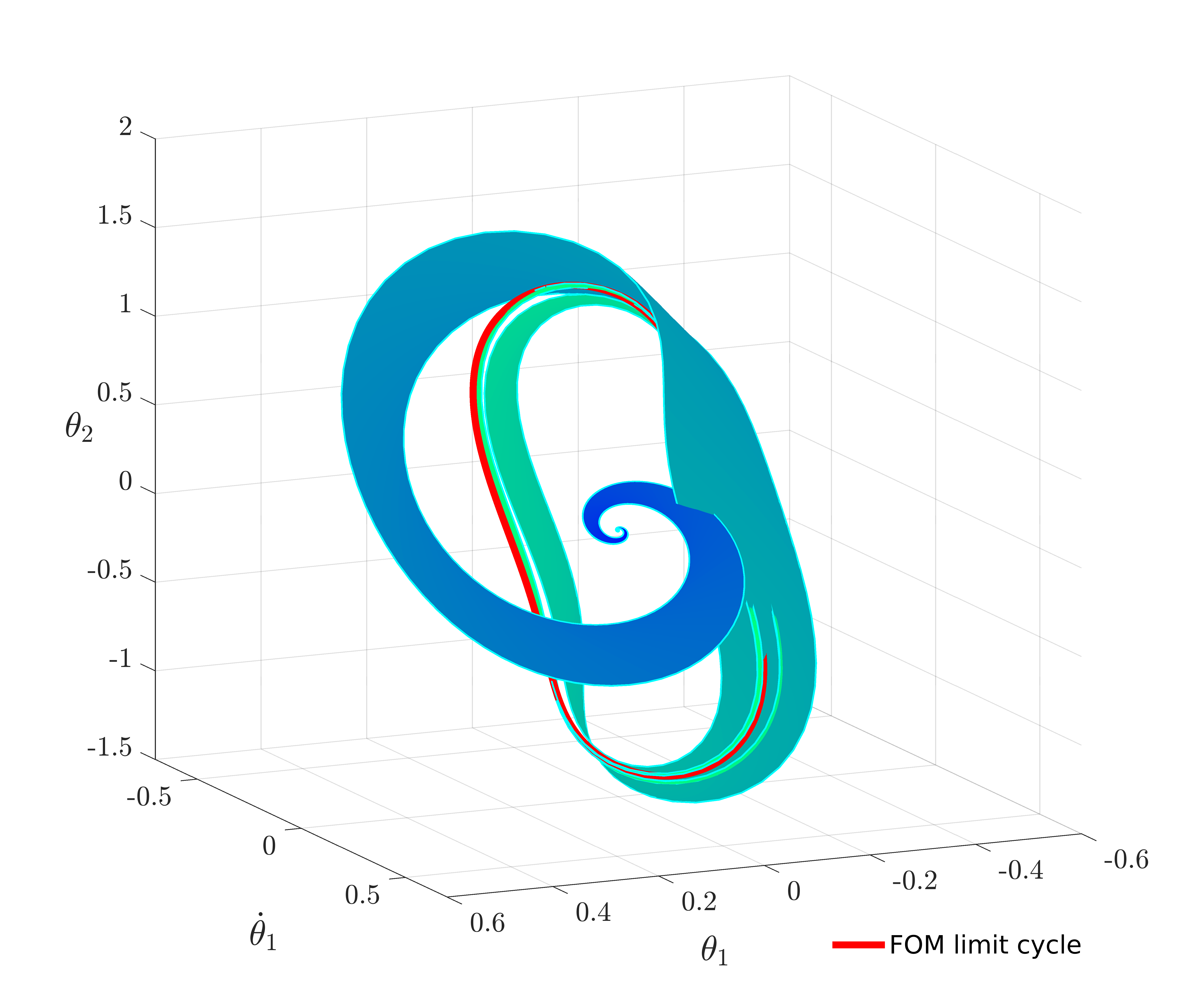}
        \caption{}
        \label{fig:Unstable_manifold_dP=0.5}
    \end{subfigure}
    \begin{subfigure}{0.49\textwidth}
        \includegraphics[width=\textwidth]{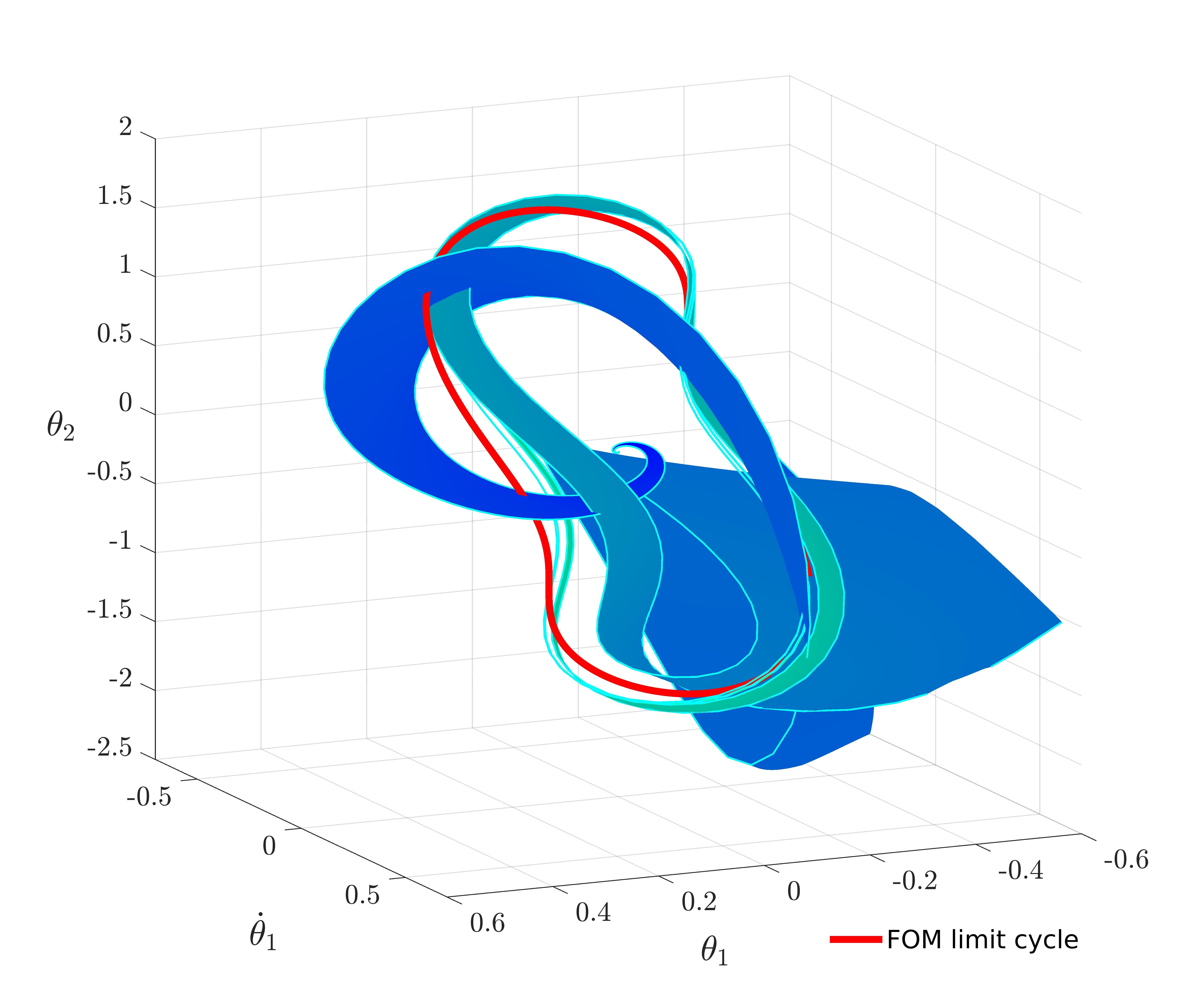}
        \caption{}
        \label{fig:Unstable_manifold_dP=0.8}
    \end{subfigure}

    \caption{Unstable manifolds for the 2-DOF Ziegler pendulum with $m_1=m_2=k_1=k_2=L=1$, $\xi_k = 0$, $\xi_m = 0.2$ and expansion about $1.05P_H$. (a) and (b) show the unstable manifolds obtained from the reduced order model for $P-P_H = 0.5$ and $P-P_H = 0.8$, respectively. The stable limit cycle emanating from the Hopf bifurcation, obtained by continuation of the full-order model (FOM), is also shown in red. The colormap ranges from blue to green with increasing time if a trajectory starting at the fixed point and going to the limit cycle is considered. The borders of the manifold, obtained by numerical integration of the ROM considering the limits of the family of initial conditions for determining the unstable manifold, are shown in cyan.}
    \label{fig:Unstable_manifolds2}
\end{figure}

This phenomenon is illustrated in \cref{fig:Unstable_manifolds2}, where the 2-DOF Ziegler pendulum with all coefficients taken equal to 1 is considered, see the results shown in \cref{fig:2DOFZiegler_BifDiagrams_UniParam}. \cref{fig:Unstable_manifolds2} considers a ROM computed after the bifurcation point: the expansion point retained is at 1.05$P_H$, corresponding to the cyan curve in \cref{fig:2DOFZiegler_BifDiagrams_UniParam_MassHigh_o9}. The unstable manifold computed with the order~9 expansion is compared to the exact limit cycle obtained from the full model in the space $(\theta_1, \dot{\theta}_1,\theta_2)$. The unstable manifold of the ROM has been computed by selecting a family of initial conditions spanned by the master eigenvectors. These initial conditions are then integrated in time with the reduced dynamics, and the nonlinear mapping \cref{eq:ziNLmap00} is then used to come back to the full phase space. By doing so, the approximated unstable manifold obtained from the ROM can be shown as a surface. A color code corresponding to the time of the trajectories is used to highlight the manifold.

Two values of the bifurcation parameter are tested: in \cref{fig:Unstable_manifold_dP=0.5}, the case $P-P_H=0.5$ is used, for which the ROM has been found to give an excellent prediction of the limit cycle. The figure shows that, as expected, the order 9 approximation of the unstable manifold exactly converges to the exact limit cycle of the full system. On the other hand, \cref{fig:Unstable_manifold_dP=0.8} considers a larger value of the bifurcation parameter: $P-P_H=0.8$. As reported in \cref{fig:2DOFZiegler_BifDiagrams_UniParam_MassHigh_o9}, in this case the prediction of the ROM departs from the full order solution, and this is retrieved in the phase space representation. As awaited, the limit cycle moved further from the unstable fixed point, and the order 9 approximation fails to exactly converge to the full-order limit cycle, as ascertained by the gap between the periodic orbit and the approximated unstable manifold in \cref{fig:Unstable_manifold_dP=0.8}.

This illustrates the gain in parametrizing the ROM after the bifurcation point, since the unstable manifolds are able to provide a correct approximation of the stable attractor of the system which is here the searched limit cycle. On the other hand, the asymptotic expansion, being a local theory, has a validity limit and fails to uniformly converge to the exact limit cycle for large values of the bifurcation parameter, when the limit cycle is too far from the unstable fixed point where the ROM is effectively computed. This bound gives a validity limit of the method that could not be outperformed by simply augmenting the expansion order of the asymptotic approximation.

\subsection{A 3-DOF Ziegler pendulum} \label{subsec:3DOFZiegler}

In this section, a 3-DOF version of the classical Ziegler pendulum is studied. The system is depicted in \cref{fig:3DOFZiegler}. It is an extension of the pendulum shown in \cref{fig:2DOFZiegler} with an additional bar CD, a spring of stiffness $k_3$ and a mass $m_3$.

A detailed derivation of the equations of motion for this system is shown in \cref{ap:3DOFZieglerEOM}. Their final form is the same as \cref{eq:EOMZiegler2DOF}, with matrices and vectors now given by
\begin{equation} \label{eq:Matrices3DOFZiegler}
\begin{gathered}
    \bfM = L^2
    \begin{bmatrix}
        m_1+m_2+m_3 & m_2+m_3 & m_3 \\
        m_2+m_3 & m_2+m_3 & m_3 \\
        m_3 & m_3 & m_3
    \end{bmatrix}
    \quad 
    \bfK =
    \begin{bmatrix}
        k_1+k_2 & -k_2 & 0 \\
        -k_2 & k_2+k_3 & -k_3 \\
        0 & -k_3 & \phantom{-}k_3
    \end{bmatrix}
    \\
    \bfKg = PL 
    \begin{bmatrix}
        -1 & \phantom{-}0 & 1 \\
        \phantom{-}0 & -1 & 1 \\
        \phantom{-}0 & \phantom{-}0 & 0
    \end{bmatrix}
    \quad
    \thetavec =
    \begin{bmatrix}
        \theta_1 \\
        \theta_2 \\
        \theta_3
    \end{bmatrix}
    \quad
    \mathbf{F}_{nl} = 
    -\frac{PL}{6}
    \begin{bmatrix}
        (\theta_1-\theta_3)^3 \\
        (\theta_2-\theta_3)^3 \\
        0
    \end{bmatrix}.
\end{gathered}
\end{equation}

\begin{figure}[h]
\centering
\begin{subfigure}{0.325\textwidth}
    \centering
    \includegraphics[width=\textwidth]{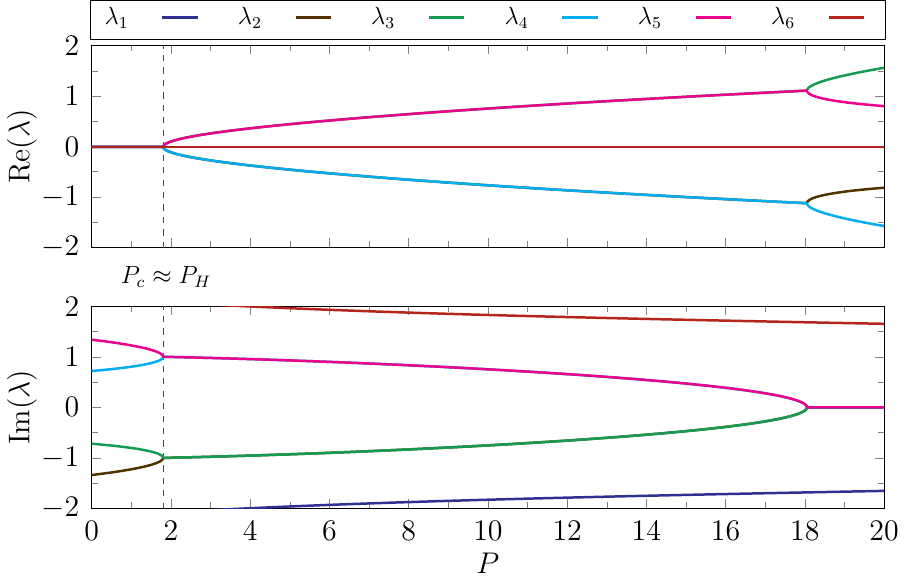}
    \caption{$\xi_m = 0.01$ and $\xi_k = 0$}
    \label{fig:3DOFZiegler_AleParam_eig_mass}
\end{subfigure}
\begin{subfigure}{0.325\textwidth}
    \centering
    \includegraphics[width=\textwidth]{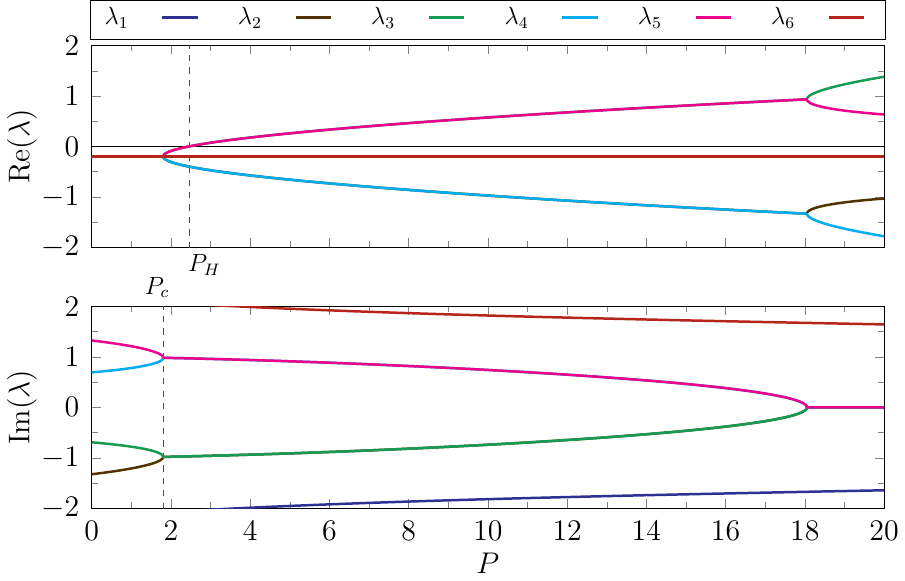}
    \caption{$\xi_m = 0.2$ and $\xi_k = 0$}
    \label{fig:3DOFZiegler_AleParam_eig_massHigh}
\end{subfigure}
\begin{subfigure}{0.325\textwidth}
    \centering
    \includegraphics[width=\textwidth]{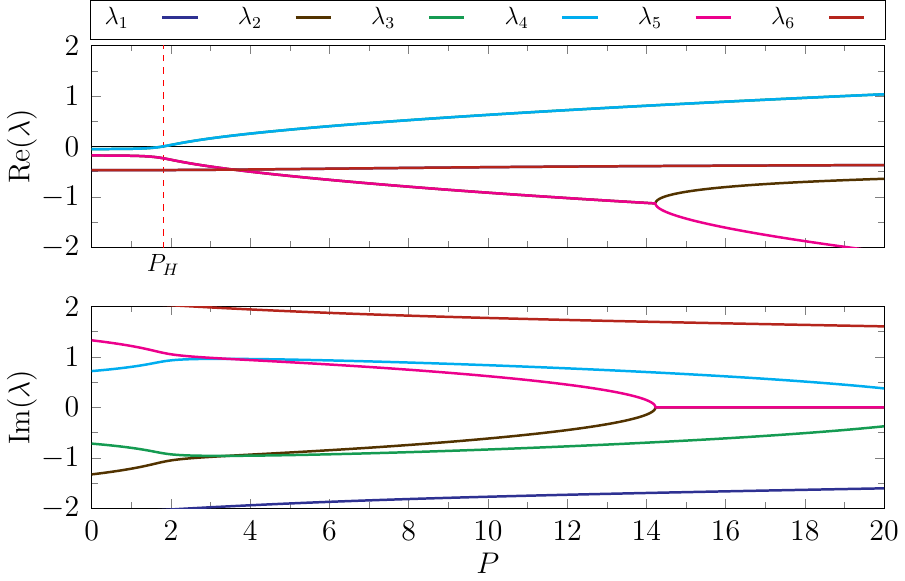}
    \caption{$\xi_m = 0$ and $\xi_k = 0.1$}
    \label{fig:3DOFZiegler_AleParam_eig_stiff}
\end{subfigure}
\caption{Eigenvalue trajectories for the 3-DOF Ziegler pendulum with $L=1$ and other parameters as given by \cref{eq:3DOFZieglerAleParam}.}
\label{fig:3DOFZiegler_AleParam_eig}
\end{figure}

Following the discussion in~\cref{subsec:2DOFZiegler}, the choice of the numerical parameters for the 3-DOF system must be such that the divergence bifurcation is repelled away from the Hopf point. For that purpose, the parameters are selected as:
\begin{equation} \label{eq:3DOFZieglerAleParam}
    k_1 = \delta^2 k_2, \quad k_2 = \delta^2 k_3, \quad m_1 = \gamma^2 m_2, \quad m_2 = \gamma^2 m_3, \quad k_3 = m_3 = 1,
\end{equation}
with $\delta$ and $\gamma$ as given by \cref{eq:2DOFZieglerAleParam}. Once again, three distinct damping scenarios are considered: two situations with mass-proportional damping ($\xi_m = 0.01$ and $\xi_m = 0.2$), and a case with stiffness-proportional damping ($\xi_k = 0.1$). A linear stability analysis is performed, whose results are depicted in \cref{fig:3DOFZiegler_AleParam_eig}. Once again, the cases with mass-proportional damping display the coalescence of the eigenfrequencies and a divergence bifurcation, both of which are not present in the situation with stiffness-proportional damping for the variation interval of the control parameter considered here ($P\in\,[0,20]$).

\begin{figure}[H]
    \centering
    \begin{subfigure}{0.49\textwidth}
        \centering
        \includegraphics[width=\textwidth]{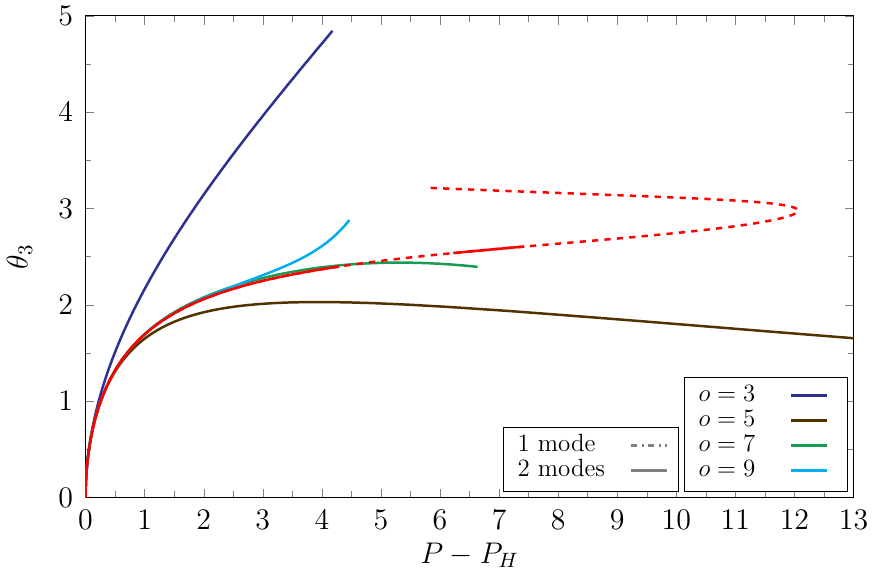}
        \caption{}
        \label{fig:3DOFZiegler_BifDiagrams_Mass_BP}
    \end{subfigure}
    \begin{subfigure}{0.49\textwidth}
        \centering
        \includegraphics[width=\textwidth]{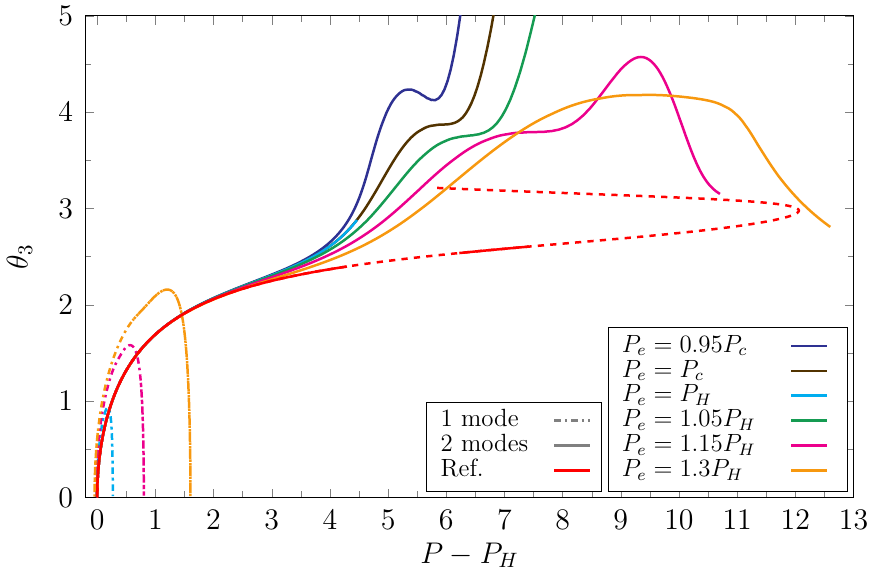}
        \caption{}
        \label{fig:3DOFZiegler_BifDiagrams_Mass_o9}
    \end{subfigure}
    \vspace{-15pt} \vphantom{a} {\scriptsize $\xi_m = 0.01, \xi_k = 0$} \vspace{6pt}

    \begin{subfigure}{0.49\textwidth}
        \centering
        \includegraphics[width=\textwidth]{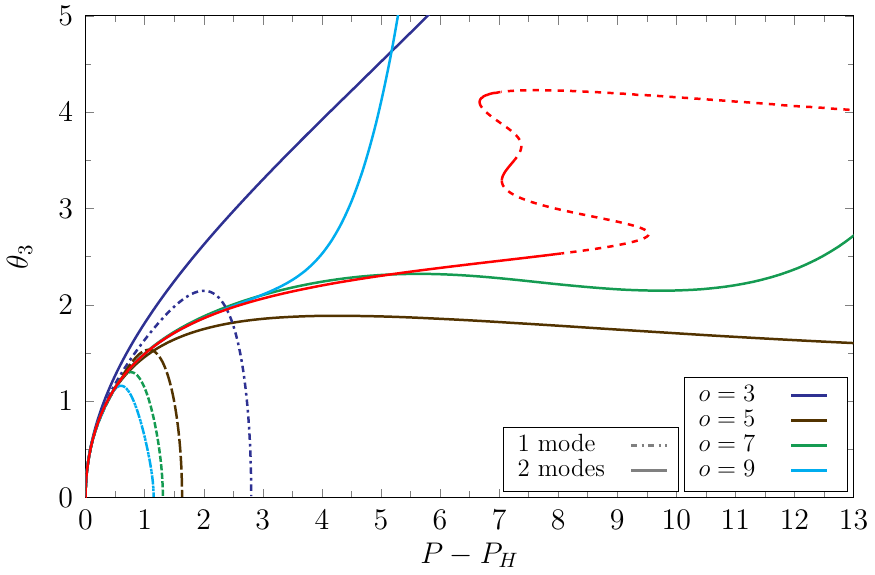}
        \caption{}
        \label{fig:3DOFZiegler_BifDiagrams_MassHigh_BP}
    \end{subfigure}
    \begin{subfigure}{0.49\textwidth}
        \centering
        \includegraphics[width=\textwidth]{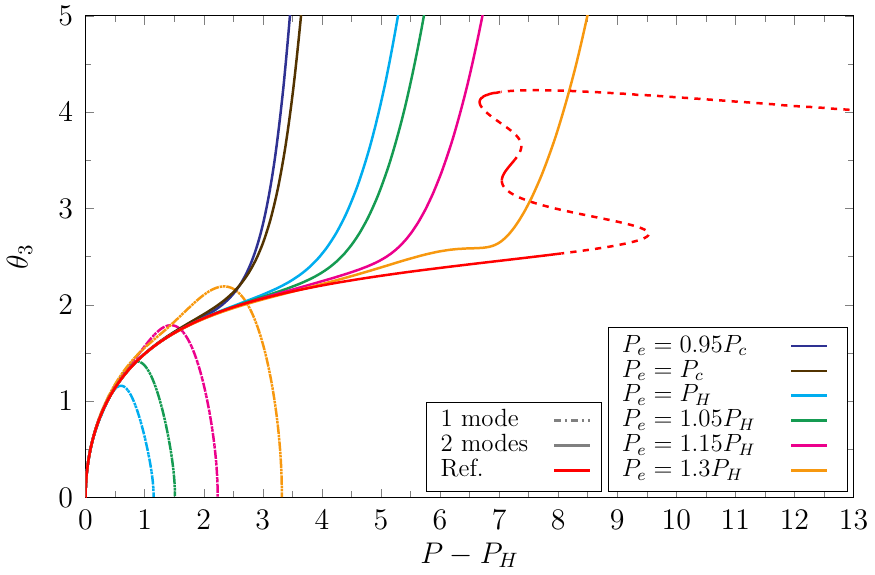}
        \caption{}
        \label{fig:3DOFZiegler_BifDiagrams_MassHigh_o9}
    \end{subfigure}
    \vspace{-15pt} \vphantom{a} {\scriptsize $\xi_m = 0.2, \xi_k = 0$} \vspace{6pt}

    \begin{subfigure}{0.49\textwidth}
        \centering
        \includegraphics[width=\textwidth]{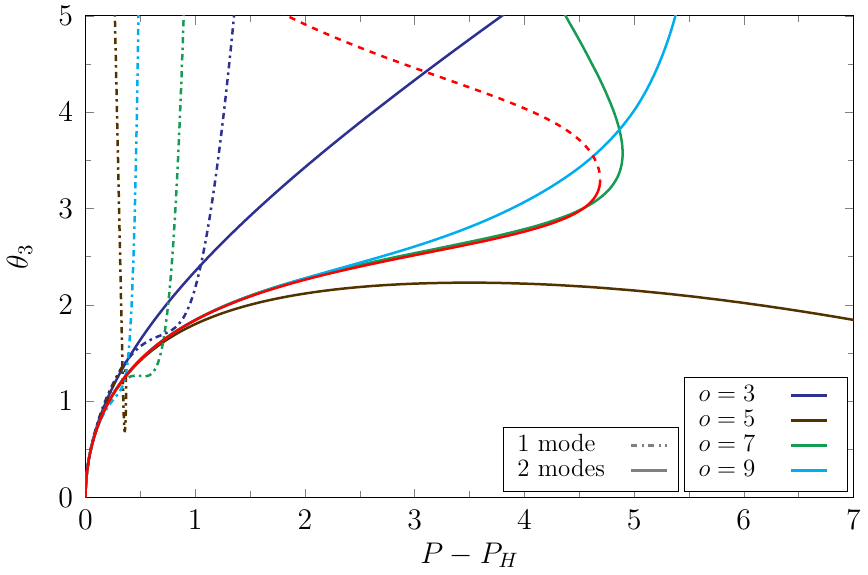}
        \caption{}
        \label{fig:3DOFZiegler_BifDiagrams_Stiff_BP}
    \end{subfigure}
    \begin{subfigure}{0.49\textwidth}
        \centering
        \includegraphics[width=\textwidth]{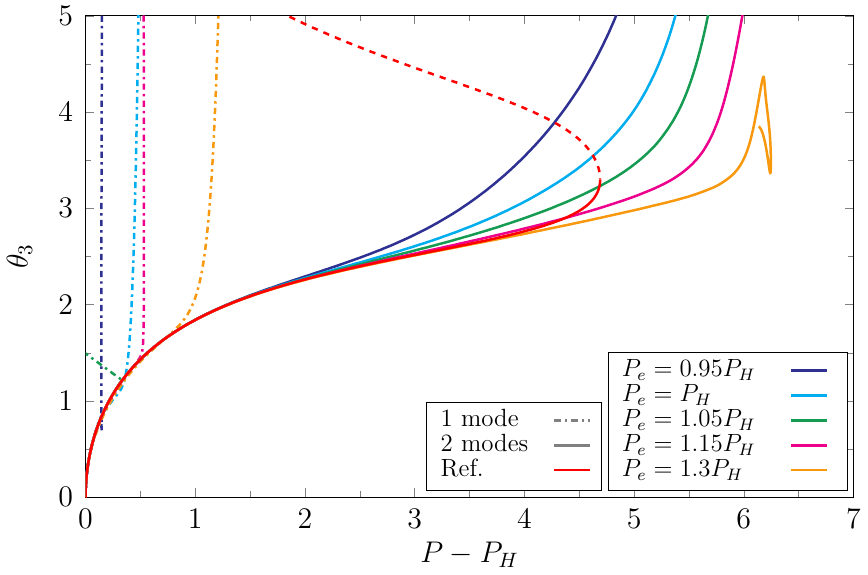}
        \caption{}
        \label{fig:3DOFZiegler_BifDiagrams_Stiff_o9}
    \end{subfigure}
    \vspace{-15pt} \vphantom{a} {\scriptsize $\xi_m = 0, \xi_k = 0.1$} \vspace{6pt}
    \caption{Bifurcation diagrams for the 3-DOF Ziegler pendulum with $L=1$ and other parameters as given by \cref{eq:3DOFZieglerAleParam}. The amplitude of the limit cycle for $\theta_3$ is given as a function of $P-P_H$. $P_e$ denotes the parametrisation point, $P_c$ the point of eigenfrequencies coalescence and $P_H$ the Hopf bifurcation point. (a), (c) and (e): solutions for increasing orders when the parametrisation is computed for $P_e=P_H$. (b), (d) and (f): solutions for a fixed order 9 of the parametrisation and different values for the expansion point $P_e$.  Each line corresponds to a different damping scenario. Full-order solutions computed by numerical continuation implemented in the package Matcont \cite{dhooge2004matcont} in red. Dashed line for unstable solutions. Stability is reported for the full-order model only.}
    \label{fig:3DOFZiegler_BifDiagrams}
\end{figure}

For the three selected cases, bifurcation diagrams exhibiting the amplitude of the limit cycle of the pendulum tip angle are presented in \cref{fig:3DOFZiegler_BifDiagrams}. They are obtained by numerical continuation of the reduced system determined by the parametrisation method, and are compared with the continuation of the full-order model, presented in red, where dashed lines depict unstable regions.

The convergence of the method relatively to the polynomial degree of the expansion can be observed in \cref{fig:3DOFZiegler_BifDiagrams_Mass_BP,fig:3DOFZiegler_BifDiagrams_MassHigh_BP,fig:3DOFZiegler_BifDiagrams_Stiff_BP}, where the parametrisation point was selected at the Hopf bifurcation point, $P_e=P_H$, and the expansion order was varied. Again, an order 9 parametrisation appears to be sufficient to retrieve the full-order results in the convergence region of the approximation.

In line with the previous case, we also remark that the two-mode strategy provides an improvement as compared to a ROM computed from the unstable mode only, in all tested cases. In some cases, see {\em e.g.}~\cref{fig:3DOFZiegler_BifDiagrams_Mass_BP}, the 1-mode ROM is even unable to locate the Hopf bifurcation point and cannot produce branches of bifurcated periodic orbits. On the other hand, the 2-mode reduction always gives satisfactory results and can predict the Hopf point.

As mentioned previously, the Jordan block strategy proposed in \cref{subsubsec:Jordan} proved to be necessary only when the expansion was performed around the point of coalescence of eigenfrequencies $P_c$. Without such adjustment, bifurcated branches are not retrieved by the ROM. For all other cases, its application does not change the results.

\cref{fig:3DOFZiegler_BifDiagrams_Mass_o9,fig:3DOFZiegler_BifDiagrams_MassHigh_o9,fig:3DOFZiegler_BifDiagrams_Stiff_o9} showcase the behaviour of the solution when the expansion point is varied and the expansion order is fixed at 9. Differently than in \cref{subsec:2DOFZiegler}, the choice of parametrisation point $P_e$ has an impact on the results even when the two-mode strategy is chosen, which can be explained by the fact that reduction is effectively performed here, since the initial problem has 3 DOF. The results are now qualitatively in line with those in \cref{subsec:2DOFZiegler} for a one-mode expansion: parametrisations performed after the bifurcation point generally produce a better approximation of the limit cycles' amplitudes on a larger range of parameter values. This phenomenon is more relevant for certain damping scenarios than others, and obviously a maximum value for the expansion point, beyond which results deteriorate, is present. In the studied situations it corresponded to approximately $P_e = 1.3P_H$. Visually, the convergence radius of the approximation seems to vary between 3 and 4 times $P_H$.

It should be noted that the full-order model bifurcation diagrams for the mass-proportional damped systems present a more complicated behaviour than those of the 2-DOF system. Specifically, Neimark-Sacker bifurcation points are present, which indicate the birth of quasi-periodic solutions. In the second case where $\xi_m=0.2$, this important change of behaviour occurs for $P-P_H\simeq 8$, a large value for which the validity limits of the ROMs are already exceeded, such that no change of stability has been found for any of the ROMs produced. For the case $\xi_m=0.01$, the Neimark-Sacker bifurcation occurs earlier, for $P-P_H\simeq 4$. In general, this bifurcation has not been retrieved by the different ROMS tested, except in a single case, but for which the location of the Neimark-Sacker point was far from the full-order solution, and has consequently not been reported in the figure for the sake of simplicity.

\section{Beck's column: a cantilever beam with a follower force}\label{sec:Beck}

\begin{figure}
    \centering
    \includegraphics[width=0.8\textwidth]{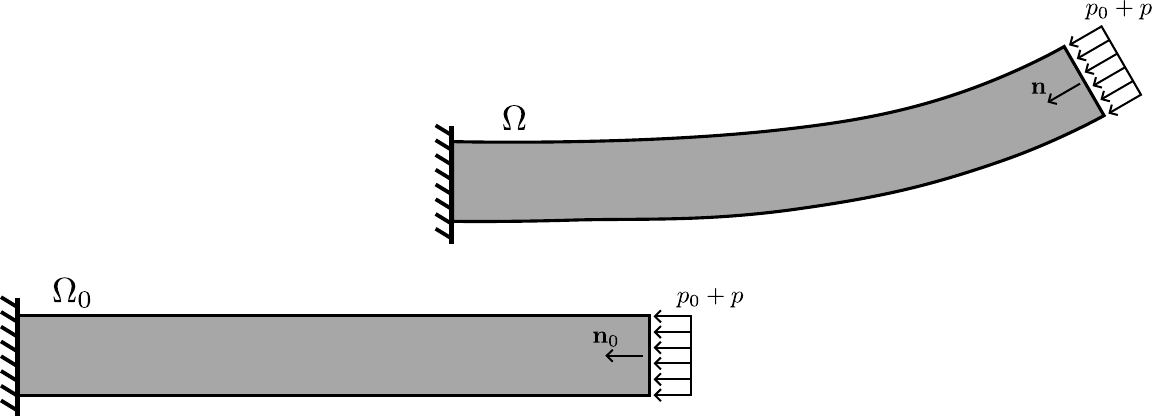}
    \caption{Beam subjected to follower force at its tip. Current and reference configurations $\Omega$ and $\Omega_0$ are shown, together with the inward normals to the section at the end of the beam, $\bfn$ and $\bfn_0$. The force is divided as the sum between a reference value $p_0$, corresponding to the expansion point for the parametrisation method, and an increment $p$, and is assumed constant along the section.}
    \label{fig:Beam_follower}
\end{figure}

In this section, the case of a continuous problem discretised by the finite element method and undergoing a Hopf bifurcation is considered, in order to highlight the application of the proposed procedure to a large dimensional finite element problem. Specifically, the cantilever beam illustrated in \cref{fig:Beam_follower}, modeled as a 2D solid in plane-strain conditions, is submitted to a non-conservative follower force at its tip. The current and reference configurations are denoted by $\Omega$ and $\Omega_0$, respectively. Additionally, the follower force is supposed to be constant throughout both the current and reference boundaries, denoted $\partial\Omega$ and $\partial\Omega_0$, and is divided into two parcels, $p_0$ and $p$, the first of which is chosen as an expansion point for the parametrisation method.

\subsection{Equations of motion and finite element discretisation}

In order to derive the equations of motion for the continuous system, the principle of virtual power is employed. It can be stated  as
\begin{equation}
 \delta \mathcal{P}_{iner} - \delta \mathcal{P}_{int} = \delta \mathcal{P}_{ext},
\end{equation} 
with $\delta \mathcal{P}_{iner}$, $\delta \mathcal{P}_{int}$ and $\delta \mathcal{P}_{ext}$ denoting the virtual powers of inertial, internal and external forces, respectively. It is assumed that the follower forces are the only external load applied to the system. Note that considering other external forces can be easily implemented with the present formalism but is not investigated here for the sake of concision. A Lagrangian approach is followed. The expressions for each of the virtual power terms are
\begin{subequations} \label{eq:virtualworks}
\begin{alignat}{2}
	&\delta \mathcal{P}_{iner} &&= \int_{\Omega_0} \tilde{\bfv} \cdot \rho_0 \ddot{\bfu} \, \dOmega_0, \label{eq:virtualworksa} \\
	&\delta \mathcal{P}_{int} &&= -\int_{\Omega_0} \left( \bfeps[\tilde{\bfv}] + \bar{\nabla}[\tilde{\bfv},\bfu] \right) : \bfS \, \dOmega_0, \label{eq:virtualworksb} \\
	&\delta \mathcal{P}_{ext} &&= \int_{\partial \Omega_0} \tilde{\bfv} \cdot \left( p_0 + p \right) \left( \bfn_0 + \frac{\bfe_3 \times \bfu_{,a}}{J_{s_0}} \right) \, \ds_0 , \label{eq:virtualworksc}
\end{alignat}
\end{subequations}
with operators $\bfeps[\tilde{\bfv}]$ and $\bar{\nabla}[\tilde{\bfv},\bfu]$ defined by
\begin{equation} \label{eq:sym_operators}
    \bfeps[\tilde{\bfv}] = \frac{1}{2} \left( \nabla \tilde{\bfv} + \nabla^T \tilde{\bfv} \right), \quad \bar{\nabla}[\tilde{\bfv},\bfu] = \frac{1}{2} \left( \nabla^T \tilde{\bfv} \cdot \nabla \bfu + \nabla^T \bfu \cdot \nabla \tilde{\bfv} \right).
\end{equation}
Additionally, $\bfu, \tilde{\bfv} \in \mathcal{C}(\bf0)$ denote the real displacement and virtual velocity fields\footnote{Note that, for simplicity, non-homogeneous Dirichlet boundary conditions are not considered even though they could be accounted for with minor technical modifications.}, $\rho_0$ the initial material density and $\bfS$ the second Piola-Kirchhoff stress tensor~\cite{Ogden1997}. In order to relate stresses and strains, we express the deformation gradient tensor as a function of displacements as
\begin{equation}
    \bfF = \bfI + \nabla \bfu,
\end{equation}
and assume a Saint Venant-Kirchhoff material, where the relation between $\bfS$ and the Green-Lagrange strain tensor
\begin{equation}
    \bfE = \frac{1}{2} \left( \bfF^T \cdot \bfF - \bfI \right) = \frac{1}{2} \left( \nabla \bfu + \nabla^T \bfu + \nabla^T \bfu \cdot \nabla \bfu \right) = \varepsilon[\bfu] + \frac{1}{2} \bar{\nabla}[\bfu,\bfu],
\end{equation}
is linear and given by
\begin{equation}
    \bfS = \fC : \bfE,
\end{equation}
with $\fC$ the fourth-order constitutive tensor, and expressions for $\varepsilon[\bfu]$ and $\bar{\nabla}[\bfu,\bfu]$ given in \cref{eq:sym_operators}.

Moreover, the operator $\nabla$ represents the gradient in the reference configuration, $a$ a curvilinear coordinate along the boundary and $J_{s_0}$ the length Jacobian. The pull-backs of \cref{eq:virtualworksa,eq:virtualworksb} from the current configuration are detailed in \cite{opreni22high,opreniPiezo}, while that of \cref{eq:virtualworksc} is given in \cref{ap:FEMdiscretization}.

Upon a finite element discretisation of the problem, the equations of motion are given by
\begin{equation}
\bfM \ddot{\bfU}_t - \left( p_0 + p \right) \bfR_0 + \left( \bfK - p_0 \bfR_u - p\bfR_u \right) \bfU_t + \bfG(\bfU_{t}, \bfU_{t}) + \bfH(\bfU_{t},\bfU_{t},\bfU_{t}) = \mathbf{0},
\end{equation} 
where $\bfM$ and $\bfK$ represent the usual finite element mass and stiffness matrices, the nonlinear quadratic and cubic tensors $\bfG$ and $\bfH$ stem from the internal forces expression, and $\bfR_0$ and $\bfR_u$ come from the follower forces term. Detailed expressions for these quantities are given in \cref{ap:FEMdiscretization}. Additionally, $\bfU_t = \bfU_0 + \bfU$ denotes the total displacement vector, with $\bfU_0$ and $\bfU$ being the displacements at the expansion point $p_0$ and its perturbation due to the parameter increment $p$. Since $\bfU_0$ is a fixed point of the system, it verifies the equilibrium equation
\begin{equation}
	\left( \bfK - p_0 \bfR_u \right) \bfU_0 + \bfG(\bfU_0,\bfU_0) + \bfH(\bfU_0,\bfU_0,\bfU_0) = p_0 \bfR_0.
\end{equation}
Eliminating the static displacements in the same way as described in \cref{subsec:Reduction}, the equations of motion can be rewritten as
\begin{equation} \label{eq:EOMfollower}
\bfM \ddot{\bfU} + \bfC \dot{\bfU} + \bfK_t \bfU - p \bfR_t + \bfG_t(\bfU, \bfU) - p \bfR_u \bfU + \bfH(\bfU,\bfU,\bfU) = \mathbf{0},
\end{equation} 
with $\bfK_t$, $\bfR_t$ and $\bfG_t$ defined by
\begin{subequations} \label{eq:KtandRt}
\begin{align}
    \bfK_t &= \bfK - p_0 \bfR_u + 2\bfG(\bfU_0,\bfI) + 3\bfH(\bfU_0,\bfU_0,\bfI), \\
    \bfR_t &= \bfR_0 - \bfR_u \bfU_0, \\
    \bfG_t(\bfU,\bfU) &= \bfG(\bfU,\bfU) + 3\bfH(\bfU_0,\bfU,\bfU),
\end{align}
\end{subequations}
and where mechanical dissipation has been added in the form of Rayleigh damping $\bfC = \alpha \bfM + \beta \bfK_t$\footnote{When comparing ROMs calculated for different expansion points, it is convenient to have a single matrix $\bfK_t$ for all of them, in order to keep a consistent damping. Therefore, in this contribution the tangent stiffness matrix calculated at the bifurcation point was chosen for all of the models.}. These equations could be further modified by transforming the system in its first-order formulation and by the addition of support variables, in order to render all nonlinearities quadratic, so that they would be in the same format as in \cref{eq:DAEaugmented}. This path is not taken, as it would unnecessarily increase the size of the system. Instead, a specific treatment considering the second-order nature of the equations, the presence of cubic nonlinearities and the fact that there are no quadratic nonlinearities on the parameter is presented in \cref{ap:2ordersimpl}.

\subsection{Numerical results}

\begin{figure}[h]
\centering
\begin{subfigure}{0.325\textwidth}
    \centering
    \includegraphics[width=\textwidth]{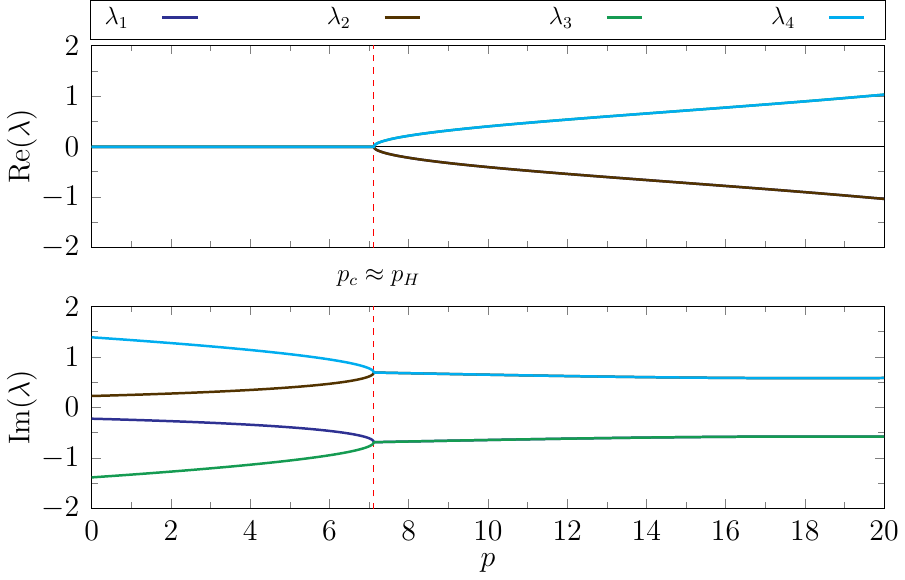}
    \caption{$\xi_m = 0.01$ and $\xi_k = 0$}
    \label{fig:Beamfollower_eig_mass}
\end{subfigure}
\begin{subfigure}{0.325\textwidth}
    \centering
    \includegraphics[width=\textwidth]{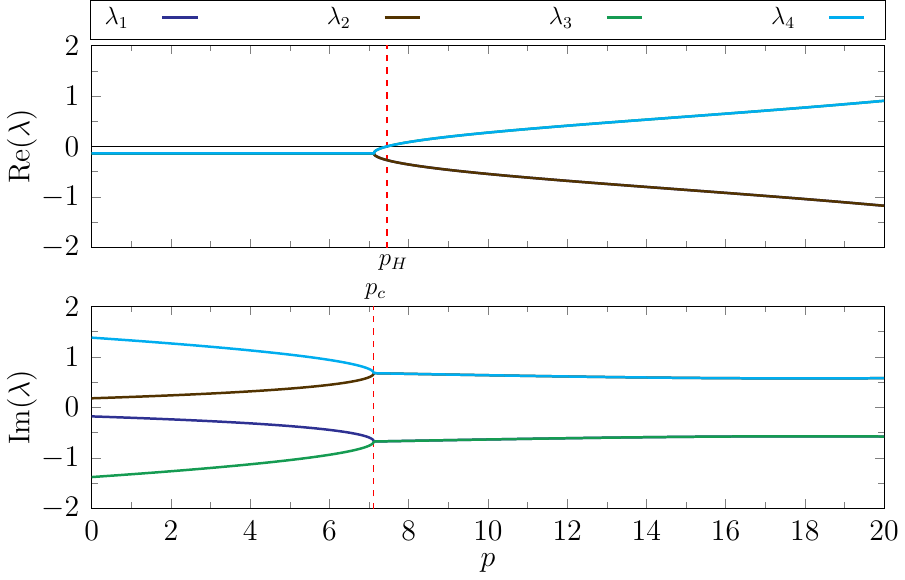}
    \caption{$\xi_m = 0.2$ and $\xi_k = 0$}
    \label{fig:Beamfollower_eig_massHigh}
\end{subfigure}
\begin{subfigure}{0.325\textwidth}
    \centering
    \includegraphics[width=\textwidth]{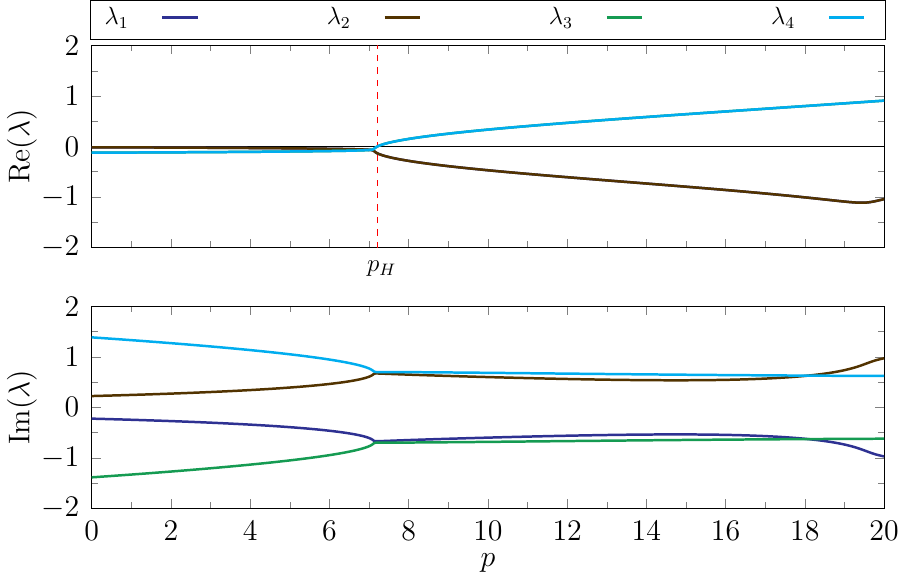}
    \caption{$\xi_m = 0$ and $\xi_k = 0.1$}
    \label{fig:Beamfollower_eig_stiff}
\end{subfigure}
\caption{Eigenvalue trajectories for Beck's column.}
\label{fig:Beamfollower_eig}
\end{figure}

The cantilever beam selected for the numerical computations has a length $L=30$ and a height $h=2$, with fictitious material properties chosen as $E=1000$, $\nu = 0$ and $\rho = 0.1$. All units are chosen to be consistent with one another. For the discretization, a mesh consisting of quadratic isoparametric triangular elements is selected, resulting in 77 elements and 153 nodes in total. The first step was to perform a stability analysis, whose results are reported in \cref{fig:Beamfollower_eig} for the following three damping scenarios: two cases with mass-proportional damping ($\xi_m = 0.01$ and $\xi_m = 0.2$), and a case with stiffness-proportional damping with $\xi_k = 0.1$. The damping ratios relate to the damping coefficients by $\alpha = \frac{\omega_0}{2}\xi_m$ and $\beta = \frac{1}{2\omega_0}\xi_k$, with $\omega_0$ a reference frequency, in this case chosen to be $\omega_0 = 0.68$, approximately the eigenfrequency at the bifurcation.

The figures underline that the main features of the stability analysis of the Ziegler pendulum are retrieved in the case of mass-proportional damping: the frequency coalescence is exactly verified, with a tendency of getting apart from the bifurcation point with increasing damping. However, the eigenvalue trajectories for stiffness-proportional damping do not resemble the ones for the Ziegler pendulum, and are much more similar to the mass-proportional cases. In this case, one can even remark on the existence of an almost coalescence just before the bifurcation point.

\begin{figure}[H]
    \raggedright
    \begin{subfigure}{0.49\textwidth}
        \centering
        \includegraphics[width=\textwidth]{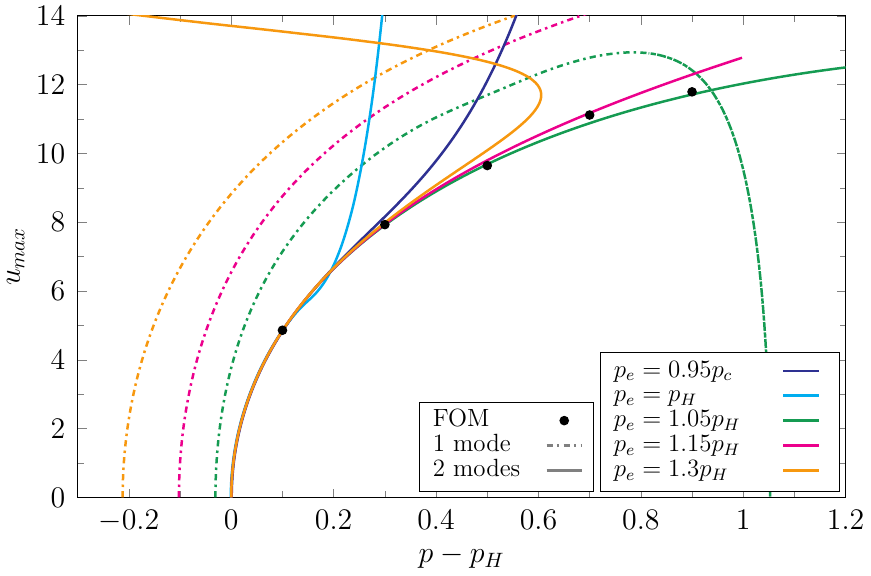}
        \caption{$\xi_m = 0.01, \xi_k = 0$}
        \label{fig:Beamfollower_BifDiagrams_Mass_o9}
    \end{subfigure}
    \begin{subfigure}{0.49\textwidth}
        \centering
        \includegraphics[width=\textwidth]{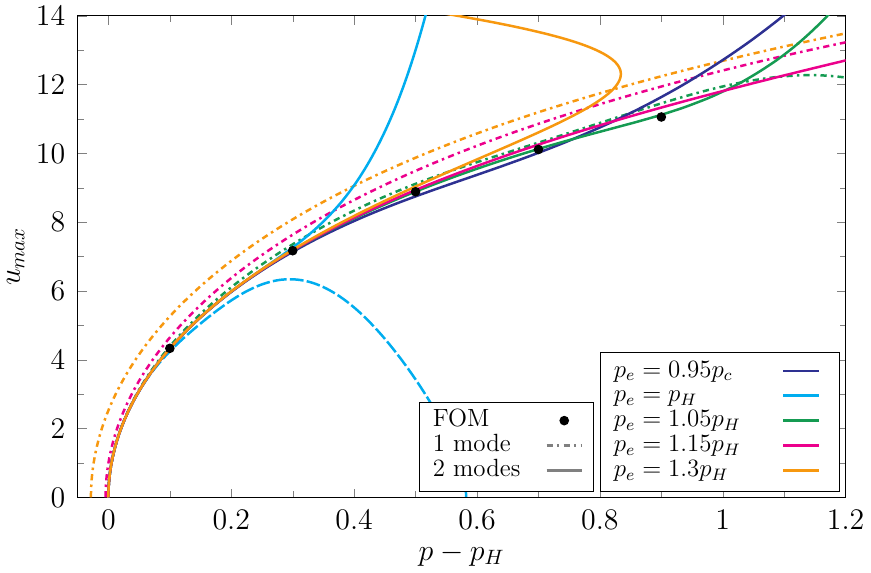}
        \caption{$\xi_m = 0.2, \xi_k = 0$}
        \label{fig:Beamfollower_BifDiagrams_MassHigh_o9}
    \end{subfigure}

    \begin{subfigure}{0.94\textwidth}
        \raggedright
        \includegraphics[width=\textwidth]{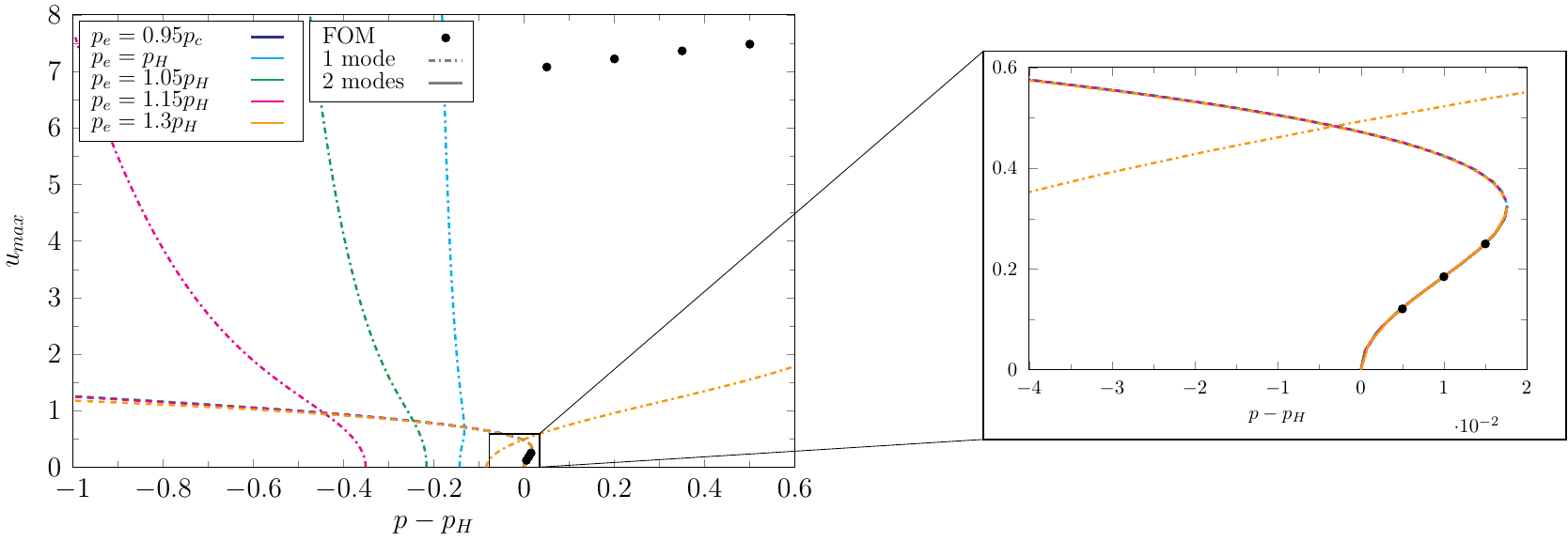}
        \caption{$\xi_m = 0, \xi_k = 0.1$}
        \label{fig:Beamfollower_BifDiagrams_Stiff_o9}
    \end{subfigure}
    \caption{Bifurcation diagrams for Beck's column. The maximum amplitude of the limit cycles of the vertical displacement of a node in the tip of the beam is given as a function of $p-p_H$, where $p_H$ denotes the Hopf bifurcation point, $p_e$ denotes the parametrisation point and $p_c$ the point of eigenfrequencies coalescence. ROMs are computed with an order 9 parametrisation and varying expansion points. Reference results obtained by direct time integration and shown with black dots. (a) mass-proportional damping with $\xi_m=0.01.$ (b) mass-proportional damping with $\xi_m=0.2.$ (c) stiffness-proportional damping with $\xi_k=0.1.$}
    \label{fig:Follower_BifDiagrams}
\end{figure}

Bifurcation diagrams showing the amplitude of limit cycle oscillations for a node at the tip of the beam are depicted in \cref{fig:Follower_BifDiagrams}. The order of the expansion in the DPIM has been set to 9 in this case, according to the previous studies, and only the effect of considering different expansion points is studied. Once again, both the one and two-mode strategies are considered, and compared to full-order model simulations obtained by direct numerical integration of the equations of motion to obtain the limit cycles' amplitudes.

The two cases with mass-proportional damping, \cref{fig:Beamfollower_BifDiagrams_Mass_o9,fig:Beamfollower_BifDiagrams_MassHigh_o9}, highlight a good agreement between the reduced order model with two master modes and the reference solution, up to a range $p - p_H \approx 1$. This corresponds to approximately 15\% of the bifurcation load, a value smaller than the corresponding one for the Ziegler pendulum. This can be explained by the higher complexity of the present system, and by the fact that this value of load corresponds to large transverse displacements, of the order of $0.4 L$, or $6 h$.
In each case, the two-mode strategy produces much more accurate results. In particular, it allows one to finely predict the location of the bifurcation point no matter the selected expansion point. In contrast, keeping only one mode in the ROM leads to an incorrect estimate of the bifurcation point when the expansion point is selected after the bifurcation, and leads to even miss the Hopf bifurcation for $p_e \leq p_H$, in \cref{fig:Beamfollower_BifDiagrams_Mass_o9}, and for $p_e \leq p_c$, in \cref{fig:Beamfollower_BifDiagrams_MassHigh_o9}, leading to the absence of a bifurcated branch in the plots. Finally, it is worth highlighting that the choice to parametrise after the bifurcation point, thus approximating the unstable manifold, provides better results with the two-mode startegy, with optimal parametrisation points located between $1.05 p_H$ and $1.15 p_H$.

The stiffness-proportional case in this situation shows an interesting behaviour, typical of a generic Hopf bifurcation scenario where the quintic terms in the normal form are dominant as compared to the cubic ones. A supercritical Hopf bifurcation is still at hand but the range over which stable limit cycles develop corresponds to a small variation of the load~$p$. For $p - p_H \simeq 0.018$, a fold bifurcation occurs giving rise to a branch of unstable limit cycles. This behaviour is correctly predicted by the ROMs computed with two master modes, for the different values of the expansion points tested and shown in Fig.~\ref{fig:Beamfollower_BifDiagrams_Stiff_o9},
with a very accurate prediction of the location of both the Hopf bifurcation point and the fold bifurcation. Interestingly, a perfect match for this very small branch of stable solutions is found with the full-order model (FOM), which has been in this case numerically integrated in time. Since no continuation procedure has been developed for the full-order model, the complete solution branch is not available. Nevertheless, stable limit cycles developing for values  $p - p_H > 0.018$ have been found with direct time integration and have a very large amplitude with $u_{max} \simeq 7$, highlighting the presence of the stable branch in this very high amplitude range. Since this branch is too far for the unstable fixed point where the parametrisation is computed, it has not been found by the different tested ROMs, that were also unable to locate the second fold bifurcation.

On the other hand, when applying the one-mode strategy, none of the expansions are able to properly capture the Hopf bifurcation point, when the expansion points are selected with negative values of $p-p_H$. Using an expansion point after the bifurcation also leads to bifurcated branches that are very badly located as compared to the full-order model, again underlining the importance of considering two master modes in the construction of the ROM.

\section{Conclusion}

This contribution is devoted to extending the applicability range of the parametrisation method for invariant manifolds to parameter-dependent systems, by focusing on the case of mechanical systems experiencing a Hopf bifurcation in the range of the parameter variation. To do so, the bifurcation parameter is included in the method's framework as an additional variable, following the so-called {\it suspension trick}~\cite{NayfehBala,Shaw94} that has already been considered in the parametrisation method in \cite{vizza2024superharm,MingwuLi2024}. The main novelty is the specific adjustements proposed in order to enhance the predictive capability of the obtained ROM. 

A first adjustment has been introduced which consists in taking into account the possible ill-conditioning of the resulting linear systems when a frequency coalescence occurs at the Hopf bifurcation point, a case generally observed for fluttering instability in conservative problems. A careful treatment of the Jordan blocks has been shown to be efficient but meaningful only in this very particular case of an exceptional point (EP). Otherwise, this strategy is not needed.

A second adjustment is to select two master modes in the ROM instead of keeping only the unstable mode. This choice has been justified by the study of the underlying conservative system, which has a center manifold of dimension four at the bifurcation point in this situation and shows either exact 1:1 resonance when there is a perfect frequency coalescence, or near 1:1 resonance when damping is added. Thus, the choice to retain the two master modes in near 1:1 resonance to construct the ROM, has been proposed. This strategy has been systematically compared to its counterpart where only one mode is retained in the parametrisation, a choice proposed in \cite{MingwuLi2024}.

In order to assess the performances, three examples were considered: Ziegler pendulums with 2- and 3-DOF and Beck's column, discretized by finite elements. Reduced order models calculated at a single expansion point in phase space were used to find bifurcation diagrams by means of numerical continuation. Different damping scenarios, and orders of expansion of the parametrisation and parametrisation points (before, at and after the bifurcation) were considered.

It was shown that the two-mode strategy is able to accurately predict the bifurcation point for all of the considered scenarios, no matter the location of the parametrisation point. This is an important finding, as it allows the use of the derived ROMs to predict bifurcations without necessarily knowing their occurrence point {\em a priori}. This is not the case for the one-mode strategy: when it is employed, if the parametrisation is performed before the bifurcation the numerical continuation is not able to find the Hopf point.

Furthermore, it has also been demonstrated that selecting two master modes for the parametrisation improves its validity range as compared to its one-mode counterpart. Considering the two Ziegler pendulums, the convergence region of the method ranged from approximately 2 to 3 times the bifurcation load. Considering Beck's column, this value is smaller, corresponding to 15\% of the bifurcation load, but is already associated with moderately large transverse displacements, of the order of $0.4L$, or $6h$.

Lastly, comparing the different choices of parametrisation points, it has been observed that parametrising after the bifurcation yields ROMs that offer the best predictions in terms of parameter variation. This result is in agreement with~\cite{MingwuLi2024}, and a phase space interpretation has been here provided. In summary, computing the ROM after the Hopf bifurcation yields correct results since the unstable manifold of the unstable fixed point connects to the stable limit cycle that is searched for. Therefore, as long as the periodic orbit is inside the convergence region of the asymptotic expansions of the parametrisation method, it can be properly reproduced by the ROM.

\section*{Acknowledgments}
Attilio Frangi acknowledges the PRIN 2022 Project ‘‘DIMIN- DIgital twins of nonlinear MIcrostructures with iNnovative model-order-reduction
strategies’’ (No. 2022XATLT2) funded by the European Union - NextGenerationEU, and Cyril Touzé acknowledges the Agence Innovation Défense (AID) who contributed to support this work through the funding
attributed to the COFLAP project (registered under the number 2023 65 0089).

\section*{Funding}
The work received no additional funding.

\section*{Conflict of interest} 
The authors declare that they have no conflict of interest.

\section*{Data availability statement}
A Julia package that implements the proposed technique is available at \url{https://github.com/MORFEproject}.

\section*{Author contributions}
\begin{itemize}
    \item Conceptualization: André de F. Stabile, Alessandra Vizzaccaro, Loïc Salles, Cyril Touzé;
    \item Methodology: André de F. Stabile, Alessandra Vizzaccaro, Alessio Colombo, Attilio Frangi;
    \item Software: André de F. Stabile, Alessio Colombo, Attilio Frangi;
    \item Formal analysis and investigation: André de F. Stabile;
    \item Writing - original draft preparation: André de F. Stabile, Cyril Touzé;
    \item Writing - review and editing: André de F. Stabile, Attilio Frangi, Cyril Touzé;
    \item Supervision: Alessandra Vizzaccaro, Loïc Salles, Attilio Frangi, Cyril Touzé.
\end{itemize}

\bibliographystyle{unsrt}
\bibliography{biblioROM}

\appendix
\gdef\thesection{\Alph{section}} 
\makeatletter
\renewcommand\@seccntformat[1]{\appendixname\ \csname the#1\endcsname.\hspace{0.5em}}
\makeatother

\section{Nonlinear tensors and their expressions in Cartesian components} \label{ap:NonlinearTensors}

In \cref{subsec:Reduction}, the identity matrix is given as an input to tensor $\bfQ_1$. To clarify the notation, the following full indicial expressions are given, where the Einstein convention of repeated indices is used:
\begin{subequations}
\begin{align}
    &\left[ \bfQ_1(\bfu,\bfv) \right]_p = (Q_1)_{pij} u_i v_j, \\
    &\left[ \bfQ_1(\bfA,\bfv) \right]_{pq} = (Q_1)_{pij} A_{iq} v_j,\\ 
    &\left[ \bfQ_1(\bfu,\bfB) \right]_{pq} = (Q_1)_{pij} u_i B_{jq}
\end{align}    
\end{subequations}
with $\bfu,\bfv \in \mathbb{C}^D$ and $\bfA,\bfB \in \mathbb{C}^{D \times D}$. It can be seen that the result of the application of $\bfQ_1$ to two vectors is a vector, while its application to a vector and a matrix gives a matrix. Analogously, for $\bfQ_2$ we define
\begin{subequations}
\begin{align}
    &\left[ \bfQ_2(\bfu,\mu) \right]_p = (Q_2)_{pi} u_i \mu, \\
    &\left[ \bfQ_2(\bfA,\mu) \right]_{pq} = (Q_2)_{pi} A_{iq} \mu,
\end{align}
\end{subequations}
with $\mu \in \mathbb{R}$. For all of the above expressions, indexes $i,j,p,q$ range from $1$ to $D$. Note in particular that these definitions allow writing
\begin{align}
    \bfQ_1(\bfy_0,\bfy) &= \bfQ_1(\bfy_0,\bfI) \bfy, \\
    \bfQ_1(\bfy,\bfy_0) &= \bfQ_1(\bfI,\bfy_0) \bfy, \\
    \bfQ_2(\bfy,\mu_0) &= \bfQ_2(\bfI,\mu_0) \bfy,
\end{align}
with $\bfI \in \mathbb{R}^{D \times D}$ the identity matrix. These expressions are used in \cref{eq:Bilinear} to pull out $\bfy$ from the argument of the tensors and sum the matrices $\bfQ_1(\bfy_0,\bfI)$, $\bfQ_1(\bfI,\bfy_0)$ and $\bfQ_2(\bfI,\mu_0)$ in order to define the tangent linear part of the dynamics, $\bfA_t$, in \cref{eq:At}.

\section{Equations of motion for the 3-DOF Ziegler pendulum} \label{ap:3DOFZieglerEOM}

This appendix is aimed at establishing the equations of motion for the 3-DOF Ziegler pendulum investigated in \cref{subsec:3DOFZiegler}, using a  Lagrangian formalism. For a non-conservative problem, the Euler-Lagrange equations write
\begin{equation} \label{eq:EulerLagrange}
    \frac{d}{dt} \left( \frac{\partial \mathcal{L}}{\partial \dot{\mathbf{q}}} \right) - \frac{\partial \mathcal{L}}{\partial \mathbf{q}} = \mathbf{Q}^{nc},
\end{equation}
with $\mathbf{q} = \thetavec = [\theta_1, \, \theta_2, \, \theta_3]$ the vector of generalized coordinates, $\mathcal{L} = T - V$ the Lagrangian, with $T$ and $V$ the kinetic and potential energies, respectively, and $\mathbf{Q}^{nc}$ the generalized non-conservative forces. The kinetic and potential energies for the system read
\begin{equation}\label{eq:Eczieg3}
    T = \frac{1}{2} \left[ m_1 \left( \dot{x}_1^2 + \dot{y}_1^2 \right) + m_2 \left( \dot{x}_2^2 + \dot{y}_2^2 \right) + m_3 \left( \dot{x}_3^2 + \dot{y}_3^2 \right) \right],
\end{equation}
and
\begin{equation}
    V = \frac{1}{2} \left[ k_1 \theta_1^2 + k_2 \left( \theta_2 - \theta_1 \right)^2 + k_3 \left( \theta_3 - \theta_2 \right)^2 \right].
\end{equation}
In~\cref{eq:Eczieg3}, cartesian coordinates are used. From the geometry shown in \cref{fig:3DOFZiegler}, and thanks to standard trigonometric relationships, one has
\begin{equation}
\begin{alignedat}{3}
    & \dot{x}_1 = L \dot{\theta}_1 \cos{\theta_1}, && \quad  \dot{x}_2 = \dot{x}_1 + L \dot{\theta}_2 \cos{\theta_2}, &&\quad \dot{x}_3 = \dot{x}_2 + L \dot{\theta}_3 \cos{\theta_3}, \\
    & \dot{y}_1 = -L \dot{\theta}_1 \sin{\theta_1}, &&\quad \dot{y}_2 = \dot{y}_1 - L \dot{\theta}_2 \sin{\theta_2}, &&\quad \dot{y}_3 = \dot{y}_2 - L \dot{\theta}_3 \sin{\theta_3}.
\end{alignedat} 
\end{equation}
The expression for the kinetic energy becomes
\begin{multline}
    T = \frac{L^2}{2} \left[ \left( m_1 + m_2 + m_3 \right) \dot{\theta}_1^2 + \left( m_2 + m_3 \right) \dot{\theta}_2^2 + m_3 \dot{\theta}_3^2 + 2\left( m_2 + m_3 \right) \dot{\theta}_1 \dot{\theta}_2 \cos{\left( \theta_1 - \theta_2 \right)} \right. \\ \left. + 2m_3 \left( \dot{\theta}_1 \dot{\theta}_3 \cos{\left( \theta_1 - \theta_3 \right)} + \dot{\theta}_2 \dot{\theta}_3 \cos{\left( \theta_2 - \theta_3 \right)} \right) \right].
\end{multline}
The different terms in \cref{eq:EulerLagrange} can be computed. The first term on the left-hand side reads

\begin{equation} \label{eq:InertialForces}
\resizebox{0.9\textwidth}{!}{
$
\begin{aligned}
    \frac{d}{dt} \left( \frac{\partial \mathcal{L}}{\partial \dot{\bftheta}} \right) &=L^2 
    \begin{bmatrix}
        m_1 + m_2 + m_3 & \left( m_2 + m_3 \right) \cos{\left( \theta_1 - \theta_2 \right)} & m_3 \cos{\left( \theta_1 - \theta_3 \right)} \\
        \left( m_2 + m_3 \right) \cos{\left( \theta_1 - \theta_2 \right)} & m_2 + m_3 & m_3 \cos{\left( \theta_2 - \theta_3 \right)} \\
        m_3 \cos{\left( \theta_1 - \theta_3 \right)} & m_3 \cos{\left( \theta_2 - \theta_3 \right)} & m_3
    \end{bmatrix}
    \ddot{
    \begin{bmatrix}
        \theta_1 \\
        \theta_2 \\
        \theta_3
    \end{bmatrix}
    } \\
    + &L^2
    \begin{bmatrix}
        0 & \left( m_2 + m_3 \right) \left( \dot{\theta}_2 - \dot{\theta}_1 \right) \sin{\left( \theta_1 - \theta_2 \right)} & m_3 \left( \dot{\theta}_3 - \dot{\theta}_1 \right) \sin{\left( \theta_1 - \theta_3 \right)} \\
        \left( m_2 + m_3 \right) \left( \dot{\theta}_2 - \dot{\theta}_1 \right) \sin{\left( \theta_1 - \theta_2 \right)} & 0 & m_3 \left( \dot{\theta}_3 - \dot{\theta}_2 \right) \sin{\left( \theta_2 - \theta_3 \right)} \\
        m_3 \left( \dot{\theta}_3 - \dot{\theta}_1 \right) \sin{\left( \theta_1 - \theta_3 \right)} & m_3 \left( \dot{\theta}_3 - \dot{\theta}_2 \right) \sin{\left( \theta_2 - \theta_3 \right)} & 0
    \end{bmatrix}
    \dot{
    \begin{bmatrix}
        \theta_1 \\[10pt]
        \theta_2 \\[10pt]
        \theta_3
    \end{bmatrix}}.
\end{aligned}
$
}
\end{equation}
In a similar fashion, the second term is found to be
\begin{equation} \label{eq:InternalForces}
\resizebox{0.9\textwidth}{!}{
$
\begin{aligned}
    \frac{\partial \mathcal{L}}{\partial \bftheta} &= - 
    \begin{bmatrix}
        k_1 + k_2 & -k_2 & 0 \\
        -k_2 & k_2 + k_3 & -k_3 \\
        0 & -k_3 & k_3
    \end{bmatrix}
    \begin{bmatrix}
        \theta_1 \\
        \theta_2 \\
        \theta_3
    \end{bmatrix} \\
    &+ L^2 
    \begin{bmatrix}
        0 & \left( m_2 + m_3 \right) \sin{\left( \theta_2 - \theta_1 \right)} & m_3 \sin{\left( \theta_3 - \theta_1 \right)} \\
        \left( m_2 + m_3 \right) \sin{\left( \theta_1 - \theta_2 \right)} & 0 & m_3 \sin{\left( \theta_3 - \theta_2 \right)} \\
        m_3 \sin{\left( \theta_1 - \theta_3 \right)} & m_3 \sin{\left( \theta_2 - \theta_3 \right)} & 0
    \end{bmatrix}
    \dot{
    \begin{bmatrix}
        \theta_1 \\
        \theta_2 \\
        \theta_3
    \end{bmatrix}
    } \odot
    \dot{
    \begin{bmatrix}
        \theta_1 \\
        \theta_2 \\
        \theta_3
    \end{bmatrix},
    }
\end{aligned}
$
}
\end{equation}
with $\odot$ denoting the Hadamard, {\em i.e. entry-wise} product. As for the generalized forces, they can found by computing
\begin{equation}
    \mathbf{Q}^{nc} = \left( \frac{\partial \mathbf{r}_3}{\partial \thetavec} \right)^T \mathbf{P},
\end{equation}
with
\begin{equation}
    \mathbf{r}_3 = L
    \begin{bmatrix}
        \sin{\theta_1} + \sin{\theta_2} + \sin{\theta_3} & \cos{\theta_1} + \cos{\theta_2} + \cos{\theta_3} & 0
    \end{bmatrix}^T
\end{equation}
the position vector of the third mass, and 
\begin{equation}
    \mathbf{P} = -P
    \begin{bmatrix}
        \sin{\theta_3} & \cos{\theta_3} & 0
    \end{bmatrix}^T
\end{equation}
the load vector. After computations, one finds
\begin{equation} \label{eq:NCForces}
    \mathbf{Q}^{nc} = PL
    \begin{bmatrix}
        \sin{\left( \theta_1 - \theta_3 \right)} \\
        \sin{\left( \theta_2 - \theta_3 \right)} \\
        0
    \end{bmatrix}.
\end{equation}

With \cref{eq:InertialForces,eq:InternalForces,eq:NCForces} at hand, it is possible to establish the equations of motion by substituting these expressions into \cref{eq:EulerLagrange}. Before that, however, in order to have a system with only polynomial nonlinearities, the sines and cosines are substituted by their series representation. Additionally, the quantities $\theta_i$, $\dot{\theta}_i$ and $\ddot{\theta}_i$, with $i = \{ 1,2,3 \}$, are assumed to be of $\mathcal{O}(\varepsilon)$, with $\varepsilon$ a small parameter, and their differences of $\mathcal{O}(\varepsilon^2)$, such that when retaining only terms up to the third order the equations of motion are found to be
\begin{multline}
    \begin{bmatrix}
        m_1 + m_2 + m_3 & m_2 + m_3 & m_3 \\
        m_2 + m_3 & m_2 + m_3 & m_3 \\
        m_3 & m_3 & m_3
    \end{bmatrix}
    \ddot{
    \begin{bmatrix}
        \theta_1 \\
        \theta_2 \\
        \theta_3
    \end{bmatrix}
    } + 
    \begin{bmatrix}
        k_1 + k_2 & -k_2 & 0 \\
        -k_2 & k_2 + k_3 & -k_3 \\
        0 & -k_3 & k_3
    \end{bmatrix}
    \begin{bmatrix}
        \theta_1 \\
        \theta_2 \\
        \theta_3
    \end{bmatrix} \\
    + PL 
    \begin{bmatrix}
        -1 & 0 & 1 \\
        0 & -1 & 1 \\
        0 & 0 & 0
    \end{bmatrix}
    \begin{bmatrix}
        \theta_1 \\
        \theta_2 \\
        \theta_3
    \end{bmatrix}
    = -\frac{PL}{6}
    \begin{bmatrix}
        \left( \theta_1 - \theta_3 \right)^3 \\
        \left( \theta_2 - \theta_3 \right)^3 \\
        0
    \end{bmatrix},
\end{multline}
where it is possible to identify the general form of \cref{eq:EOMZiegler2DOF} and the expressions given in \cref{eq:Matrices3DOFZiegler}.

\section{Inclusion of Jordan blocks: calculation details} \label{app:camillejordan}

In this appendix a numerical procedure to determine a properly conditioned eigenspace associated to the coalescing eigenvalues of an exceptional point is described. The main objective is to guarantee that the two eigenvectors corresponding to this block, obtained by standard numerical routines, are not almost aligned, and thus properly generate a two-dimensional eigenspace, without the risk of incurring in numerical conditioning problems. To do that, a Jordan block will be artificially imposed in the eigenvalues matrix, and new eigenvectors spanning the desired eigenspace will be determined. For this, consider the original eigenvalue problem that is solved numerically:
\begin{subequations} \label{eq:eigenproblem}
\begin{align}
    \bfB \bfY \bfD &= \bfA_t \bfD, \\
    \bfD \bfX^* \bfB &= \bfX^* \bfA_t,
\end{align}
\end{subequations}
with $^*$ denoting the Hermitian transpose, and where, following the discussion at the end of the introduction of \cref{subsec:Reduction}, only the first $d$ eigenvalues are considered. The triplet $(\bfD,\bfY,\bfX)$ found by solving the eigenproblem consists of the eigenvalues matrix $\bfD$, with size $d \times d$, and of the left and right eigenvector matrices $\bfX$ and $\bfY$ with size $D \times d$. They are obtained by employing a traditional numerical eigenvalue solver to the problem. Since numerically two eigenvalues are never exactly the same, matrix $\bfD$ will be diagonal, but such that, in the presence of an exceptional point, $D_{ii} \approx D_{jj}$ for the two indices $i$ and $j$ related to the degenerate eigenvalues, with $i<j$, and the eigenvectors corresponding to them will be nearly aligned. Our goal is then to find a triplet $(\bfLambda,\bfYbar,\bfXbar)$ that also satisfies the eigenproblem, but to impose the desired Jordan block at the new eigenvalues matrix $\bfLambda$, such that it reads
\begin{subequations}
\begin{gather}
    \text{diag}(\bfLambda) = \text{diag}(\bfD) \\
    \text{if }D_{ii} \approx D_{jj} \text{ and } i<j: \bfLambda_{ij} = \tau_{ij},
\end{gather}
\end{subequations}
which, for the block relating to the indexes $i$ and $j$, yields
\begin{equation}
    \bfLambda_{
    \tiny
    \begin{bmatrix}
        ii & ij \\ ji & jj
    \end{bmatrix}
    }
    =
    \begin{bmatrix}
        \lambda_i & \tau_{ij} \\
        0 & \lambda_j
    \end{bmatrix},
\end{equation}
with $\tau_{ij}$ not assumed equal to 1 for the sake of generality. Since $(\bfLambda,\bfYbar,\bfXbar)$ solves the eigenproblem given by \cref{eq:eigenproblem}, it follows that
\begin{equation}
    \label{eq:reduction}
    \bfLambda \bfXbar^* \bfB \bfYbar = \bfXbar^* \bfA_t \bfYbar =
    \bfXbar^* \bfB \bfYbar \bfLambda,
\end{equation}
and since $\bfLambda$ is imposed to be non-diagonal, the equality of the leftmost hand side and the rightmost hand side imposes $\bfXbar^* \bfB \bfYbar$ and $\bfXbar^* \bfA_t \bfYbar$ to be diagonal. In particular, one can choose
\begin{subequations}
\begin{align}
    \bfXbar^* \bfB \bfYbar  &= \bfI, \\
    \bfXbar^* \bfA_t \bfYbar  &= \bfLambda.
\end{align}    
\end{subequations}
    
In order not to change the span of the subspaces generated by $\bfY$ and $\bfX$, $\bfYbar$ and $\bfXbar$ are sought as linear combinations of them, namely $\bfYbar = \bfY \bfmu$ and $\bfXbar = \bfX \bfnu$, with $\bfmu$ and $\bfnu$ two $d \times d$ matrices. By inputting these definitions into \cref{eq:reduction}, it is possible to see that they have to respect 
\begin{subequations}
\begin{align}
    & \bfmu\bfLambda = \bfD\bfmu \label{eq:mu}, \\
    & \bfnu = ((\bfX^* \bfB \bfY \bfmu)^{*})^{-1}.
\end{align}  
\end{subequations}

One notices that the definition of $\bfnu$ stems directly from that of $\bfmu$, and that $\bfmu$ can be defined with an arbitrary amplitude since \cref{eq:mu} has the structure of an eigenvalue problem. For each non-degenerate eigenvalue, say the $l$-th, \cref{eq:mu} is solved by imposing the $l$-th column of $\bfmu$ to be equal to the unit vector $\bfe_l$. However, for the pair of two degenerate eigenvalues $i$ and $j$
\begin{equation}
    \begin{bmatrix}
        \mu_{ii} & \mu_{ij} \\
        \mu_{ji} & \mu_{jj} 
    \end{bmatrix}
    \begin{bmatrix}
        \lambda_i & \tau_{ij} \\
        0 & \lambda_j
    \end{bmatrix}
    =
    \begin{bmatrix}
        \lambda_i & 0 \\
        0 & \lambda_j
    \end{bmatrix}
    \begin{bmatrix}
        \mu_{ii} & \mu_{ij} \\
        \mu_{ji} & \mu_{jj} 
    \end{bmatrix},
\end{equation}
which is solved by
\begin{align}
    &
    \mu_{ij}= 
    \dfrac{\tau_{ij}}{\lambda_i-\lambda_j}\mu_{ii},
    \\
    &
    \mu_{ji}=0,
\end{align}
where $\tau_{ij}$, $\mu_{ii}$, and $\mu_{jj}$ are still arbitrary. In order to fix a choice, we impose $\mu_{ii}=1$ so that the $i$-th eigenvector stays in the basis unchanged: $\bfYbar_i = \bfY_i$. The value of entry $\mu_{jj}$ remains to be determined, and is done by imposing a condition ensuring that the new eigenvector, $\bfYbar_j$, has as norm similar to the other ones, as follows. At present, the expressions for the eigenvectors are given by
\begin{equation}
    \bfYbar_i = \bfY_i, \quad \bfYbar_j = \dfrac{\tau_{ij}}{\lambda_i-\lambda_j} \bfY_i +
    \mu_{jj} \bfY_j, \quad \bfXbar_i = \bfX_i - \dfrac{\tau_{ij}}{\lambda_i-\lambda_j}\dfrac{1}{\mu_{jj}} \bfX_j, \quad \bfXbar_j = \dfrac{1}{\mu_{jj}} \bfX_j.
\end{equation}
One should then notice that $\bfY_j$ can always be expressed as
\begin{equation}
    \bfY_j = \gamma_{ij}\bfY_i - \Delta\bfY_{ij},
\end{equation}
with $\Delta\bfY_{ij}$ orthogonal to $\bfY_i$, and $\gamma_{ij}$ a normalization constant that can be found from
\begin{equation}
    \gamma_{ij} = \dfrac{\bfY_i^* \bfY_j}{\bfY_i^* \bfY_i}.
\end{equation}
Now, since the difference between $\lambda_i$ and $\lambda_j$ is small and the same can be said of the norm of $\Delta\bfY_{ij}$, as $\bfY_i$ and $\bfY_j$ are almost perfectly aligned, the new eigenvector $\bfYbar_j$ can be chosen as $\bfYbar_j \propto \dfrac{\Delta\bfY_{ij}}{\lambda_i-\lambda_j}$, so that the ratio of two small quantities should give a finite valued vector orthogonal to $\bfYbar_i$. This is achieved by choosing 
\begin{equation}
    \mu_{jj} = -\dfrac{\tau_{ij}}{\gamma_{ij}(\lambda_i-\lambda_j)},
\end{equation}
such that the expression for the new eigenvector is finally found as
\begin{equation}
    \bfYbar_j = \dfrac{\tau_{ij}}{\lambda_i-\lambda_j}
    \left(
    \bfY_i - \frac{\bfY_j}{\gamma_{ij}}
    \right).
\end{equation}

The parameter $\tau_{ij}$ serves as a scaling factor, and will affect the norm of the new eigenvectors. It can then be chosen to maintain all of them with lengths of the same order of magnitude. However, an automatic procedure to do so is not simple, and its optimal choice might often depend on previous knowledge of the problem at hand. In this contribution, we always choose $\tau_{ij}=1$, as is the standard for the definition of Jordan blocks. This value was verified not to cause numerical instabilities for the problems at hand.

\paragraph{Additional terms in the homological equations}\mbox{}

The introduction of off-diagonal entries due to the presence of Jordan blocks in the matrix of eigenvalues will create new terms in the homological equations that are solved recursively in the parametrisation method. This development aims at focusing on the new terms as compared to the regular one, already treated and developed in~\cite{vizza2024superharm}. For a complete understanding of all the terms involved in the process, the interested reader is thus referred to~\cite{vizza2024superharm}. Herein, only the new terms are detailed.  At order $p$, the homological equation stemming from \cref{eq:DAEaugmented} reads:
\begin{equation}
    \bfB \P{\nabla_{\bfztil} \bfW(\bfztil) \bff(\bfztil)} = \bfA_t \P{\bfW(\bfztil)} + \P{\bfQ(\bfW(\bfztil),\bfW(\bfztil))}.
\end{equation}
In this equation, the only term involving the reduced dynamics is $\P{\nabla_{\bfztil} \bfW(\bfztil) \bff(\bfztil)}$, being therefore the sole affected by the off-diagonal terms. Its expression can be rewritten as
\begin{equation}
\footnotesize
    \P{\nabla_{\bfztil} \bfW(\bfztil) \bff(\bfztil)} = \sum_{s=1}^{d+1} \left[ \bfW^{(1,s)} \P{f_s(\bfztil)} \vphantom{\left( \sum_{j=1}^{d+1} f_s^{(1,j)} \widetilde{z}_j \right)} + \left( \sum_{j=1}^{d+1} f_s^{(i,j)} \widetilde{z}_j \right) \frac{\partial \P{\bfW(\bfztil)}}{\partial \widetilde{z}_s} + \P{\frac{\partial \P[\dblnk]{\bfW(\bfztil)}}{\partial \widetilde{z}_s} \P[\dblnk]{f_s(\bfztil)}} \right].
\end{equation}
Out of the three parcels in the sum, which will be called respectively $\bfN_1(\bfztil)$, $\bfN_2(\bfztil)$ and $\bfN_3(\bfztil)$, following the same notation as in \cite{vizza2024superharm}, only $\bfN_2(\bfztil)$ is modified, such that the other terms are not further expanded here. It will be responsible for new terms in the homological equations, that can be calculated as
\begin{subequations} \label{eq:N2Jordan}
\begin{gather}
    \bfNtil_2(\bfztil) = \sum_{s=1}^{d} \sum_{j=s+1}^{d+1} \sum_{k_W=1}^{m_p} \alpha_s(p,k_W) f_s^{(1,j)} \bfW^{(p,k_W)} \bfztil^{\alphavec(p,k_W)-\bfe_s+\bfe_j}, \\
    \alphavec(p,k) = \alphavec(p,k_W) -\bfe_s+\bfe_j.
\end{gather}
\end{subequations}
When solving the homological equation for the monomial numbered $\boldsymbol{\alpha}(p,k)$, the nonlinear mapping coefficient $\mathbf{W}^{(p,k_W)}$ is involved, and should be known at this stage of the solution process provided that a suitable ordering of the monomials is chosen, due to the upper-triangular structure of the $\bff^{(1)}$ matrix. $\bfNtil_2(\bfztil)$ thus contributes to the RHS of the equations, and can be simply summed to the expression of $\bfN_2(\bfztil)$ given in \cite{vizza2024superharm}.

\section{Weak from of the power of external forces in the reference configuration and finite element discretization} 
\label{ap:FEMdiscretization}

In this appendix, a derivation of the pull-back of the virtual power of the external (follower) forces, \cref{eq:virtualworksc}, from the current to the reference configuration is detailed. This renders explicit geometric nonlinearities concealed in the fact that in the former case the integrals are evaluated in the current, evolving configuration of the body. In the end, explicit expressions for the matrices related to this term on the finite element equations of motion are given.

In current configuration, the virtual power of the external forces is given by
\begin{equation}
    \delta \mathcal{P}_{ext} = \int_{\partial \Omega} \tilde{\bfv} \cdot \left( p_0 + p \right) \bfn \, \ds,
\end{equation}
with $\bfn$ an inward unit normal vector in the current configuration, as depicted in \cref{fig:Beam_follower}. In order to simplify calculations, a curvilinear coordinate $a$, attached to the boundary where the follower force is applied\footnote{At this level, $a$ can be whatever desired curvilenear coordinate, but in what follows it will be the one stemming from the finite element discretization without a change of notation.}, is introduced, such that
\begin{equation}
    \ds = J_s \da, \quad \ds_0 = J_{s_0} \da,
\end{equation}
with $J_s = \norm{\bfy_{,a}}$ and $J_{s_0} = \norm{\bfx_{,a}}$, where $\bfy$ and $\bfx$ denote the current and reference positions, respectively. Then, the vector
\begin{equation}
    \bft = \bfy_{,a} = \bfx_{,a} + \bfu_{,a},
\end{equation}
is tangent to the boundary, whose positive orientation is assumed to be counterclockwise. It can be rotated of 90° counterclockwise and normalized in order to find
\begin{equation}
    \bfn = \frac{\bfe_3 \times \bfy_{,a}}{\norm{\bfy_{,a}}} = \frac{1}{J_s} \bfe_3 \times \left( \bfx_{,a} + \bfu_{,a} \right).
\end{equation}
Thus 
\begin{equation}
    \bfn \, \ds = \frac{1}{J_s} \bfe_3 \times \left( \bfx_{,a} + \bfu_{,a} \right) J_s \, \da = \frac{1}{J_{s_0}} \bfe_3 \times \left( \bfx_{,a} + \bfu_{,a} \right) \, \ds_0 = \left( \bfn_0 + \frac{\bfe_3 \times \bfu_{,a}}{J_{s_0}} \right) \, \ds_0,
\end{equation}
where $\bfn_0 = \frac{\bfe_3 \times \bfx_{,a}}{J_{s_0}}$ is the inward unit normal vector in the reference configuration. With this, the expression for the virtual power contribution finally becomes
\begin{equation} \label{eq:weak_form_ext}
\begin{aligned}
    \delta \mathcal{P}_{ext} &= \left( p + p_0 \right) \int_{\partial \Omega_0} \tilde{\bfv} \cdot \bfn_0 \, \ds_0 + \left(p + p_0 \right) \int_{\partial \Omega_0} \tilde{\bfv} \cdot \frac{\bfe_3 \times \bfu_{,a}}{J_{s_0}} \, \ds_0,
\end{aligned}
\end{equation}
where the two integrals above will contribute to $\bfR_0$ and $\bfR_u$ in \cref{eq:EOMfollower,eq:KtandRt}, respectively.

\subsection*{Finite element discretization}
In this section, the finite element discretisation of the power of external forces will be performed. In order to do so, the real displacement and virtual velocity fields are written as a function of the nodal parameters by using shape functions matrix $\bfN$:
\begin{equation}
    \bfu = \bfN \bfU, \quad  \tilde{\bfv} = \bfN \tilde{\bfV},
\end{equation}
with $\bfU$ and $\tilde{\bfV}$ denoting the nodal displacements and virtual velocities for a numbering of degrees of freedom such that directions $x$ and $y$ alternate for increasing node numbers. It should be noted that, since the integrals to be computed are defined only on the boundary of the domain, the elements that will be considered are line elements on the boundary, such that the shape functions above depend only on a parametric coordinate $a$ associated to those elements. The derivatives of the displacements and positions (the former appearing implicitly in the vector $\bfn_0$) in \cref{eq:weak_form_ext} can thus be written as
\begin{align}
    \bfu_{,a} &= \bfN_{,a} \bfU, \\
    \bfx_{,a} &= \bfN_{,a} \bfX.
\end{align}

Expressions for the matrices $\bfN$ and $\bfN_{,a}$, are explicitly defined as
\begin{gather}
    \bfN =
    \begin{bmatrix}
        N_{1} & 0 & \cdots & N_{n} & 0 \\
        0 & N_{1} & \cdots & 0 & N_{n}
    \end{bmatrix}, \quad
    \bfN_{,a} =
    \begin{bmatrix}
        N_{1,a} & 0 & \cdots & N_{n,a} & 0 \\
        0 & N_{1,a} & \cdots & 0 & N_{n,a}
    \end{bmatrix}, 
\end{gather}
with $n$ the number of nodes per line element. Moreover, vector $\bfX$ collects the nodal coordinates in the reference configuration. The cross products on \cref{eq:weak_form_ext} can be represented with the help of matrix $\bfE_3 = \begin{bmatrix}
    0 & -1 \\
    1 & \phantom{-}0
\end{bmatrix}$, such that
\begin{align}
    \bfe_3 \times \bfu_{,a} &= \bfE_3 \bfN_{,a} \bfU, \\
    \bfn_0 &= \frac{1}{J_{s_0}} \bfE_3 \bfN_{,a} \bfX.
\end{align}

With these definitions, the integrals in \cref{eq:weak_form_ext} can be computed for a generic line finite element $e$ as
\begin{align}
    \int_{\partial \Omega_0} \tilde{\bfv} \cdot \bfn_0 \, \ds_0 &= \tilde{\bfV}^T \int_{\hat{\Omega}_e} \bfN^T \bfE_3 \bfN_{,a} \bfX \, \da = \tilde{\bfV}^T \bfR_0^e, \\
    \int_{\partial \Omega_0} \frac{\bfe_3 \times \bfu_{,a}}{J_{s_0}} \, \ds_0 &= \tilde{\bfV}^T \left( \int_{\hat{\Omega}_e} \bfN^T \bfE_3 \bfN_{,a} \da
    \right)\bfU = \tilde{\bfV}^T \bfR_u^e \bfU,
\end{align}
where now the integration is performed over a reference line element of domain $\hat{\Omega}_e$, and in practice is computed via numerical integration. The quantities $\bfR_0^e$ and $\bfR_u^e$ are the elemental level counterparts of $\bfR_0$ and $\bfR_u$ in \cref{eq:EOMfollower}, and can be obtained through the standard finite element procedure of assembly, see {\em e.g.} \cite{Belytschko2013,Bonnet2014}.

\section{Simplification for second-order mechanical systems} \label{ap:2ordersimpl}

In order to transform the equations of motion stemming from vibrating structures from second to first-order in time to fit the framework of \cref{subsec:Reduction}, the velocities of each node are introduced as auxiliary variables to the problem, such that a set of additional equations is appended to the system, doubling its size and increasing the time needed to build the reduced-order model. These equations, however, are trivial, and this fact can be taken into account by substituting the velocity mappings as a function of the displacement ones, halving the size of the system to be solved. Additionally, even though the treatment of the equations proposed in \cref{subsec:Reduction} is general, as all analytic nonlinearities can be transformed into quadratic ones, it is once again not the most computationally efficient one in the present situation, as the introduction of auxiliary variables is needed in order to apply quadratic recast to the cubic nonlinearities present in \cref{eq:EOMfollower}. This can be quite expensive for high-dimensional finite element systems. Therefore, this section is concerned with the specialization of the DPIM algorithm to mechanical systems of the form given by \cref{eq:EOMfollower}. Specifically, quadratic terms $\bfQ_3(\mu,\mu)$ from \cref{eq:DAEaugmented} will be neglected, since this corresponds to all of the situations treated in this paper. Departing from \cref{eq:EOMfollower}, it can be put into first-order format as
\begin{equation}
\small
    \begin{bmatrix}
        \bfM & \bf0 & \bf0 \\
        \bf0 & \bfM & \bf0 \\
        \bf0 & \bf0 & 1
    \end{bmatrix}
    \dot{
    \begin{bmatrix}
        \bfU \\
        \bfV \\
        p
    \end{bmatrix}
    }
    =
    \begin{bmatrix}
        \phantom{-}\bf0 & \phantom{-}\bfM & \bf0 \\
        -\bfK_t & -\bfC & \bfR_t \\
        \phantom{-}\bf0 & \phantom{-}\bf0 & 0
    \end{bmatrix}
    \begin{bmatrix}
        \bfU \\
        \bfV \\
        p
    \end{bmatrix}
    +
    \begin{bmatrix}
        \bf0 \\
        \bfG_t(\bfU,\bfU) - p\bfR_u\bfU \\
        0
    \end{bmatrix}
    +
    \begin{bmatrix}
        \bf0 \\
        \bfH(\bfU,\bfU,\bfU)\\
        0
    \end{bmatrix},
\end{equation}
where the first equation connects velocities $\bfV$ and displacements $\bfU$. In what follows, the nonlinear mappings for the unknowns are expanded as
\begin{subequations}
\begin{align}
    \bfU &= \bfU(\bfz) = \sum_{p=1}^{o} \P{\bfU(\bfz)} = \sum_{p=1}^{o} \sum_{k=1}^{m_p} \bfU^{(p,k)} \bfz^{\alphavec (p,k)},\\
    \bfV &= \bfV(\bfz) = \sum_{p=1}^{o} \P{\bfV(\bfz)} = \sum_{p=1}^{o} \sum_{k=1}^{m_p} \bfV^{(p,k)} \bfz^{\alphavec (p,k)}.
\end{align}
\end{subequations}

\subsection*{Eigenproperties of the second-order system}

In order to particularize the eigenvalue problem to the second-order scenario, the left and right eigenvectors are divided into their parts corresponding to displacements, velocities and to the pressure (external bifurcation parameter), as
\begin{equation}
    \bfY_s = 
    \begin{bmatrix}
        \bfY_s^U \\
        \bfY_s^V \\
        \text{Y}_s^p
    \end{bmatrix}, \quad
    \bfX_s = 
    \begin{bmatrix}
        \bfX_s^U \\
        \bfX_s^V \\
        \text{X}_s^p
    \end{bmatrix}, \quad
    s = 1,\ldots,2d+1.
\end{equation}

When $s \in [1,\ldots,2d ]$, by employing \cref{eq:usualEigenvalues1storder} the usual expressions for the left and right eigenproblems stemming from a second-order mechanical system are recovered:
\begin{subequations}
\begin{align}
    \bfY_s^V &= \lambda_s \bfY_s^U \\
    \left( \lambda_s^2 \bfM + \lambda_s \bfC + \bfK_t \right) \bfY_s^U &= \bf0,
\end{align}
\end{subequations}
and
\begin{subequations}
\begin{align}
    (\bfX_s^U)^* \bfM &= (\bfX_s^V)^* \left( \lambda_s \bfM + \bfC \right) \\
    (\bfX_s^V)^* \bfK_t &= -\lambda_s (\bfX_s^U)^* \bfM,
\end{align}
\end{subequations}
with $\text{Y}^p_s = 0$ and $\text{X}^p_s = 0$. These equations reflect the fact that the velocities are trivially related to the displacements as their derivative in time, and thus impose simple relationships between the parts of the eigenvalues corresponding to these two quantities, as already observed {\em e.g.} in ~\cite{vizza21high,vizza2024superharm}. Additionally, when considering the part of the eigenproblem associated with the bifurcation parameter $p$, such that $s=2d+1$, \cref{eq:pEigenvalue1storder} can be specialized to this case, yielding
\begin{equation}
    \begin{bmatrix}
        \bf0 & \bfM \\
        -\bfK_t & -\bfC
    \end{bmatrix}
    \begin{bmatrix}
        \bfY_{2N+1}^U \\
        \bfY_{2N+1}^V 
    \end{bmatrix}
    = 
    \begin{bmatrix}
        \bf0 \\
        -\bfR_t
    \end{bmatrix}
    \quad \Rightarrow \quad
    \bfY_{2N+1}^U = \bfK_t^{-1} \bfR_t, \quad \bfY_{2N+1}^V = \bf0,
\end{equation}
and $\text{Y}^p_{2d+1} = 1$. This equation reveals a clear physical interpretation for the eigenvalue related to the control parameter. It represents a static deformation under the scalar load $p$ for the linearized problem about the equilibrium state corresponding to $p_0$. 

\subsection*{Order 1 parametrisation}
The order 1 homological equation for this case writes
\begin{equation} \label{eq:Homological1}
    \begin{bmatrix}
        \bfB & \bf0 \\
        \bf0 & 1 \\
    \end{bmatrix}
    \begin{bmatrix}
        \bar{\bfW}^{(1)} \bar{\bff}^{(1)} & \bar{\bfW}^{(1)} \bff^{(1,d+1)} \\
        \bf0 & 0
    \end{bmatrix}
    \begin{bmatrix}
        \bfz \\
        p
    \end{bmatrix}
    =
    \begin{bmatrix}
        \bfA_t & \bar{\bfR}_t \\
        \bf0 & 0 \\
    \end{bmatrix}
    \begin{bmatrix}
        \bar{\bfW}^{(1)} & \bfW^{(1,d+1)} \\
        \bf0 & 1
    \end{bmatrix}
    \begin{bmatrix}
        \bfz \\
        p
    \end{bmatrix},
\end{equation}
with the linear part of the nonlinear mappings divided as
\begin{equation}
    \left[ \bfW(\bfztil) \right] = \bfW^{(1)} \bfztil =
    \begin{bmatrix}
        \bar{\bfW}^{(1)} & \bfW^{(1,d+1)}
    \end{bmatrix}
    \begin{bmatrix}
        \bfz \\
        p
    \end{bmatrix},
\end{equation}
and where the augmented (with the addition of the velocity variables) external forces vector $\bar{\bfR}_t = \begin{bmatrix} \bf0 & \bfR_t \end{bmatrix}^T$. It should be noted that the last line of \cref{eq:Homological1} is a tautology. This is a consequence of the fact that the bifurcation parameter is excluded from the nonlinear mappings and that its reduced dynamics is known, which in turn makes the last equation redundant. Splitting now the first line into its parts corresponding to usual normal variables and to the parameter, we obtain for the part related to $\bfz$:
\begin{equation}
    \bfB \bar{\bfW}^{(1)} \bar{\bff}^{(1)} = \bfA_t \bar{\bfW}^{(1)},
\end{equation}
where we recognize the usual structure of the right eigenvalue problem, such that
\begin{subequations}
\begin{align}
    \bar{\bfW}^{(1)} &= \bfY \\
    \bar{\bff}^{(1)} &= \bfD.
\end{align}
\end{subequations}

In the above expressions, it is assumed that no Jordan blocks are present in the linear parts of the dynamics. When this is the case, one simply has to impose $\bar{\bfW}^{(1)} = \bfYbar$ and $\bar{\bff}^{(1)} = \bfLambda$. Considering now the part of \cref{eq:Homological1} related to the parameter $p$, we have
\begin{equation} \label{eq:Homological1p}
    \bfA_t \bfW^{(1,d+1)} = \bfB \bfY \bff^{(1,d+1)} - \bar{\bfR}_t.
\end{equation}
This equation is underdetermined, since both $\bfW^{(1,d+1)}$ and $\bff^{(1,d+1)}$ are unknowns of the problem, and thus several different choices are possible for the values of these parameters. Specifically, if one chooses to employ the normal form style \cite{vizza21high,opreni22high}, as will be done for the developments that follow, we have
\begin{subequations} \label{eq:paramOrder1p}
\begin{align}
    \bfW^{(1,d+1)} &= -\bfA_t^{-1} \bar{\bfR}_t, \label{eq:paramOrder1pa} \\
    \bff^{(1,d+1)} &= \bf0. \label{eq:paramOrder1pb}
\end{align}
\end{subequations}

At this point, a few clarifying comments are in order. In the first place, it is interesting to highlight that \cref{eq:Homological1p} corresponds exactly to Eq. (39) of \cite{vizza2024superharm} with $\tilde{\lambda} = 0$ and $\bfC = \bar{\bfR}_t$. This is so because the treatment of the bifurcation parameter in this contribution and of the forcing in \cite{vizza2024superharm} is analogous, such that arriving at the same kinds of equations is natural. In particular, the results in \cref{eq:paramOrder1p} could have been obtained by solving Eq. (42) of the previously mentioned paper. Indeed, the fact that $\tilde{\lambda} = 0$ means that no resonances are present at order one (since it is assumed that all modes have non-vanishing eigenvalues), making the kernel of $\tilde{\lambda}\bfB-\bfA_t$ of dimension 0 and the set $\mathcal{R}$ defined in \cite{vizza2024superharm} empty, enabling \cref{eq:paramOrder1p} to be retrieved. Finally, comparing \cref{eq:paramOrder1pa} with \cref{eq:pEigenvalue1storder}, it can be noted that $\bfW^{(1,d+1)} = \bfY_{d+1}$, as in this scenario $\bfA_0 = \bfRtil_t$.

From this remark, using the splitting between the displacement and velocity parts of the eigenvalues and of the mappings, the following equations hold:
\begin{subequations} \label{eq:eigenvecs_simplif}
\begin{align}
    \bfU^{(1)} &= \bfY^U \\
    \bfV^{(1)} &= \bfD \bfY^U,
\end{align}
\end{subequations}
and thus the linear part of the mappings is determined only from the right eigenvalues associated to displacements. 

\subsection*{Order s homological equation}

The order $s$ homological equation assumes the following modified format, already neglecting the last line, related to the bifurcation parameter, as it will once again result in a tautology:
\begin{equation} \label{eq:homologicalp}
\bfB \P[s]{\nabla_{\bfztil} \bfW(\bfztil) \bff(\bfztil)} = \bfA_t \P[s]{\bfW (\bfztil)} + \P[s]{\bfGtil_t(\bfW(\bfztil),\bfW(\bfztil))} - \bfRtil_u \P[s]{p\bfW(\bfztil)} + \P[s]{\bfHtil(\bfW(\bfztil),\bfW(\bfztil),\bfW(\bfztil))}.
\end{equation}

Then, by introducing expansions
\begin{subequations} \label{eq:expansions}
\begin{align}
    \P[s]{\bfW(\bfztil)} &= \sum_{k=1}^{m_p} \bfW^{(p,k)} \bfztil^{\alphavec(p,k)}, \label{eq:Wexpansion} \\
    \P[s]{\bfGtil_t(\bfW(\bfztil),\bfW(\bfztil))} &= \sum_{k=1}^{m_p} \bfGtil_t^{(p,k)} \bfztil^{\alphavec(p,k)}, \label{eq:Gtilexpansion} \\
    \P[s]{\bfHtil(\bfW(\bfztil),\bfW(\bfztil),\bfW(\bfztil)} &= \sum_{k=1}^{m_p} \bfHtil^{(p,k)} \bfztil^{\alphavec(p,k)}. \label{eq:Htilexpansion}
\end{align}
\end{subequations}
where $m_p$ denotes the number of order $p$ monomials and $\alphavec(p,k)$ is a multi-index indicating the normal variables exponents and the term in \cref{eq:Gtilexpansion} can be computed in exactly the same fashion as the quadratic term in \cite{vizza2024superharm}.

Assuming that no Jordan blocks are present in the linear part of the dynamics, the left-hand side part is given by
\begin{equation} \label{eq:nablaWfexpansion}
\P{\nabla_{\bfz} \bfW(\bfztil) \bff(\bfztil)} = \bfN_1(\bfztil) + \bfN_2(\bfztil) + \bfN_3(\bfztil),
\end{equation}
with expressions for $\bfN_1(\bfztil)$, $\bfN_2(\bfztil)$ and $\bfN_3(\bfztil)$ defined in \cite{vizza2024superharm}. In order to consider a non-diagonal linear dynamics, term $\bfNtil_2(\bfztil)$ given in \cref{eq:N2Jordan} has to be included in the calculations. Additionally, for the cubic term in \cref{eq:Htilexpansion}, its coefficients can be determined by
\begin{align}
    \bfHtil^{(p,k)} &= \sum_{p_1=1}^{p-2} \sum_{p_2=1}^{p-p_1-1} \sum_{k_1,k_2,k_3=1}^{m_{p_1},m_{p_2},m_{p_3}} \bfHtil(\bfW^{(p_1,k_1)},\bfW^{(p_2,k_2)},\bfW^{(p_3,k_3)}) \\
    &p_3: \quad p_3 = p - p_1 - p_2 \nonumber \\
    &k: \quad \alphavec(p,k) = \alphavec(p_1,k_1) + \alphavec(p_2,k_2) + \alphavec(p_3,k_3) \nonumber
\end{align}

Finally, the parcel depending on the parameter in \cref{eq:homologicalp} can be calculated by noticing that all terms $p \bfW(\bfz)$ of order $p$ are obtained by having $\bfW(\bfz)$ of order $p-1$. Thus, its expression is of the form
\begin{align} \label{eq:newterm}
\P{p \bfW(\bfz)} &= \sum_{k_W=1}^{m_{p-1}} \bfW^{(k_W,p-1)} p \bfztil^{\alphavec(k_W,p-1)} = \sum_{k_W=1}^{m_{p-1}} \bfW^{(k_W(k),p-1)} \bfztil^{\alphavec(k,p)}, \\
&k: \quad \alphavec(p,k) = \alphavec(p-1,k_W) + \mathbf{e}_{d+1}. \nonumber
\end{align}
With $k_W(k)$ uniquely defined if $k,k_W$ verify the second line in \cref{eq:newterm} and $\bfW^{(k_W(k),p-1)}$ takes as null otherwise. It should be noted that all the quantities in \cref{eq:newterm} are known when solving for order $p$, and therefore this term contributes to the right-hand side of the homological equation.

With this, the system of equations to be solved in order to determine the vectors of unknowns corresponding to the $k-$th monomial of order $p$ is the one given by Eq. (79) of \cite{vizza2024superharm}:
\begin{equation} \label{eq:systemOrig}
\begin{bmatrix}
    \sigma^{(p,k)} \bfB - \bfA_t & \bfB \bfY_{\mathcal{R}} & \bf0 \\
    \bfX_{\mathcal{R}}^* \bfB & \bf0 & \bf0 \\
    \bf0 & \bf0 & \bfI
\end{bmatrix}
\begin{bmatrix}
    \bfW^{(p,k)} \\
    \bff^{(p,k)}_{\mathcal{R}} \\
    \bff^{(p,k)}_{\cancel{\mathcal{R}}}
\end{bmatrix}
=
\begin{bmatrix}
    \bfR^{(p,k)} \\
    \bf0 \\
    \bf0
\end{bmatrix},
\end{equation}
where $\sigma^{(p,k)} = \alphavec(p,k) \cdot \text{diag}(\bfLambda)$ is defined in order to check resonance conditions, $\mathcal{R}$ is the set of monomials resonant to monomial $(p,k)$, {\em i.e.} such that $\lambda_r \approx \sigma^{(p,k)}$, and $\bff^{(p,k)}_{\mathcal{R}}$ and $\bff^{(p,k)}_{\cancel{\mathcal{R}}}$ are the reduced dynamics coefficients of resonant and non-resonant monomials. Also, in order to employ the above equation in the present context, the expression for $\bfR^{(p,k)}$ needs to include the newly defined terms:
\begin{equation} \label{eq:newRHS}
\bfR^{(p,k)} = \bfGtil_t^{(p,k)} + \bfHtil^{(p,k)} - \bfB \left(\bfN_2^{(p,k)} + \bfNtil_2^{(p,k)} + \bfN_3^{(p,k)} \right) + \bfRtil_u \bfW^{(k_W(k),p-1)}.
\end{equation}

Finally, to take advantage of the trivial relationship between displacements and velocities for the mechanical problem, the right-hand side vector can be divided into two parts according to
\begin{equation}
    \bfR^{(p,k)} =
    \begin{bmatrix}
        \bfM \bfmu^{(p,k)} \\
        \bfnu^{(p,k)}
    \end{bmatrix}.
\end{equation}
Where it is possible to factor out the mass matrix in the first line of the equation due to \cref{eq:newRHS} and the structure of $\bfB$, $\bfGtil_t^{(p,k)}$, $\bfHtil^{(p,k)}$ and $\bfRtil_u$. With this, by substituting $\bfW^{(p,k)} = \begin{bmatrix} \bfU^{(p,k)} & \bfV^{(p,k)} \end{bmatrix}^T$, employing \cref{eq:eigenvecs_simplif} and the definitions of matrices $\bfA_t$ and $\bfB$ the first line of \cref{eq:systemOrig} can be split into
\begin{subequations}
\begin{align}
    \bfV^{(p,k)} = \sigma \bfU^{(p,k)} + &\sum_{r \in \mathcal{R}} f_r^{(p,k)} \bfY_r^U - \bfmu^{(p,k)} \label{eq:Vmappings} \\
    \left( \sigma^2 \bfM + \sigma \bfC + \bfK \right) \bfU^{(p,k)} + &\sum_{r \in \mathcal{R}} f_r^{(p,k)} \left[ \left( \sigma + \lambda_r \right) \bfM + \bfC \right] \bfY_r^U = \bfXi^{(p,k)},
\end{align}
\end{subequations}
with
\begin{equation}
    \bfXi^{(p,k)} = \bfnu^{(p,k)} + (\sigma \bfM + \bfC) \bfmu^{(p,k)},
\end{equation}
and where the superscrpit $(p,k)$ has been dropped in $\sigma$ in order to lighten notation. Analogously, the second line in \cref{eq:systemOrig} is transformed into
\begin{equation}
    \bfX_r^{V^*} \left[ \left( \sigma + \lambda_r \right) \bfM + \bfC \right] \bfU^{(p,k)} + \sum_{s \in \mathcal{R}} f_s^{(p,k)} \bfX_r^{V^*} \bfM \bfY_s^U = \bfX_r^{V^*} \bfM \bfmu^{(p,k)}, \quad \forall r \in \mathcal{R},
\end{equation}
such that the final, halved in size, system of equations to be solved is given by
\begin{equation} \label{eq:systemRed}
\begin{bmatrix}
    \sigma^2 \bfM + \sigma \bfC + \bfK & \left[ \left( \sigma \bfI_{\mathcal{R}} + \bfLambda_{\mathcal{R}} \right) \bfM + \bfC \right] \bfY_{\mathcal{R}}^U & \bf0 \\
    \bfX_{\mathcal{R}}^{V^*} \left[ \left( \sigma \bfI_{\mathcal{R}} + \bfLambda_{\mathcal{R}} \right) \bfM + \bfC \right] & \bfX_{\mathcal{R}}^{V^*} \bfM \bfY_{\mathcal{R}}^U & \bf0 \\
    \bf0 & \bf0 & \bfI
\end{bmatrix}
\begin{bmatrix}
    \bfU^{(p,k)} \\
    \bff^{(p,k)}_{\mathcal{R}} \\
    \bff^{(p,k)}_{\cancel{\mathcal{R}}}
\end{bmatrix}
=
\begin{bmatrix}
    \bfXi^{(p,k)} \\
    \bfX_{\mathcal{R}}^{V^*} \bfM \bfmu^{(p,k)} \\
    \bf0
\end{bmatrix},
\end{equation}
with the velocity mappings $\bfV^{(p,k)}$ found {\em a posteriori} by use of \cref{eq:Vmappings}.

\end{document}

%% file: newcommands.tex
\renewcommand{\P}[2][p]{{\left[#2\right]_{#1}}}
\newcommand{\dblnk}{\begin{minipage}{12pt}\fontsize{6pt}{3pt}\selectfont{$>1$\\$<p$}\end{minipage}}

\newcommand{\norm}[1]{\left\lVert#1\right\rVert}

\newcommand{\dOmega}{{\text d}\Omega}

\newcommand{\bfm}[1]{\mathbf{#1}}
\newcommand{\bsn}[1]{\boldsymbol{#1}}

\newcommand{\bfe}{\bfm{e}}
\newcommand{\bff}{\bfm{f}}

\newcommand{\bfn}{\bfm{n}}
\newcommand{\bfu}{\bfm{u}}

\newcommand{\bfy}{\bfm{y}}
\newcommand{\bfytil}{\mathbf{\widetilde{y}}}

\newcommand{\bft}{\bfm{t}}
\newcommand{\bfv}{\bfm{v}}
\newcommand{\bfx}{\bfm{x}}
\newcommand{\bfz}{\bfm{z}}
\newcommand{\bfztil}{\bfm{\widetilde{z}}}

\newcommand{\alphavec}{\boldsymbol \alpha}

\newcommand{\thetavec}{\boldsymbol \theta}
\newcommand{\bfM}{\bfm{M}}
\newcommand{\bfK}{\bfm{K}}
\newcommand{\bfKg}{\bfm{K_g}}
\newcommand{\bfH}{\bfm{H}}
\newcommand{\bfHtil}{\bfm{\widetilde{H}}}
\newcommand{\bfD}{\bfm{D}}
\newcommand{\bfA}{\bfm{A}}
\newcommand{\bfAtil}{\bfm{\widetilde{A}}}
\newcommand{\bfB}{\bfm{B}}
\newcommand{\bfBtil}{\bfm{\widetilde{B}}}
\newcommand{\bfC}{\bfm{C}}
\newcommand{\bfE}{\bfm{E}}
\newcommand{\bfF}{\bfm{F}}
\newcommand{\bfG}{\bfm{G}}
\newcommand{\bfGtil}{\bfm{\widetilde{G}}}
\newcommand{\bfI}{\bfm{I}}
\newcommand{\bfX}{\bfm{X}}

\newcommand{\bfXbar}{\bfm{\bar{X}}}
\newcommand{\bfY}{\bfm{Y}}

\newcommand{\bfYbar}{\bfm{\bar{Y}}}
\newcommand{\bfN}{\bfm{N}}
\newcommand{\bfNtil}{\bfm{\widetilde{N}}}

\newcommand{\bfQ}{\bfm{Q}}
\newcommand{\bfQtil}{\bfm{\widetilde{Q}}}
\newcommand{\bfR}{\bfm{R}}
\newcommand{\bfRtil}{\bfm{\widetilde{R}}}

\newcommand{\bfS}{\bfm{S}}
\newcommand{\bfV}{\bfm{V}}
\newcommand{\bfU}{\bfm{U}}
\newcommand{\bfW}{\bfm{W}}

\newcommand{\bftheta}{\bsn{\theta}}

\newcommand{\bfeps}{\bsn{\varepsilon}}

\newcommand{\fC}{\mathfrak{C}}

\newcommand{\da}{{\text d}a}

\newcommand{\ds}{{\text d}s}

\newcommand{\bfLambda}{\pmb{\Lambda}}

\newcommand{\bfmu}{\pmb{\mu}}
\newcommand{\bfnu}{\pmb{\nu}}
\newcommand{\bfXi}{\pmb{\Xi}}


%% file: main_firstdraft.bbl
\begin{thebibliography}{10}

\bibitem{PONSIOEN2018}
S.~Ponsioen, T.~Pedergnana, and G.~Haller.
\newblock Automated computation of autonomous spectral submanifolds for
  nonlinear modal analysis.
\newblock {\em Journal of Sound and Vibration}, 420:269 -- 295, 2018.

\bibitem{PONSIOEN2020}
S.~Ponsioen, S.~Jain, and G.~Haller.
\newblock Model reduction to spectral submanifolds and forced-response
  calculation in high-dimensional mechanical systems.
\newblock {\em Journal of Sound and Vibration}, 488:115640, 2020.

\bibitem{BreunungHaller18}
T.~Breunung and G.~Haller.
\newblock Explicit backbone curves from spectral submanifolds of forced-damped
  nonlinear mechanical systems.
\newblock {\em Proceedings of the Royal Society A: Mathematical, Physical and
  Engineering Sciences}, 474(2213):20180083, 2018.

\bibitem{artDNF2020}
A.~Vizzaccaro, Y.~Shen, L.~Salles, J.~Blahos, and C.~Touz{\'e}.
\newblock Direct computation of nonlinear mapping via normal form for
  reduced-order models of finite element nonlinear structures.
\newblock {\em Computer Methods in Applied Mechanics and Engineering},
  284:113957, 2021.

\bibitem{JAIN2021How}
S.~Jain and G.~Haller.
\newblock How to compute invariant manifolds and their reduced dynamics in
  high-dimensional finite-element models.
\newblock {\em Nonlinear Dynamics}, 107:1417--1450, 2022.

\bibitem{vizza21high}
A.~Vizzaccaro, A.~Opreni, L.~Salles, A.~Frangi, and C.~Touz\'e.
\newblock High order direct parametrisation of invariant manifolds for model
  order reduction of finite element structures: application to large amplitude
  vibrations and uncovering of a folding point.
\newblock {\em Nonlinear Dynamics}, 110:525--571, 2022.

\bibitem{opreni22high}
A.~Opreni, A.~Vizzaccaro, C.~Touz\'e, and A.~Frangi.
\newblock High order direct parametrisation of invariant manifolds for model
  order reduction of finite element structures: application to generic forcing
  terms and parametrically excited systems.
\newblock {\em Nonlinear Dynamics}, 111:5401--5447, 2023.

\bibitem{Mereles:rotor}
A.~Mereles, D.~Stuani Alves, and K.~Lucchesi Cavalca.
\newblock Model reduction of rotor-foundation systems using the approximate
  invariant manifold method.
\newblock {\em Nonlinear Dynamics}, 111:10743–10768, 2023.

\bibitem{martin2019}
A.~Martin and F.~Thouverez.
\newblock Dynamic analysis and reduction of a cyclic symmetric system subjected
  to geometric nonlinearities.
\newblock {\em Journal of Engineering for Gas Turbines and Power}, 141:041027,
  2019.

\bibitem{li2021periodic}
M.~Li, S.~Jain, and G.~Haller.
\newblock Nonlinear analysis of forced mechanical systems with internal
  resonance using spectral submanifolds -- part {I}: {P}eriodic response and
  forced response curve.
\newblock {\em Nonlinear Dynamics}, 110:1005--1043, 2022.

\bibitem{opreniPiezo}
A.~Opreni, G.~Gobat, C.~Touz\'e, and A.~Frangi.
\newblock Nonlinear model order reduction of resonant piezoelectric
  micro-actuators: an invariant manifold approach.
\newblock {\em Computers and Structures}, 289:107154, 2023.

\bibitem{Frangi:electromech}
A.~Frangi, A.~Colombo, A.~Vizzaccaro, and C.~Touzé.
\newblock Reduced order modelling of fully coupled electro-mechanical systems
  through invariant manifolds with applications to microstructures.
\newblock {\em International Journal for Numerical Methods in Engineering},
  submitted, 2023.

\bibitem{Pinho:shells}
F.~A.~X.~Carneiro Pinho, M.~Amabili, Z.~J.~G.~N.~Del Prado, and F.~M.~Alves
  da~Silva.
\newblock Nonlinear forced vibration analysis of doubly curved shells via the
  parameterization method for invariant manifold.
\newblock {\em Nonlinear Dynamics}, 112:20677–20701, 2024.

\bibitem{Buza:NS}
G.~Buza.
\newblock Spectral submanifolds of the {N}avier-{S}tokes equation.
\newblock {\em SIAM Journal on Applied Dynamical Systems}, 23(2):1052--1089,
  2024.

\bibitem{LeBihan_2017}
B.~Le Bihan, J.~J. Masdemont, G.~G{\'o}mez, and S.~Lizy-Destrez.
\newblock Invariant manifolds of a non-autonomous quasi-bicircular problem
  computed via the parameterization method.
\newblock {\em Nonlinearity}, 30(8):3040, 2017.

\bibitem{ShawPierre91}
S.~W. Shaw and C.~Pierre.
\newblock Non-linear normal modes and invariant manifolds.
\newblock {\em Journal of Sound and Vibration}, 150(1):170--173, 1991.

\bibitem{Shaw94}
S.~W. Shaw.
\newblock An invariant manifold approach to nonlinear normal modes of
  oscillation.
\newblock {\em Journal of Nonlinear Science}, 4:419--448, 1994.

\bibitem{touze03-NNM}
C.~Touz\'e, O.~Thomas, and A.~Chaigne.
\newblock Hardening/softening behaviour in non-linear oscillations of
  structural systems using non-linear normal modes.
\newblock {\em Journal of Sound and Vibration}, 273(1-2):77--101, 2004.

\bibitem{TOUZE:JSV:2006}
C.~Touz\'e and M.~Amabili.
\newblock Non-linear normal modes for damped geometrically non-linear systems:
  application to reduced-order modeling of harmonically forced structures.
\newblock {\em Journal of Sound and Vibration}, 298(4-5):958--981, 2006.

\bibitem{Haller2016}
G.~Haller and S.~Ponsioen.
\newblock Nonlinear normal modes and spectral submanifolds: existence,
  uniqueness and use in model reduction.
\newblock {\em Nonlinear Dynamics}, 86(3):1493--1534, 2016.

\bibitem{ReviewROMGEOMNL}
C.~Touz{\'e}, A.~Vizzaccaro, and O.~Thomas.
\newblock Model order reduction methods for geometrically nonlinear structures:
  a review of nonlinear techniques.
\newblock {\em Nonlinear Dynamics}, 105:1141--1190, 2021.

\bibitem{TouzeCISM2}
C.~Touz{\'e} and A.~Vizzaccaro.
\newblock Nonlinear normal modes as invariant manifolds for model order
  reduction.
\newblock In C.~Touz{\'e} and A.~Frangi, editors, {\em Model Order Reduction
  for Design, Analysis and Control of Nonlinear Vibratory Systems}, pages
  59--116, New York, NY, 2024. Springer Series CISM courses and lectures, vol.
  614.

\bibitem{LamarqueUP}
C.~H. Lamarque, C.~Touz{\'e}, and O.~Thomas.
\newblock An upper bound for validity limits of asymptotic analytical
  approaches based on normal form theory.
\newblock {\em Nonlinear Dynamics}, 70(3):1931--1949, 2012.

\bibitem{vizza2024superharm}
A.~Vizzaccaro, G.~Gobat, A.~Frangi, and C.~Touzé.
\newblock Direct parametrisation of invariant manifolds for forced
  non-autonomous systems including superharmonic resonances.
\newblock {\em Nonlinear Dynamics}, 112:6255--6290, 2024.

\bibitem{Mereles:bif}
A.~Mereles, D.~Stuani Alves, and K.~Lucchesi Cavalca.
\newblock Bifurcations and limit cycle prediction of rotor systems with
  fluid-film bearings using center manifold reduction.
\newblock {\em Nonlinear Dynamics}, 111:17749–17767, 2023.

\bibitem{NayfehBala}
A.~H. Nayfeh and B.~Balachandran.
\newblock {\em Applied Nonlinear Dynamics: Analytical, Computational, and
  Experimental Methods}.
\newblock Wiley, New-York, 1995.

\bibitem{Sinou2003}
J.-J. Sinou, F.~Thouverez, and L.~Jezequel.
\newblock Analysis of friction and instability by the centre manifold theory
  for a non-linear sprag-slip model.
\newblock {\em Journal of Sound and Vibration}, 265:527--559, 2003.

\bibitem{ShawShock}
S.~W. Shaw, C.~Pierre, and E.~Pesheck.
\newblock Modal analysis-based reduced-order models for nonlinear structures:
  An invariant manifold approach.
\newblock {\em The Shock and Vibration Digest}, 31(1):3--16, 1999.

\bibitem{Hsu:NFa}
L.~Hsu.
\newblock Analysis of critical and post-critical behaviour of non-linear
  dynamical systems by the normal form method, part {I}: normalization
  formulae.
\newblock {\em Journal of Sound and Vibration}, 89(2):169--181, 1983.

\bibitem{Hsu:NFb}
L.~Hsu.
\newblock Analysis of critical and post-critical behaviour of non-linear
  dynamical systems by the normal form method, part {II}: divergence and
  flutter.
\newblock {\em Journal of Sound and Vibration}, 89(2):183--194, 1983.

\bibitem{IoossForcing}
C.~Elphick, G.~Iooss, and E.~Tirapegui.
\newblock Normal form reduction for time-periodically driven differential
  equations.
\newblock {\em Physics Letters A}, 120(9):459--463, 1987.

\bibitem{Jezequel91}
L.~J{\'e}z{\'e}quel and C.~H. Lamarque.
\newblock Analysis of non-linear dynamical systems by the normal form theory.
\newblock {\em Journal of Sound and Vibration}, 149(3):429--459, 1991.

\bibitem{JIANG2005H}
D.~Jiang, C.~Pierre, and S.W. Shaw.
\newblock Nonlinear normal modes for vibratory systems under harmonic
  excitation.
\newblock {\em Journal of Sound and Vibration}, 288(4):791 -- 812, 2005.

\bibitem{Thurnher:resparam}
T.~Thurnher, G.~Haller, and S.~Jain.
\newblock {Nonautonomous spectral submanifolds for model reduction of nonlinear
  mechanical systems under parametric resonance}.
\newblock {\em Chaos: An Interdisciplinary Journal of Nonlinear Science},
  34(7):073127, 2024.

\bibitem{Martin:rotation}
A.~Martin, A.~Opreni, A.~Vizzaccaro, M.~Debeurre, L.~Salles, A.~Frangi,
  O.~Thomas, and C.~Touz{\'e}.
\newblock Reduced order modeling of geometrically nonlinear rotating structures
  using the direct parametrisation of invariant manifolds.
\newblock {\em Journal of Theoretical, Computational and Applied Mechanics},
  10430, 2023.

\bibitem{SilkeGlas:MorMan}
P.~Buchfink, S.~Glas, B.~Haasdonk, and B.~Unger.
\newblock Model reduction on manifolds: A differential geometric framework.
\newblock {\em Physica D: Nonlinear Phenomena}, 468:134299, 2024.

\bibitem{oulghelou:grassmann}
M.~Oulghelou and C.~Allery.
\newblock Non intrusive method for parametric model order reduction using a
  bi-calibrated interpolation on the {G}rassmann manifold.
\newblock {\em Journal of Computational Physics}, 426:109924, 2021.

\bibitem{Barnett2022}
J.~Barnett and C.~Farhat.
\newblock Quadratic approximation manifold for mitigating the kolmogorov
  barrier in nonlinear projection-based model order reduction.
\newblock {\em Journal of Computational Physics}, 464:111348, 2022.

\bibitem{Bolotin1964}
V.~V. Bolotin.
\newblock {\em The dynamic stability of elastic systems}.
\newblock Holden-Day, 1964.

\bibitem{Dowell}
E.~H. Dowell.
\newblock {\em A modern course in aeroelasticity}.
\newblock Kluwer Academic, 2004.

\bibitem{Dimitriadis}
G.~Dimitriadis.
\newblock {\em Introduction to nonlinear aeroelasticity}.
\newblock John Wiley and sons, 2017.

\bibitem{GregoryPaidoussis}
R.W. Gregory and M.P. Pa{\"i}doussis.
\newblock Unstable oscillation of tubular cantilevers conveying fluids, {I.}
  theory.
\newblock {\em Proccedings of the Royal Society London Series A}, page 02936,
  1966.

\bibitem{DOARE2002}
O.~Doaré and E.~{de Langre}.
\newblock The flow-induced instability of long hanging pipes.
\newblock {\em European Journal of Mechanics - A/Solids}, 21(5):857--867, 2002.

\bibitem{Ziegler1952}
H.~Ziegler.
\newblock Die stabilit{\"a}tskriterien der elastomechanik.
\newblock {\em Ingenieur-Archiv}, 20(1):49--56, 1952.

\bibitem{Bentvelsen2018}
B.~Bentvelsen and A.~Lazarus.
\newblock {Modal and stability analysis of structures in periodic elastic
  states: application to the Ziegler column}.
\newblock {\em {Nonlinear Dynamics}}, 91(2):1349--1370, 2018.

\bibitem{Luongo:minZieg}
A.~Luongo and F.~D’Annibale.
\newblock {A paradigmatic minimal system to explain the Ziegler paradox}.
\newblock {\em {Continuum Mechanics and Thermodynamics}}, 27:211--222, 2015.

\bibitem{Beck1952}
M.~Beck.
\newblock Die knicklast des einseitig eingespannten, tangential gedr{\"u}ckten
  stabes.
\newblock {\em Z. Angew. Math. Phys.}, 3(3):225--228, 1952.

\bibitem{Migli:Beck}
G.~Migliaccio and F.~D’Annibale.
\newblock {On the role of different nonlinear damping forms in the dynamic
  behavior of the generalized Beck’s column}.
\newblock {\em {Nonlinear Dynamics}}, 112:13733--13750, 2024.

\bibitem{MingwuLi2024}
M.~Li and L.~Wang.
\newblock Parametric model reduction for a cantilevered pipe conveying fluid
  via parameter-dependent center and unstable manifolds.
\newblock {\em International Journal of Non-Linear Mechanics}, 160:104629,
  2024.

\bibitem{Haro}
A.~Haro, M.~Canadell, J.-L. Figueras, A.~Luque, and J.-M. Mondelo.
\newblock {\em The parameterization method for invariant manifolds. From
  rigorous results to effective computations}.
\newblock Springer, Switzerland, 2016.

\bibitem{COCHELIN2009}
B.~Cochelin and C.~Vergez.
\newblock A high order purely frequency-based harmonic balance formulation for
  continuation of periodic solutions.
\newblock {\em Journal of Sound and Vibration}, 324(1):243 -- 262, 2009.

\bibitem{KARKAR2013}
S.~Karkar, B.~Cochelin, and C.~Vergez.
\newblock A high-order, purely frequency based harmonic balance formulation for
  continuation of periodic solutions: The case of non-polynomial
  nonlinearities.
\newblock {\em Journal of Sound and Vibration}, 332(4):968--977, 2013.

\bibitem{Guillot2019}
L.~Guillot, B.~Cochelin, and C.~Vergez.
\newblock A generic and efficient {T}aylor series-based continuation method
  using a quadratic recast of smooth nonlinear systems.
\newblock {\em International Journal for Numerical Methods in Engineering},
  119(4):261--280, 2019.

\bibitem{Luongo2015}
A.~Luongo and F.~D'Annibale.
\newblock Linear and nonlinear damping effects on the stability of the
  {Z}iegler column.
\newblock In Mohamed Belhaq, editor, {\em Structural Nonlinear Dynamics and
  Diagnosis}, pages 335--352, Cham, 2015. Springer International Publishing.

\bibitem{Rocard49}
Y.~Rocard.
\newblock {\em Dynamique g{\'e}n{\'e}rale des vibrations}.
\newblock Masson, 1949.

\bibitem{fung2008introduction}
Y.C. Fung.
\newblock {\em An Introduction to the Theory of Aeroelasticity}.
\newblock Dover Books on Aeronautical Engineering. Dover Publications, 2008.

\bibitem{Kirillov:stab}
O.~N. Kirillov.
\newblock {\em Nonconservative stability problems of modern physics}.
\newblock De Gruyter Studies in Mathematical Physics 14. De Gruyter,
  Berlin/Boston, 2013.

\bibitem{Heiss2000}
W.~D. Heiss.
\newblock Repulsion of resonance states and exceptional points.
\newblock {\em Phys. Rev. E}, 61:929--932, 2000.

\bibitem{Seyranian2005}
A.~P. Seyranian, O.~N. Kirillov, and A.~A. Mailybaev.
\newblock Coupling of eigenvalues of complex matrices at diabolic and
  exceptional points.
\newblock {\em Journal of Physics A: Mathematical and General},
  38(8):1723--1740, 2005.

\bibitem{Wiggins}
S.~Wiggins.
\newblock {\em Introduction to applied nonlinear dynamical systems and chaos}.
\newblock Springer-Verlag, New-York, 2003.
\newblock Second edition.

\bibitem{Stabile:morfesym}
A.~de~Figueiredo~Stabile, C.~Touzé, and A.~Vizzaccaro.
\newblock Normal form analysis of nonlinear oscillator equations with automated
  arbitrary order expansions.
\newblock {\em Journal of Theoretical, Computational and Applied Mechanics},
  submitted, 2024.

\bibitem{dhooge2004matcont}
A.~Dhooge, W.~Govaerts, and Y.~A. Kuznetsov.
\newblock Matcont: a matlab package for numerical bifurcation analysis of odes.
\newblock {\em ACM SIGSAM Bulletin}, 38(1):21--22, 2004.

\bibitem{gucken83}
J.~Guckenheimer and P.~Holmes.
\newblock {\em Nonlinear oscillations, dynamical systems and bifurcations of
  vector fields}.
\newblock Springer-Verlag, New-York, 1983.

\bibitem{DauchotMann}
O.~Dauchot and P.~Manneville.
\newblock Local versus global concepts in hydrodynamics stability theory.
\newblock {\em Journal de Physique {II} France}, 7(2):371--389, 1997.

\bibitem{DOEDEL2011222}
E.~J. Doedel, B.~Krauskopf, and H.~M. Osinga.
\newblock Global invariant manifolds in the transition to preturbulence in the
  {L}orenz system.
\newblock {\em Indagationes Mathematicae}, 22(3):222--240, 2011.
\newblock Devoted to: Floris Takens (1940–2010).

\bibitem{BowenRuelle}
R.~Bowen and D.~Ruelle.
\newblock The ergodic theory of axiom-{A} flows.
\newblock {\em Inventiones Math.}, 29:181--202, 1975.

\bibitem{EckmanRuelle}
J.P. Eckman and D.~Ruelle.
\newblock Ergodic theory of chaos and strange attractors.
\newblock {\em Review of Modern Physics}, 57:617--656, 1985.

\bibitem{Shilnikov}
D.~V. Turaev and L.~P. Shil'nikov.
\newblock An example of a wild strange attractor.
\newblock {\em Sbornik: Mathematics}, 189(2):291--314, 1998.

\bibitem{Ogden1997}
R.W. Ogden.
\newblock {\em Non-linear Elastic Deformations}.
\newblock Dover Civil and Mechanical Engineering. Dover Publications, 1997.

\bibitem{Belytschko2013}
T.~Belytschko, W.K. Liu, B.~Moran, and K.~Elkhodary.
\newblock {\em Nonlinear Finite Elements for Continua and Structures}.
\newblock Wiley, 2013.

\bibitem{Bonnet2014}
M.~Bonnet, A.~Frangi, and C.~Rey.
\newblock {\em {The finite element method in solid mechanics}}.
\newblock {McGraw Hill Education}, 2014.

\end{thebibliography}
